\documentclass[12pt]{article}
\usepackage{lmodern}
\usepackage{setspace}

\usepackage{graphicx}
\usepackage{times}
\usepackage{epsfig}
\usepackage{amsmath}
\usepackage{amssymb}
\usepackage{algorithm}
\usepackage{algorithmic}
\usepackage{color}

\def\argmax{\mathop{\rm argmax}}
\def\argmin{\mathop{\rm argmin}}

\newcommand{\s}{\ensuremath{\mathbb{S}}}
\newcommand{\real}{\ensuremath{\mathbb{R}}}

\newcommand{\ltwo}{\ensuremath{\mathbb{L}^2}}

\newcommand{\inner}[2]{\left\langle #1,#2 \right\rangle}
\newcommand{\innerd}[2]{\left\langle \left\langle#1,#2 \right\rangle \right\rangle}

\newtheorem{corollary}{Corollary}
\newtheorem{lemma}{Lemma}
\newtheorem{theorem}{Theorem}
\newtheorem{definition}{Definition}

\title{Registration of Functional Data Using Fisher-Rao Metric}
\author{A. Srivastava$^{\dagger}$, W. Wu$^{\dagger}$, S. Kurtek$^{\dagger}$,  E. Klassen$^{\ddagger}$, and J. S. Marron$^*$\\
$^{\dagger}$ Dept. of Statistics \&
$^{\ddagger}$ Dept. of Mathematics, Florida State University,\\
$^*$ Dept.  of Statistics, University of North Carolina
}

\date{}

\begin{document}
\tabcolsep 0pt

\maketitle
\begin{abstract}
We introduce a novel geometric framework for separating the phase
and the amplitude variability in functional data of the type
frequently studied in growth curve analysis. This framework uses the
Fisher-Rao Riemannian metric to derive a proper distance on the
quotient space of functions modulo the time-warping group. A
convenient square-root velocity function (SRVF) representation
transforms the Fisher-Rao metric into the standard $\ltwo$ metric,
simplifying the computations.  This distance is then used to define
a Karcher mean template and warp the individual functions to align
them with the Karcher mean template. The strength of this framework
is demonstrated by deriving a consistent estimator of a signal
observed under random warping, scaling, and vertical translation. These
ideas are demonstrated using both simulated and real data from
different application domains: the Berkeley growth study,
handwritten signature curves, neuroscience spike trains, and gene
expression signals. The proposed method is empirically shown to be
be superior in performance to several recently published methods for
functional alignment.
\end{abstract}

\section{Introduction}
\label{sec:introduction}

The problem of statistical analysis in function spaces is important
in a wide variety of applications arising in nearly every branch of
science, ranging from speech processing to geology, biology and
chemistry. One can easily encounter a problem where the observations
are real-valued functions on an interval, and the goal is
to perform their statistical analysis. By statistical analysis we
mean {\it to compare, align, average,  and model} a collection of
such random observations. These problems can, in principle,  be
addressed using tools from functional analysis, e.g. using the $\ltwo$ Hilbert
structure of the function spaces, where one can compute
$\ltwo$ distances, cross-sectional (i.e. point-wise)  means and variances, and principal
components of the observed functions \cite{ramsay-silverman-2005}.
However, a serious challenge arises when functions are observed with
flexibility or domain warping along the $x$ axis. This
warping may come either from an uncertainty in the measurement
process, or may simply denote an inherent
variability in the underlying process itself that needs to be separated from the
variability along the $y$ axis (or the vertical axis), such as variations in maturity
in the context of growth curves. As another
possibility, the warping may be introduced as a tool to horizontally align the
observed functions,  reduce  their variance and increase
parsimony in the resulting model. We will call these functions {\it
elastic functions}, keeping in mind that we allow only the $x$-axis
(the domain) to be warped and the $y$-values to change only
consequentially.

\begin{figure}
\begin{center}
\includegraphics[height=2.5in]{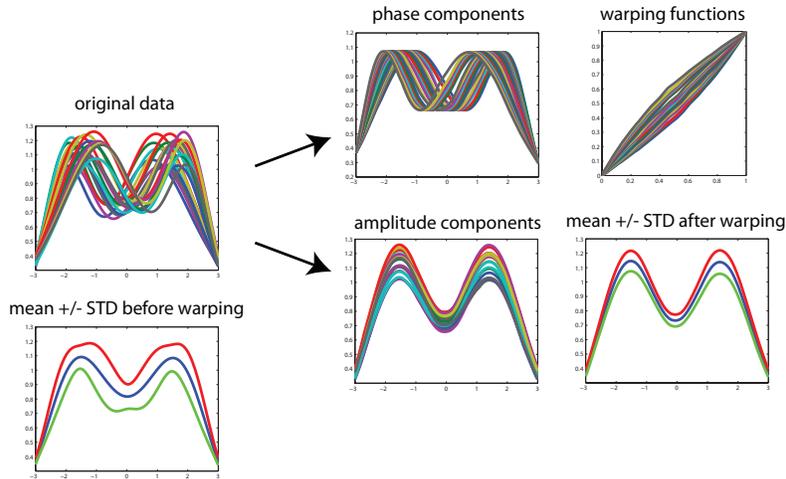}
\caption{Separation of phase and amplitude variability in function data.} \label{fig:motivate}
\end{center}
\end{figure}

Consider the set of functions shown in the top-left panel of Fig.
\ref{fig:motivate}. These functions differ from each other in both
heights and locations of their peaks and valleys. One would like to
separate the variability associated with the heights, called the
{\it amplitude} variability, from the variability associated
with the locations, termed the {\it phase} variability.
Extracting the amplitude variability implies temporally aligning the given
functions using nonlinear time warping, with the result shown in the bottom
right. The corresponding set of warping functions, shown in
the top right, represent the phase variability. The phase
component can also be illustrated by applying these warping
functions to the same function, also shown in the top right.
The main reason for separating functional data into these
components is to better preserve the structure of the observed data, since
a separate modeling of amplitude and phase variability will be more natural,
parsimonious and efficient.

As another, more practical, example we consider the height evolution of subjects in the famous
Berkeley growth data\footnote{http://www.psych.mcgill.ca/faculty/ramsay/datasets.html}.
Fig. \ref{fig:growth-results} shows the time derivatives of the growth
curves, for female and male subjects, to highlight periods of faster growth.
Although the growth rates associated with different individuals are different, it is of
great interest to discover broad common patterns underlying the growth data, particularly
after aligning functions using time warping. Thus,
one would like an automated  technique for alignment of functions.
Section \ref{sec:experiment} shows
examples of data sets from the other applications studied in this paper, including
handwriting curves, gene expression signals, and neuroscience spike trains.

In some applications it may be relatively easy to decide how to warp
functions for a proper alignments. For instance, there may be some temporal
landmarks that have to be aligned across observations. In that case the warping functions
can be piecewise smooth (e.g. linear) functions that ensure that the landmarks
are strictly aligned. This situation requires a manual specification
of landmarks which can be a cumbersome process,
especially for large datasets. In some other cases there may be some natural models that can
be adopted for the warping functions. However, in general, one does not have such landmarks
or natural warping functions, and needs a comprehensive framework where the
alignment of observed functions is performed automatically in an unsupervised fashion.
We seek a principled framework that will automatically
estimate domain warpings of the observed functions in order to optimally align them.
The two main goals of this paper are:

\begin{enumerate}

\item {\bf Joint Alignment and Comparison} (Section 3):
There are two distinct steps in the analysis of elastic functions: (1) warpings or registration of functions and (2) their comparison.
An important requirement in our framework is that these two processes,
warping and comparison, are performed  in a single,
unified framework, i.e. under a single objective function,
as for example was done in \cite{kneip-ramsay:2008}.  A fundamental idea is to avoid treating warping as a {\it pre-processing}
step where the individual functions are warped according to an objective function that is different from
the metric used to compare them.

\item {\bf Signal Estimation Under Random Scales, Translations, and Warpings} (Section 4):
An application of this framework is in estimation of a signal under
the following observation model. Let $f_i$ be an observation of a
function $g$ under random scaling, random vertical translation, and
random warping, and we seek an estimator for $g$ using $\{f_i,
i=1,2,\dots,n\}$. We will use this estimator for performing the
alignment mentioned in the previous item.
\end{enumerate}
Before we introduce our framework that achieves these goals, we present a brief summary of some
past methods, and their strengths and limitations.

\subsection{Past Techniques}
There exists a large literature on statistical analysis of functions, in part due to
the pioneering efforts of Ramsay and Silverman
 \cite{ramsay-silverman-2005},  Kneip and Gasser \cite{kneip-gasser-annals:92},
and several others \cite{muller-JASA:2004,muller-biometrika:2008}.
When restricting to the analysis of elastic functions, the
literature is relatively recent and limited
\cite{ramsay-li-RSSB:98,gervini-gasser-RSSB:04,muller-JASA:2004,muller-biometrika:2008,kneip-ramsay:2008}.
There are basically two categories of the past papers on this
subject. One set treats the problem of functional alignment or
registration as a pre-processing step. Once the functions are
aligned, they are analyzed using the standard tools from functional
analysis, e.g the cross-sectional mean and covariance computation
and PCA. The second set of papers study both comparison and analysis
jointly, using energy-minimization procedures. Although the latter generally
provides better results due to a joint solution, the choice
of the energy function deserves careful scrutiny.

As an example for the first set, in \cite{muller-JASA:2004}, the
authors use warping functions that are convex combinations of
functions of the type: $\gamma_i(t) = \left( { \int_0^t
|f_i^{(\nu)}(s)|^p ds \over  \int_0^1 |f_i^{(\nu)}(s)|^p ds}
\right)^{1/p}$, where $\nu$ and $p$ are two parameters, with the
recommended values being $\nu = 0$ and $p=1$. Then, the warped
functions $\{f_i \circ \gamma_i\}$ are analyzed using standard
statistical techniques under the Hilbert structure of
square-intergable functions. Similarly, James \cite{james:10} uses
moment-based matching for aligning functions, followed up by the
standard FPCA. The main problem with this approach is that the
objective function for alignment is unrelated to the metric for
comparing aligned functions. The two steps are conceptually disjoint
and a change in the objective function for alignment may change the
subsequent results.

We introduce some additional notation.  Let $\Gamma$ be the set of orientation-preserving diffeomorphisms of the
unit interval $[0,1]$:
$\Gamma = \{\gamma: [0,1] \to [0,1] | \gamma(0) = 0,\ \gamma(1)=1,\ \ \gamma\ \ \mbox{is a diffeo}\}$.
Elements of $\Gamma$ form a group, i.e.
 (1) for any $\gamma_1, \gamma_2 \in \Gamma$,  their
composition $\gamma_1 \circ \gamma_2 \in \Gamma$; and (2) for any
$\gamma \in \Gamma$, its inverse $\gamma^{-1} \in \Gamma$, where the
identity is the self-mapping $\gamma_{id}(t) = t$. The role of $\Gamma$ in elastic function
analysis is paramount. Why?
For a function $f \in {\cal F}$, where ${\cal F}$ is an appropriate space of
functions on $[0,1]$ (defined later), the composition $f \circ \gamma$ denotes the
{\it re-parameterization} or a {\it domain warping} of $f$ using $\gamma$. Therefore, $\Gamma$
is also referred to as the re-parameterization or the warping group.
In this paper we will use $\| f \|$ to denote $(\int_0^1 |f(t)|^2 dt)^{1/2}$, i.e., the standard $\ltwo$ norm on the
space of real-valued functions on $[0,1]$.
A majority of past methods study the problem of registration and comparisons of
functions, either separately or jointly, by solving:
 \begin{equation}
 \inf_{\gamma \in \Gamma} \| f_1 - (f_2 \circ \gamma) \|\  \label{eq:quantity1}
 \end{equation}
The use of this quantity is problematic because it is not symmetric. The optimal
alignment of $f_1$ to $f_2$ gives a different minimum, in general, when compared to the
optimal alignment of $f_2$ to $f_1$.
One can enforce a symmetry in Eqn. \ref{eq:quantity1} using a double optimization, i.e. by
 seeking a solution to the problem
$ \inf_{(\gamma_1, \gamma_2) \in \Gamma \times \Gamma} \| (f_1\circ
\gamma_1) - (f_2 \circ \gamma_2) \|$. However, this is a degenerate problem. Another way of ensuring symmetry is to solve:
$ \inf_{\gamma \in \Gamma} \| f_1 - (f_2 \circ \gamma) \| +  \inf_{\gamma \in \Gamma} \| f_2 - (f_1 \circ \gamma) \|$.
While this is symmetric, it still does not lead to  a proper distance on the space of functions.

The basic quantity in Eqn. \ref{eq:quantity1} is commonly
used to form objective functions of the type:
\begin{equation}
E_{\lambda,i}[\nu] =   \inf_{\gamma_i \in \Gamma}  \left(\|  (f_i \circ \gamma_i)  - \nu\|^2  + \lambda\ \ {\cal R}( \gamma_i) \right)\ ,\ \ i=1,2,\dots,n\ , \label{eq:energy-min}
\end{equation}
where ${\cal R}$ is a smoothness penalty on the $\gamma_i$s to keep
them close to $\gamma_{id}(t) = t$. The optimal $\gamma_i^*$ are
then used to align the $f_i$s, followed by a cross-sectional
analysis of the aligned functions. This procedure, once again,
suffers from the problem of separation between the registration and
the comparison steps. Another issue here is: What should $\nu$ be?
It seems natural to use the cross-sectional mean of $f_i$s but that choice is
problematic both empirically and conceptually (more on that later).
Tang and M\"uller
\cite{muller-biometrika:2008} use $\nu = f_j$, obtain a set of
pairwise warping functions $\gamma_{ij}$ for each $i$, and average
them to form the warping function for $f_i$. Kneip and Ramsay
\cite{kneip-ramsay:2008}  take a template-based approach and use a
different $\nu$ for each $i$, given by $\nu_i = \sum_{j=1}^p c^i_j
\xi_j $. Here, the $\xi_j$s are certain basis elements that are also
estimated from the data and, in turn, relate to the principal
components of the observations. Although this formulation has the
nice property of solving for the registration and the principal
components simultaneously,  it implicitly uses the quantity in Eqn.
\ref{eq:quantity1} to compute the residuals.

\subsection{Proposed Approach}
We are going to take a differential geometric approach that provides a natural and fundamental
framework for alignment of elastic functions.
This approach is motivated by recent developments in shape analysis of
parametrized curves \cite{younes-elastic-distance,srivastava_etal_PAMI:10}.
The use of elastic functions for analysis of variance and clustering has also been
studied in \cite{kaziska-FANOVA:10} and for analysis of spike train
data in \cite{wu-srivastava-neuro:10}.

It is problematic to use the cross-sectional mean of $\{f_i\}$
in Eqn. \ref{eq:energy-min} for finding optimal alignments. To understand this issue,
consider the following estimation problem. Let $f_i = c_i (g \circ \gamma_i) + e_i$,
$i=1,2,\dots,n$, represent observations of a signal $g \in {\cal F}$ under random warpings $\gamma_i
\in \Gamma$, scalings $c_i \in \real_+$ and vertical translations $e_i \in \real$, and we seek an estimator for $g$ given $\{ f_i\}$.
Note that estimation of $g$ is equivalent to the alignment of $f_i$s since, given $g$, one
can estimate $\gamma_i$s and compute $f_i \circ \gamma_i^{-1}$ to align them. So, we
focus on deriving an estimator for $g$. In this context, it is easy to see that
the cross-sectional mean for $\{f_i\}$ is not an estimator of $g$. In fact, we claim that
to derive an estimator for $g$ it is more natural to work in the quotient space
${\cal F}/\Gamma$ rather than ${\cal F}$ itself. This quotient space is the set of orbits
of the types $[f] = \{ (f \circ \gamma) | \gamma \in \Gamma\}$. We will show that
the Karcher mean of the orbits $\{ [f_i]\}$ is a consistent estimator of the orbit of
$g$ and that a specific element of that mean orbit, selected using a pre-determined criterion,
is a consistent estimator of $g$.

Now, the definition of Karcher mean requires a proper
distance on ${\cal F}/\Gamma$.  The quantity in Eqn. \ref{eq:quantity1} cannot be used
since
$\| f_1 - f_2\| \neq \| (f_1 \circ \gamma) - (f_2 \circ \gamma)\|$ for general
$f_1, f_2 \in {\cal F}$ and $\gamma \in \Gamma$. (This point was also
noted by Vantini \cite{vantini-functional:2009} although the solution proposed \cite{sangalli-CSDA:2010},
restricting to only the linear warpings, is not for general use.)
Instead, we use $d_{FR}$, the distance  resulting from the Fisher-Rao Riemannian metric, since the action of
$\Gamma$ is by isometries under that metric.
That is,
$d_{FR}(f_1, f_2) = d_{FR}(f_1 \circ \gamma, f_2 \circ \gamma)$, for all $f_1$, $f_2$, and $\gamma$.
Fisher-Rao Riemannian metric was introduced in
1945 by C. R. Rao \cite{rao:45} where he used the Fisher information
matrix  to compare different probability distributions. This metric
was studied rigorously in the 70s and 80s by Amari \cite{amari85}, Efron
\cite{efron-annals:75}, Kass \cite{Kass97}, Cencov \cite{Cencov82},
and others. While those earlier efforts were focused on analyzing
parametric families,  we use the {\it nonparametric} version of the
Fisher-Rao Riemannian metric  in this paper. (This nonparametric form has found an
important use in shape analysis of curves \cite{srivastava_etal_PAMI:10}.)
An important attribute of this
metric is that it is preserved under warping, and
Cencov \cite{Cencov82} showed that it is the only metric with this attribute.
It is difficult to compute the distance $d_{FR}$ directly under this metric but
Bhattacharya
\cite{bhattacharya-43} introduced a square-root representation that greatly
simplifies this calculation. We will modify this square-root
representations for use with more general functions.

\section{Function Representation and Metric}
\label{sec:shaperep} In order to develop a natural and efficient framework for aligning
elastic functions,
we introduce a square-root representation of functions.

\subsection{Representation Space of Functions}
Let $f$ be a real-valued function on the interval $[0,1]$. We are going to restrict to
those $f$ that are absolutely continuous on $[0,1]$; let ${\cal F}$ denote the
set of all such functions. We define a mapping: $Q:
\real \to \real$ according to:
$Q(x) \equiv  \left\{ \begin{array}{cc}
x/\sqrt{|x|} & \mbox{if}\  |x| \neq 0 \\
0 &\ \ \mbox{otherwise}
\end{array} \right.$.
Note that $Q$ is a continuous map. For the purpose of
studying the function $f$, we will represent it using a square-root velocity function (SRVF)
defined as $q: [0,1] \to \real$, where
$q(t) \equiv Q(\dot{f}(t)) = {\dot{f}(t)}/{\sqrt{|\dot{f}(t) |}}$.
This representation includes those functions whose parameterization can become
singular in the analysis.
It can be shown that if the function $f$ is absolutely continuous, then the
resulting SRVF is square integrable. Thus, we will define $\ltwo([0,1],\real)$ (or simply $\ltwo$) to be the set of all SRVFs.
For every $q \in \ltwo$ there
exists a function $f$ (unique up to a constant, or a vertical translation) such that the given $q$
is the SRVF of that $f$. In fact, this function can be obtained precisely
using the equation: $f(t)  = f(0) + \int_0^t q(s) |q(s)| ds$. Thus, the representation $f \Leftrightarrow (f(0), q)$ is
invertible.

If we warp a function $f$ by $\gamma$, how does its SRVF change? The
SRVF of $f \circ \gamma$ is given by:
$\tilde{q}(t) =  {{d \over dt}(f\circ \gamma)(t) \over \sqrt{|{d \over dt}(f \circ \gamma)(t) |} } = (q \circ \gamma)(t) \sqrt{\dot{\gamma}(t)}$.
We will denote this transformation by $(q, \gamma) = (q
\circ \gamma) \sqrt{\dot{\gamma}}$. The motivations for using SRVF
for functional analysis are many and to understand these merits we
first present the relevant metric.

\subsection{Elastic Riemannian Metric}
In this paper we will use the Fisher-Rao Riemannian metric for
analyzing functions.   We remind the reader that a Riemmanian
metric is a smoothly-varying inner product defined on
the tangent spaces of the manifold.
\begin{definition} \label{defn:FR}
For any $f \in {\cal F}$ and $v_1, v_2 \in T_f({\cal F})$, where $T_f({\cal F})$ is
the tangent space to ${\cal F}$ at $f$,
the Fisher-Rao Riemannian metric is defined as the inner product:
\begin{equation}
\innerd{v_1}{v_2}_{f} = {1 \over 4} \int_0^1 \dot{v}_1(t) \dot{v_2}(t) {1 \over |\dot{f}(t)|} dt\ . \label{eq:def-FR}
\end{equation}
\end{definition}
In case we are dealing only with functions such that $\dot{f}(t) \geq 0$, e.g. cumulative
distribution functions or growth curves, then we obtain a more classical
version of the Fisher-Rao metric. Thus, the above definition is a more general
form of the Fisher-Rao metric, the one that deals with signed functions instead of just
density functions.

This metric has many fundamental advantages, including the fact that
it is the only Riemannian metric that is invariant to the domain warping \cite{Cencov82}, and has played an important role in
information geometry.
This metric is somewhat complicated since it
changes from point to point on ${\cal F}$, and it is not
straightforward to derive equations for computing geodesics  in ${\cal F}$.
For instance, the geodesic distance between any two points $f_1, f_2 \in {\cal F}$
is based on finding a geodesic path between them under the F-R metric. This minimization is non-trivial and only some numerical
algorithms are known to attempt this problem.
Once we find a geodesic path connecting $f_1$ and $f_2$ in ${\cal F}$, its
length becomes the geodesic distance $d_{FR}$.
However, a
small transformation provide an enormous simplification of this task. This motivates the use of SRVFs for representing
and aligning elastic functions.
\begin{lemma} \label{lemma:transform}
Under the SRVF representation, the Fisher-Rao Riemannian metric becomes the
standard $\ltwo$ metric.
\end{lemma}
Proof is given in the appendix.
This result can be used to compute the distance $d_{FR}$ between any
two functions as follows. Simply compute the $\ltwo$ distance between the
corresponding SRVFs and set $d_{FR}$ to that value:
$d_{FR}(f_1, f_2)  = \| q_1 - q_2\|$.
The next question is: What is the effect of warping on $d_{FR}$? This
is answered by the following result.
\begin{lemma}
For any two SRVFs $q_1, q_2 \in \ltwo$ and $\gamma \in \Gamma$,
$\| (q_1, \gamma) - (q_2, \gamma)\| = \|q_1 - q_2\|$. \label{lemma:isometry}
\end{lemma}
See the appendix for the proof.
In the case of functions with the non-negativity constraint (that is, $\dot{f} \geq 0$), this transformation was used by
Bhattacharya \cite{bhattacharya-43}.

\subsection{Elastic Distance on Quotient Space}
So far we have defined the Fisher-Rao distance on ${\cal F}$ and
have found a simple way to compute it using SRVFs. But we have not
involved any warping function in the distance calculation and thus
it represents a non-elastic comparison of functions. The next step
is to define an elastic distance between functions as follows. The
orbit of an SRVF $q \in \ltwo$ is given by: $[q] =
\mbox{closure}\{(q, \gamma) | \gamma \in \Gamma \}\ =
\mbox{closure}\{(q \circ \gamma)\sqrt{\dot{\gamma}}) | \gamma \in
\Gamma \}$. It is the set of SRVFs associated with all the
warpings of a function, and their limit points. Any two elements of $[q]$
represent functions which have the same $y$ variability but
different $x$ variability. Let ${\cal S}$ denote the set of all such
orbits. To compare any two orbits we need a metric on ${\cal S}$. 
We will use the Fisher-Rao distance to induce a distance between orbits, and 
we can do that only because under this the action of $\Gamma$  is 
by isometries. 
\begin{definition} \label{def:elast-dist}
For any two functions $f_1,\ f_2 \in {\cal F}$ and the corresponding SRVFs, $q_1, q_2 \in \ltwo$, we
define the elastic distance $d$ on the quotient space ${\cal S}$ to be:
$d([q_1],[q_2]) =
 \inf_{\gamma \in \Gamma}  \|q_1 -  (q_2, \gamma) \|$.
  \end{definition}
Note that the distance $d$ between a function and its  domain-warped
version is zero. However, it can be shown that if two SRVFs belong
to different orbits, then the distance between them is non-zero.
Thus, this distance $d$ is a proper distance (i.e. it satisfies
non-negativity, symmetry, and the triangle inequality) on ${\cal S}$ but
not on $\ltwo$ itself, where it is only a pseudo-distance.

Table 1 provides a quick summary of relationships between the Fisher-Rao metric
and ${\cal F}$ on one hand, and the $\ltwo$ metric and the space of SRVFs on the other.
\begin{table} \label{tab:1}
\begin{center}
\title{Table 1. Bijective Relationship Between Function Space ${\cal F}$ and SRVF space $\ltwo$ }\\
\begin{tabular}{|c||c|c|}
\hline
Item & Function Space ${\cal F}$ & SRVF Space $\ltwo$\\
\hline
Representation & $f$& $(q, f(0))$\\
\hline
Transformation & $f(t) = f(0) + \int_0^t q(s) |q(s)| ds$ & $q(t) =  {\dot{f}(t)}/{\sqrt{|\dot{f}(t) |}}$\\
\hline
Metric & Fisher-Rao Metric & $\ltwo$ Metric \\
& $\innerd{v_1}{v_2}_{\cal F} = \int_0^1\dot{v}_1(t) \dot{v}_2(t) {1 \over |\dot{f}(t)|} dt$ & $\inner{w_1}{w_2} =  \int_0^1 w_1(t) w_2(t) dt $ \\
\hline
Distance & $d_{FR}(f_1, f_2)$ & $ \|q_1 - q_2\|$ \\
\hline
Isometry & $d_{FR}(f_1, f_2) = d_{FR}(f_1 \circ \gamma, f_2 \circ \gamma)$ & $\| q_1 - q_2\| = \| (q_1, \gamma) - (q_2, \gamma)\|$ \\
\hline
Geodesic & Numerical Solution & Straight Line  \\
\hline
Elastic Distance & $d = \inf_{\gamma \in \Gamma} d_{FR}(f_1, f_2 \circ \gamma)$ &
$d =  \inf_{\gamma \in \Gamma} \left(  \|q_1 -  (q_2 \circ \gamma)\sqrt{\dot{\gamma}}) \| \right)$ in ${\cal S}$\\
between $f_1$ and $f_2$ & & Solved Using Dynamic Programming \\
\hline
\end{tabular}
\end{center}
\end{table}

\section{Karcher Mean and Function Alignment} \label{sec:KM}
An important goal of this warping framework  is to
align the functions so as to improve the matching of features
(peaks and valleys) across functions. A natural idea is to compute a cross-sectional mean
of the given functions and then align the given functions to this mean template.
The problem is that we do not have a proper distance function
on $\ltwo$, invariant to time warpings,
that can be used to define a mean.
But we have a distance function on the quotient space ${\cal S}$, so we will use
a mean on that space to derive a template for function alignment. We will do so in
two steps:
First, for a given collection of functions $f_1, f_2, \dots, f_n$, and their SRVFs $q_1, q_2, \dots,
q_n$, we compute the mean of the corresponding orbits $[q_1], [q_2], \dots, [q_n]$
in the quotient space ${\cal S}$; we will call it $[\mu]_n$.
Next, we compute an appropriate element of this mean orbit to define a template $\mu_n$ in $\ltwo$.
Then, the alignment of individual
functions comes from warping their SRVFs to match the template ${\mu}_n$ under the elastic distance.

We remind the reader that if $dist$ denotes the geodesic distance between
points on a Riemannian manifold  $M$, and $\{p_i, i=1,2,\dots,n\}$ is a collection of
points on $M$, then a local minimizer of the cost function $p \mapsto
\sum_{i=1}^n dist(p,p_i)^2$ is defined as the Karcher mean of those points \cite{karcher:77}. It is also 
known by other names such as the intrinsic mean or the Fr\'{e}chet mean.
The algorithm for computing a Karcher mean is based on gradients and has become a standard
procedure in statistics on nonlinear manifolds (see, for example \cite{le-mean-shape}). We will not
present the general procedure but will describe its use in our problem.

\subsection{Karcher Mean of Points in $\Gamma$}

In this section we will define a Karcher mean of a set of warping functions
$\{ \gamma_i\}$, under the Fisher-Rao metric,  using the differential geometry of $\Gamma$.
Analysis on $\Gamma$ is not straightforward because
 it is a nonlinear manifold.
To understand its geometry, we will
represent an element $\gamma \in \Gamma$ by the square-root of its derivative $\psi = \sqrt{\dot{\gamma}}$.
Note that this is the same as the SRVF defined earlier for elements of ${\cal F}$, except that
$\dot{\gamma} > 0$ here.
The identity element $\gamma_{id}$ maps to a constant function with value $\psi_{id}(t) = 1$.
Since $\gamma(0) = 0$, the mapping
from $\gamma$ to $\psi$ is a bijection and one can reconstruct $\gamma$ from $\psi$ using
$\gamma(t) = \int_0^t \psi(s)^2 ds$.
An important advantage of this transformation is that since
$\| \psi\|^2 = \int_0^1 \psi(t)^2 dt = \int_0^1 \dot{\gamma}(t) dt = \gamma(1) - \gamma(0) = 1$,
the set of all such $\psi$s is $\s_{\infty}$, the unit sphere in the Hilbert space $\ltwo$. In other words, the
square-root representation simplifies the complicated geometry of $\Gamma$ to the unit sphere.
Recall that the distance between any two points on the unit sphere, under the Euclidean metric, is
simply the length of the shortest arc of a great circle connecting them on the sphere.
Using Lemma \ref {lemma:transform},  the Fisher-Rao distance between
any two warping functions is found to be
$d_{FR}(\gamma_1, \gamma_2)  =  \cos^{-1}(\int_0^1 \sqrt{\dot{\gamma}_1(t)} \sqrt{\dot{\gamma}_2(t)} dt )$.
Now that we have a proper distance on $\Gamma$, we can define a Karcher mean as follows.
\begin{definition} \label{def:mean-Gamma}
For a given set of warping functions $\gamma_1, \gamma_2, \dots,
\gamma_n \in \Gamma$, define their Karcher mean to be
$\bar{\gamma}_n = \argmin_{\gamma \in \Gamma} \sum_{i=1}^n d_{FR}(\gamma, \gamma_i)^2$.
\end{definition}
The search for this minimum is performed
using Algorithm 1 as follows: \\
\title{\bf Algorithm 1: Karcher Mean of $\{\gamma_i\}$ Under $d_{FR}$}: \\
Let $\psi_i = \sqrt{\dot{\gamma}_i}$ be the SRVFs for the given warping functions. Initialize $\mu_{\psi}$ to be one of the $\psi_i$s or
use $w/\|w\|$, where $w = {1 \over n}\sum_{i=1}^n \psi_i$.
\begin{enumerate}
\item For $i=1,2,\dots,n$, compute the shooting vector
$v_i = { \theta_i \over \sin(\theta_i)}(\psi_i - \cos(\theta_i)
\mu_{\psi})$, where $\theta_i =
\cos^{-1}(\int_0^1\mu_{\psi}(t)\psi_i(t)dt)$.

\item Compute the average $\bar{v} = {1 \over n} \sum_{i=1}^n v_i$.
\item If $\|\bar{v}\|$ is small, then stop.
Else, update $\mu_{\psi} \mapsto \cos(\epsilon \| \bar{v}\|)\mu_{\psi} + \sin(\epsilon \|\bar{v}\|){ \bar{v} \over \|\bar{v}\|}$,
for a small step size $\epsilon> 0$ and return to Step 1.

\item Compute the mean warping function using
$\bar{\gamma}_n = \int_0^t \mu_{\psi}(s)^2 ds$.
\end{enumerate}

\subsection{Karcher Mean of Points in ${\cal S} = \ltwo/\Gamma$ }

Next we consider the problem of finding means of points in the quotient space ${\cal S}$. Since we
already have a well-defined distance on ${\cal S}$ (given in Definition \ref{def:elast-dist}), the
definition of the Karcher mean follows.
\begin{definition} \label{def:mean1}
Define the Karcher mean $[\mu]_n$ of the given SRVF orbits $\{[q_i]\}$ in the space ${\cal S}$
as a local minimum of the
sum of squares of elastic distances:
\begin{equation}
[\mu]_n = \argmin_{[q] \in {\cal S}} \sum_{i=1}^n d([q],[q_i])^2\ . \label{eq:Karcher-min}
\end{equation}
\end{definition}
We emphasize that the Karcher mean $[\mu]_n$ is actually an orbit
of functions, rather than a function.
That is, if $\mu_0$ is a minimizer of the cost function in Eqn. \ref{eq:Karcher-min}, then
so is $(\mu_0, \gamma)$ for any $\gamma$.
The full algorithm for computing the Karcher
mean in ${\cal S}$ is given next.\\

\noindent
\title{\bf Algorithm 2: Karcher Mean of $\{[q_i]\}$ in ${\cal S}$} \label{algo:mean-svrf}
\begin{enumerate}
\item Initialization Step:
Select  $\mu = q_j$, where $j$ is any index in $\argmin_{1\le i \le
n} || q_i - \frac{1}{n} \sum_{k=1}^n q_k||$.

\item For each $q_i$ find $\gamma_i^*$ by solving:
$\gamma_i^* = \argmin_{\gamma \in \Gamma} \| \mu - (q_i \circ
\gamma) \sqrt{\dot{\gamma}}\|$. The solution to this optimization
comes from a dynamic programming algorithm. In cases where a solution 
does not exist in $\Gamma$, the dynamic programming algorithm still provides an 
approximation in $\Gamma$.

\item Compute the aligned SRVFs using $\tilde{q}_i  \mapsto (q_i \circ \gamma_i^*) \sqrt{\dot{\gamma_i^*}}$.
\item  If the increment $\| {1 \over n} \sum_{i=1}^n \tilde{q}_i - \mu\|$ is small, then stop.
Else, update the mean using $\mu \mapsto {1 \over n} \sum_{i=1}^n
\tilde{q}_i$ and return to step 2.

\end{enumerate}
The iterative update in Steps 2-4 is based on the gradient of the
cost function given in Eqn. \ref{eq:Karcher-min}. Although we prove
its convergence next, its convergence to a global minimum is not
guaranteed. Denote the estimated mean in the $k$th iteration by
$\mu^{(k)}$.  In the $k$th iteration, let $\gamma_i^{(k)}$ denote
the optimal domain warping from $q_i$ to $\mu^{(k)}$ and let $\tilde
q_i^{(k)} = (q_i \circ \gamma_i^{(k)})\sqrt{\dot \gamma_i^{(k)}}$.
Then, $\sum_{i=1}^n d([\mu^{(k)}], [q_i])^2 = \sum_{i=1}^n \|
\mu^{(k)}-\tilde q_i^{(k)} \|^2 \ge \sum_{i=1}^n \|
\mu^{(k+1)}-\tilde q_i^{(k)} \|^2
         \ge \sum_{i=1}^n d([\mu^{(k+1)}], [q_i])^2$.
Thus, the cost function decreases iteratively and as zero is a
natural lower bound, $\sum_{i=1}^n d([\mu^{(k)}], [q_i])^2$ will
always converge.

\subsection{Center of an Orbit}
The remaining task is to find
a particular element of this mean orbit so that it can be used as a template to align
the given functions. Towards this purpose, we will define the center of an orbit using a condition
similar to past papers,
see e.g. \cite{muller-biometrika:2008}, which says that the mean of the warping
functions should be the identity. A major difference here is that
we use the Karcher mean and not the cross-sectional mean as was done in the past.

\begin{definition}\label{def:mean2}
For a given set of SRVFs $q_1, q_2, \dots, q_n$ and $q$, define an element $\tilde{q}$ of $[q]$
as the center of $[q]$ with respect to the set $\{q_i\}$ if the warping functions $\{\gamma_i\}$,
where $\gamma_i = \argmin_{\gamma \in \Gamma} \| \tilde{q} - (q_i, \gamma)\|$,
have the Karcher mean $\gamma_{id}$.
\end{definition}
We will prove the existence of such an element by construction.

\begin{figure}
\begin{center}
\includegraphics[height=2.0in]{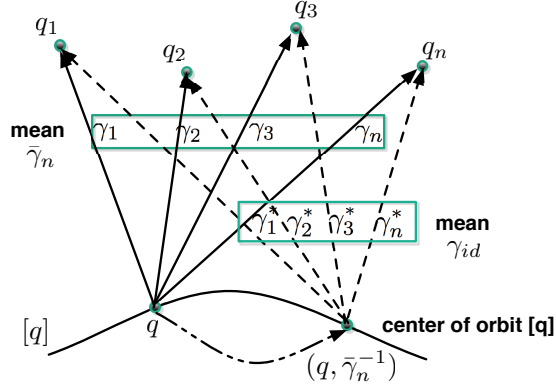}
\caption{Finding center of the orbit $[q]$ with respect to the set $\{q_i\}$.}
\label{fig:template-select}
\end{center}
\end{figure}

\noindent
\title{\bf Algorithm 3: Finding Center of an Orbit} \label{algo:find-mean}:
WLOG, let $q$ be any element of the orbit $[q]$.
\begin{enumerate}
\item For each $q_i$ find $\gamma_i$ by solving:
$\gamma_i = \argmin_{\gamma \in \Gamma}  \left( \| q - (
q_i\circ \gamma) \sqrt{\dot{\gamma}}\| \right)$.
\item Compute the mean $\bar{\gamma}_n$ of all $\{\gamma_i\}$ using Algorithm 1.
The center of $[q]$ {\it wrt} $\{q_i\}$ is given by $\tilde{q} =  (q,  \bar{\gamma}_n^{-1})$.
\end{enumerate}
This algorithm is depicted pictorially in Fig. \ref{fig:template-select}
We need to show that $\tilde{q}$ resulting from Algorithm 3
satisfies the mean condition in Definition \ref{def:mean2}.  Note
that $\gamma_i$ is chosen to minimize $\| q - (q_i, \gamma) \|$, and
also that $\| \tilde{q} - (q_i, \gamma) \|  =  \| (q,
\bar{\gamma}_n^{-1}) - (q_i, \gamma) \|  = \| q - (q_i, \gamma \circ
\bar{\gamma}_n) \|$. Therefore, $\gamma_i^* = \gamma_i \circ
\bar{\gamma}_n^{-1}$ minimizes $\| \tilde{q} - (q_i, \gamma) \|$.
That is, $\gamma_i^*$ is a warping that aligns $q_i$ to $\tilde{q}$.
To verify the Karcher mean of $\gamma_i^*$, we compute the sum of
squared distances $ \sum_{i=1}^n d_{FR}( \gamma, \gamma_i^*)^2 =
\sum_{i=1}^n d_{FR}(\gamma,\gamma_i \circ \bar{\gamma}_n^{-1})^2 =
\sum_{i=1}^n d_{FR}(\gamma \circ \bar{\gamma}_n, \gamma_i)^2$.  As
$\bar{\gamma}_n$ is already the mean of $\gamma_i$, this sum of
squares is minimized when $\gamma = \gamma_{id}$.  That is, the mean
of $\gamma_i^*$ is $\gamma_{id}$.

We will apply this setup in our problem by finding the center of $[\mu]_n$
with respect to the given SRVFs $\{q_i\}$.

\subsection{Complete Alignment Algorithm}
Now we can utilize the three algorithms, Algorithm 1-3, to present the full
procedure for finding a template $\mu_n$ that is used to align the individual functions.

\noindent {\bf Complete Alignment Algorithm}: Given a set of functions $f_1, f_2, \dots f_n$ on $[0,1]$,
let $q_1, q_2, \dots, q_n$ denote their SRVFs, respectively.
\begin{enumerate}
\item Computer the Karcher mean of $[q_1],\ [q_2],\dots, [q_n]$ in ${\cal S}$ using Algorithm 2.
Denote it by $[\mu]_n$.

\item Find the center of $[\mu]_n$ {\it wrt} $\{q_i\}$ using Algorithm 3; call it $\mu_n$. (Note that this
algorithm requires a step for computing the Karcher mean of warping functions using Algorithm 1).

\item For $i=1,2,\dots,n$, find $\gamma_i^*$  by solving:
$\gamma_i^* = \argmin_{\gamma \in \Gamma} \| \mu_n - (q_i, \gamma)\|$.

\item Compute the aligned SRVFs $\tilde{q}_i = (q_i, \gamma_i^*)$
and aligned functions $\tilde{f}_i = f_i \circ \gamma_i^*.$

\item Return the template $\mu_n$, the warping functions $\{ \gamma_i^*\}$, and
the aligned functions $\{ \tilde{f}_i\}$.

\end{enumerate}

\subsection{Simulation Results}

To illustrate this method we use a number of simulated datasets. 
Although our framework is developed for functions on $[0,1]$, it
can easily be adapted to an arbitrary interval using a linear transformation.

\begin{enumerate}

\item {\bf Simulated Data 1}: As the first example, we study a
set of simulated functions used previously in
\cite{kneip-ramsay:2008}. The individual functions are given by:
$y_i(t) = z_{i,1} e^{-(t-1.5)^2/2} + z_{i,2}e^{-(t+1.5)^2/2}$,
$i=1,2,\dots, 21$, where $z_{i,1}$ and $z_{i,2}$ are {\it i.i.d}
normal with mean one and standard deviation $0.25$. Each of these
functions is then warped according to: $\gamma_i(t) =
6({e^{a_i(t+3)/6} -1 \over e^{a_i} - 1}) - 3$ if  $a_i \neq 0$,
otherwise $\gamma_i = \gamma_{id}$, where $a_i$ are equally spaced
between $-1$ and $1$,  and the observed functions are computed using
$f_i(t) = y_i(\gamma_i(t))$. A set of 21 such functions forms the
original data and is shown in the left panel of  Fig.
\ref{fig:mean-result-sim}, and the remaining panels show the results
of our method. The second panel presents the resulting aligned
functions $\{ \tilde{f}_i\}$ and the third panel plots the
corresponding warping functions $\{ \gamma_i^*\}$. The remaining panels show the
cross-sectional mean and mean $\pm$ standard deviations of $\{f_i\}$
and $\{\tilde{f}_i\}$, respectively.
\begin{figure}
\begin{center}
\begin{tabular}{ccccc}
$\{f_i\}$ & $\{ \tilde{f}_i\}$ & $\{ \gamma_i^*\}$ &
mean $\pm$ std, before  & mean $\pm$ std, after  \\
\includegraphics[height=1.0in]{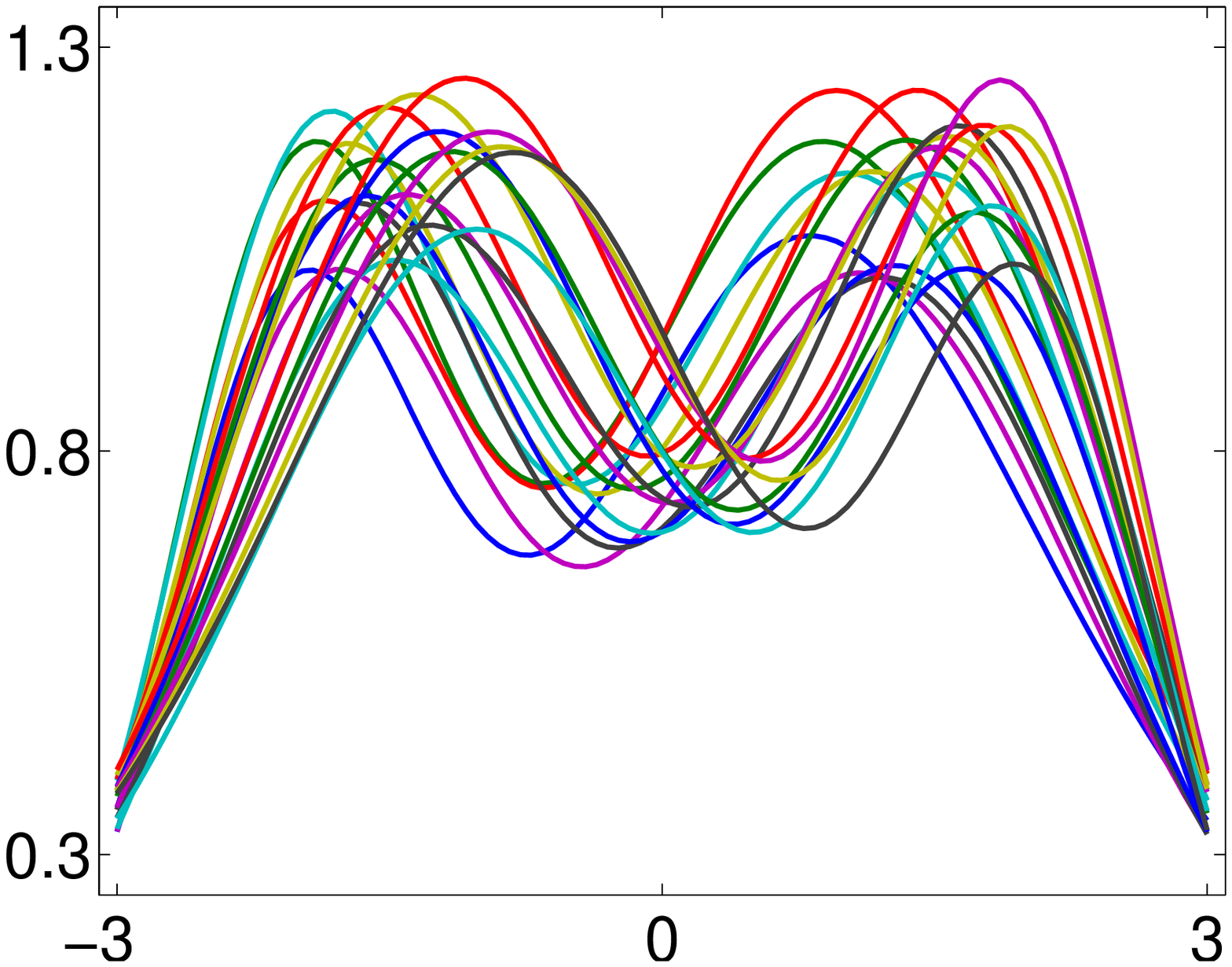} \ &
\includegraphics[height=1.0in]{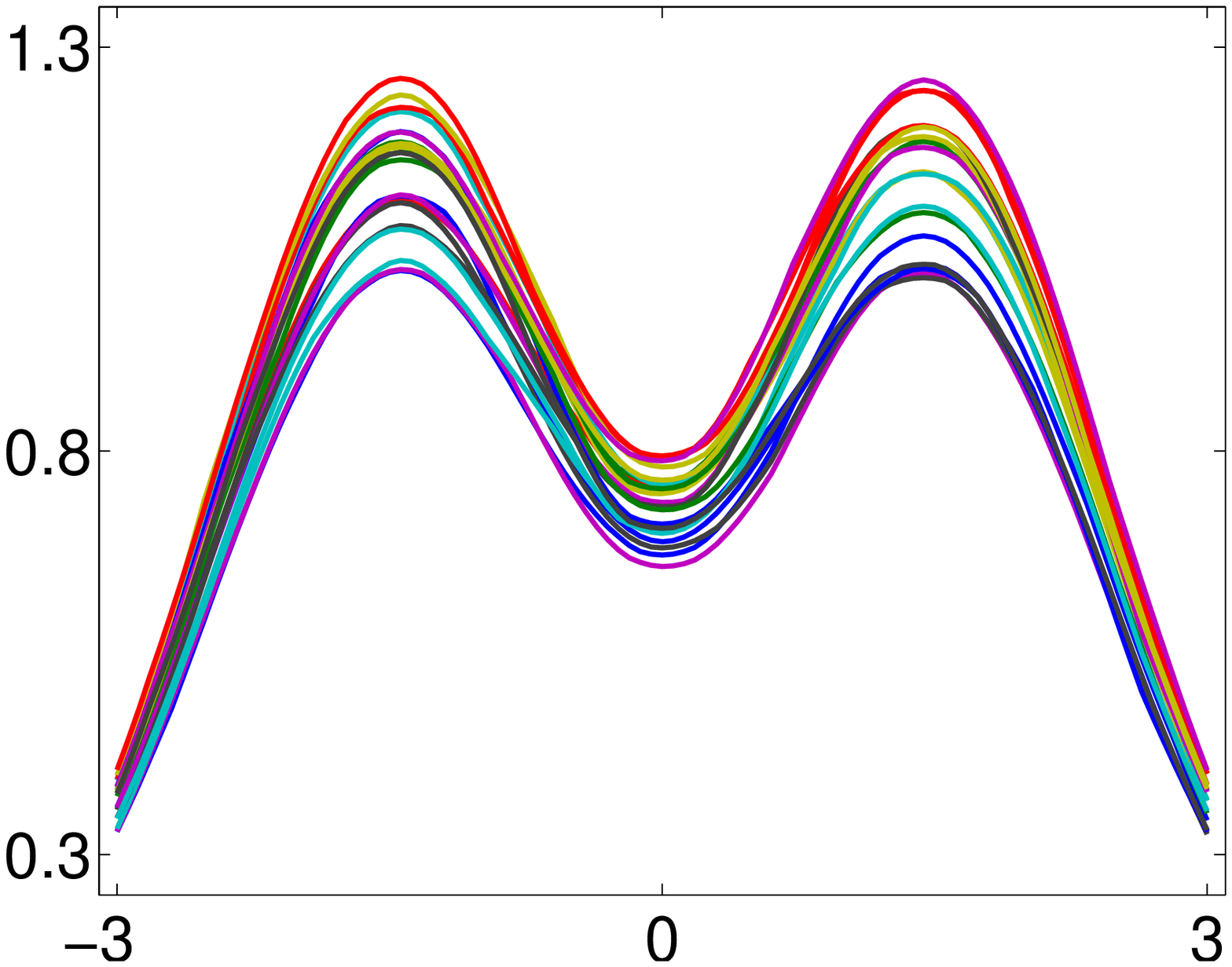} \ &
\includegraphics[height=1.0in]{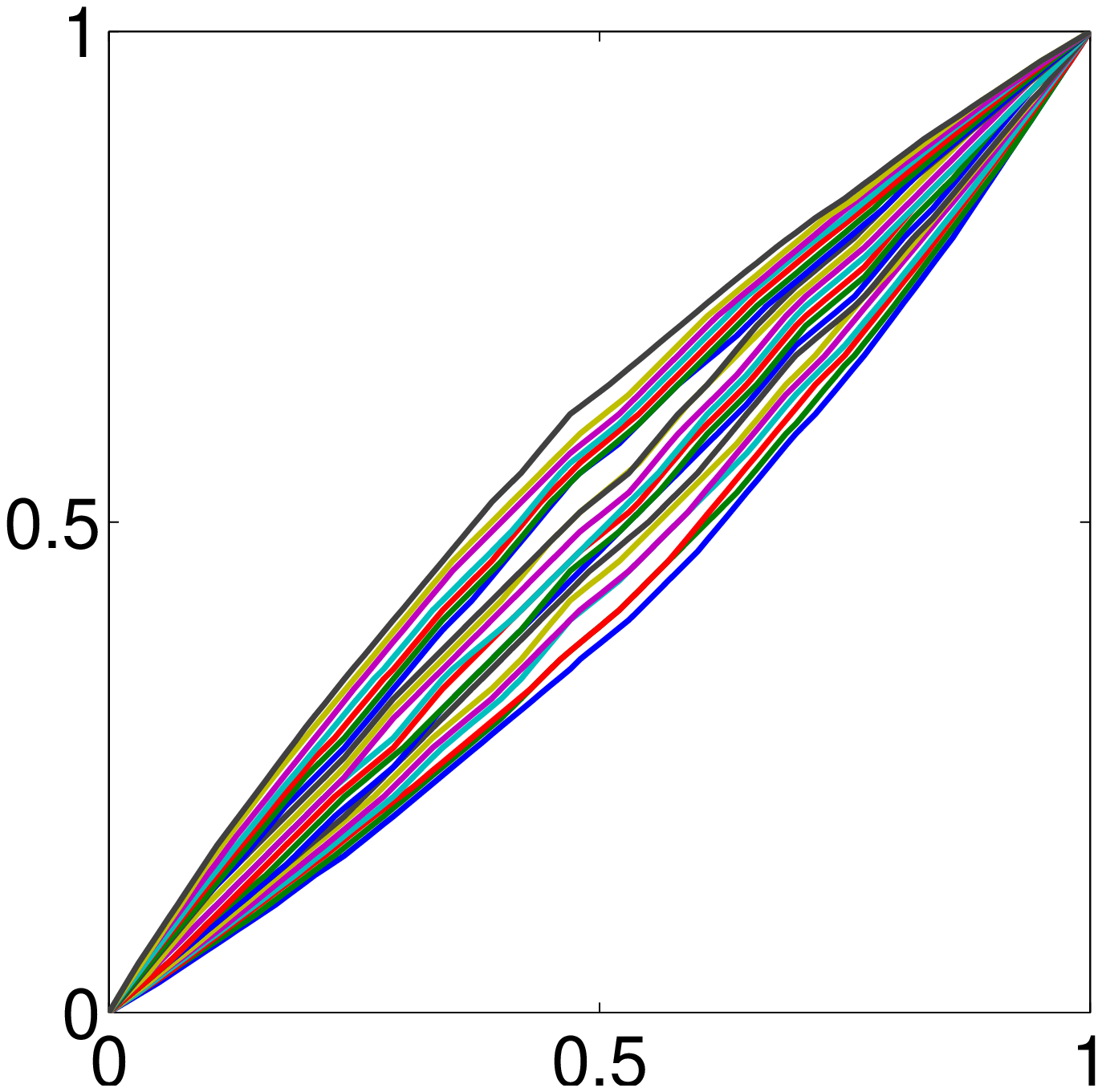} \ &
\includegraphics[height=1.0in]{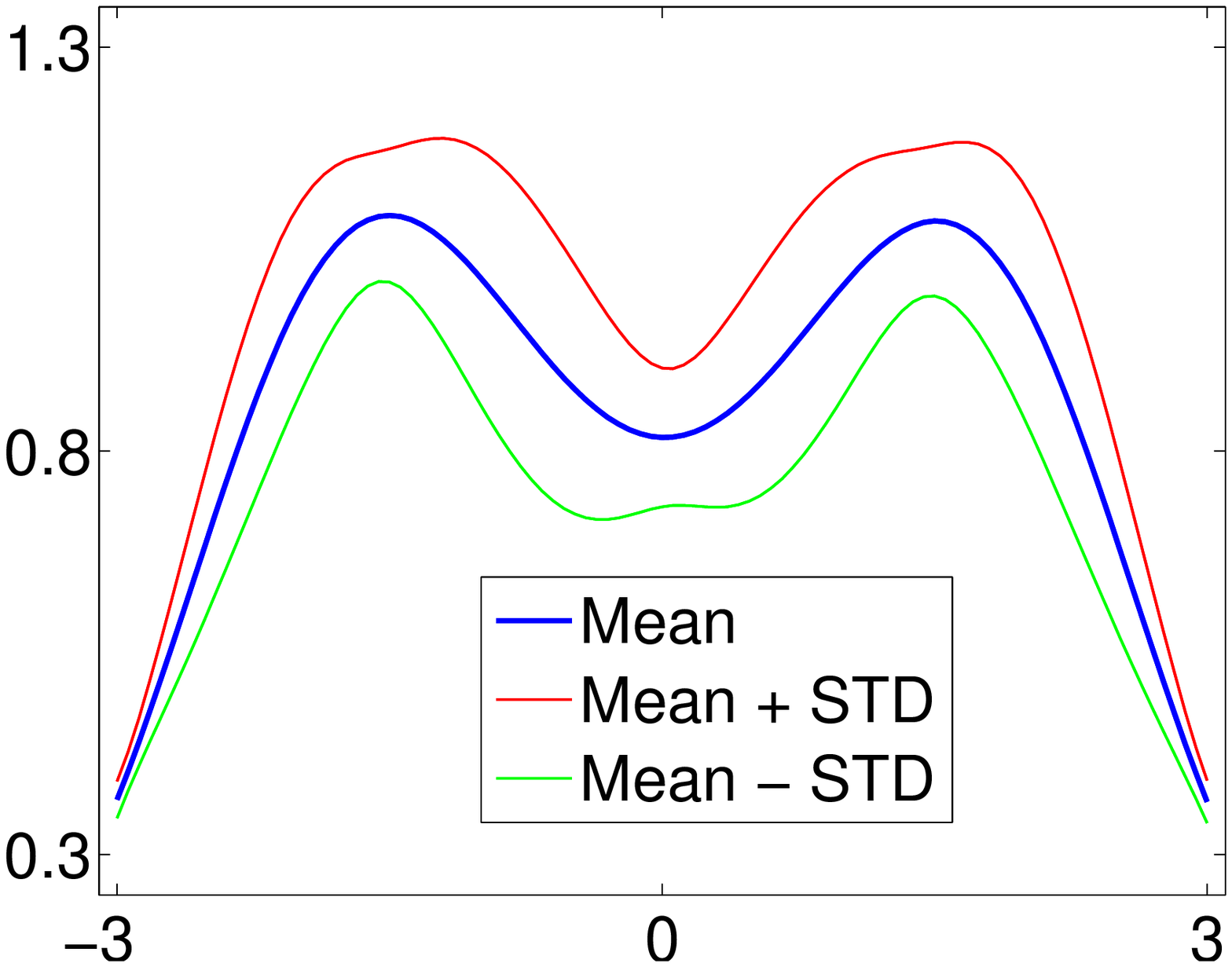} \ &
\includegraphics[height=1.0in]{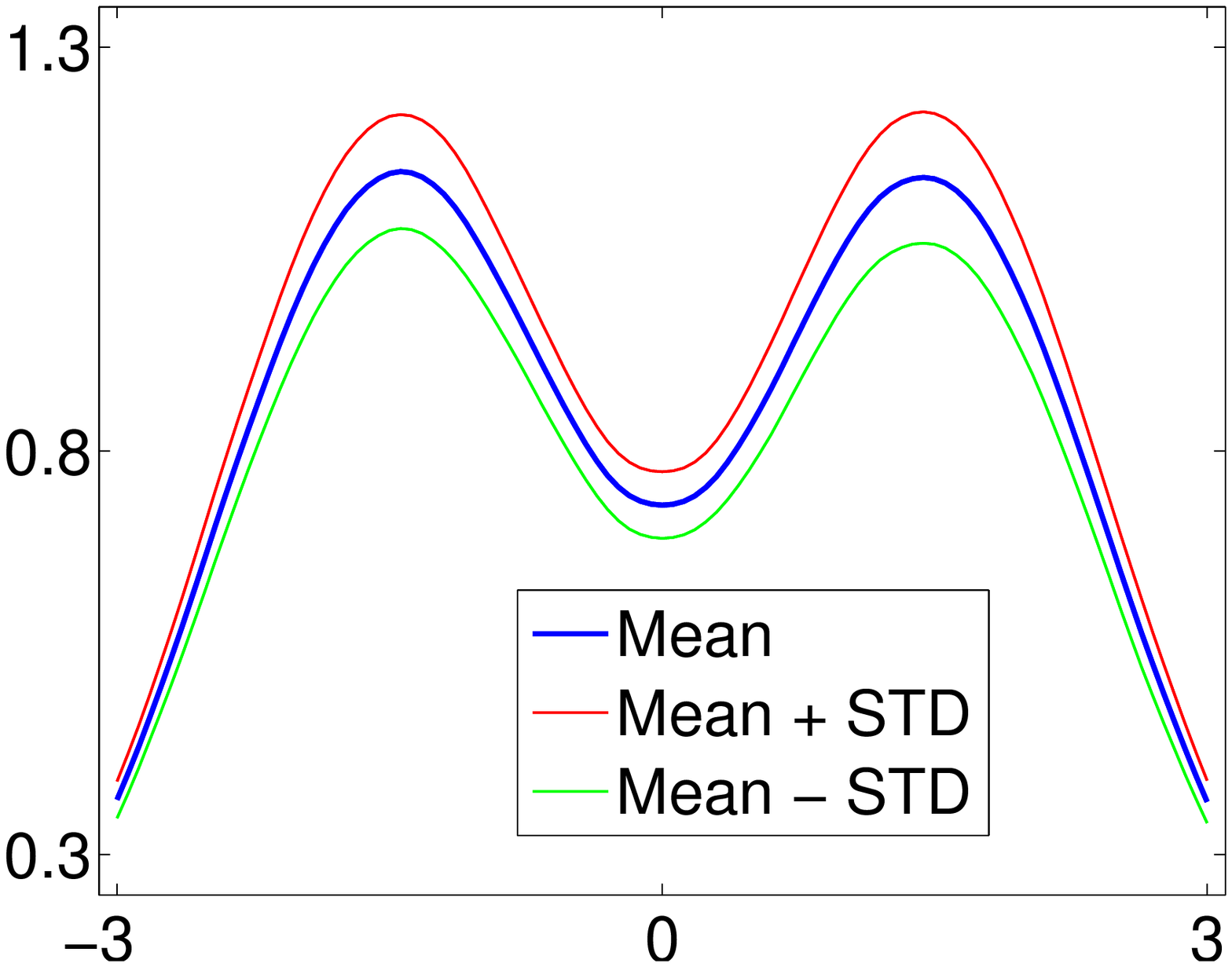}
\end{tabular}
\caption{Results on simulated data set 1.}
\label{fig:mean-result-sim}
\end{center}
\end{figure}
The plot of $\{\tilde{f}_i\}$ shows a tighter alignment of functions
with sharper peaks and valleys. The two peaks are at $-1.5$ and
$1.5$ which is exactly what we expect. This means that the effects
of warping generated by the $\gamma_i$s have been completely removed
and only the randomness from the $y_i$s remains. Also, the plot of
mean $\pm$ standard deviation shows a thinning of bands around the
mean due to the alignment.

\item {\bf Simulated Data 2}: As a simple test of our method we analyze a set of functions with no
underlying phase variability. To do this, we take $\{y_i\}$, as above, but this time we do not warp them
at all; these functions are shown in the left panel of Figure
\ref{fig:result-sim2}. Note that, by construction, the two peaks in
these functions are always aligned, only their amplitudes are
different. There is a slight misalignment in the valleys between the two peaks due to
differing mixture weights.
The result of the alignment process
is shown in the remaining panels. The second panel shows that the aligned
functions are very similar to the original data, except for
a better alignment of the valleys. The next panel shows
the estimated warping functions which are very close to the
identity. The last panel shows the means of the original and the
aligned functions and
 they are practically identical.
\begin{figure}
\begin{center}
\begin{tabular}{cccc}
$\{f_i\}$ & $\{ \tilde{f}_i\}$ & $\{ \gamma_i^*\}$ & mean, before
and after\\
\includegraphics[height=1.0in]{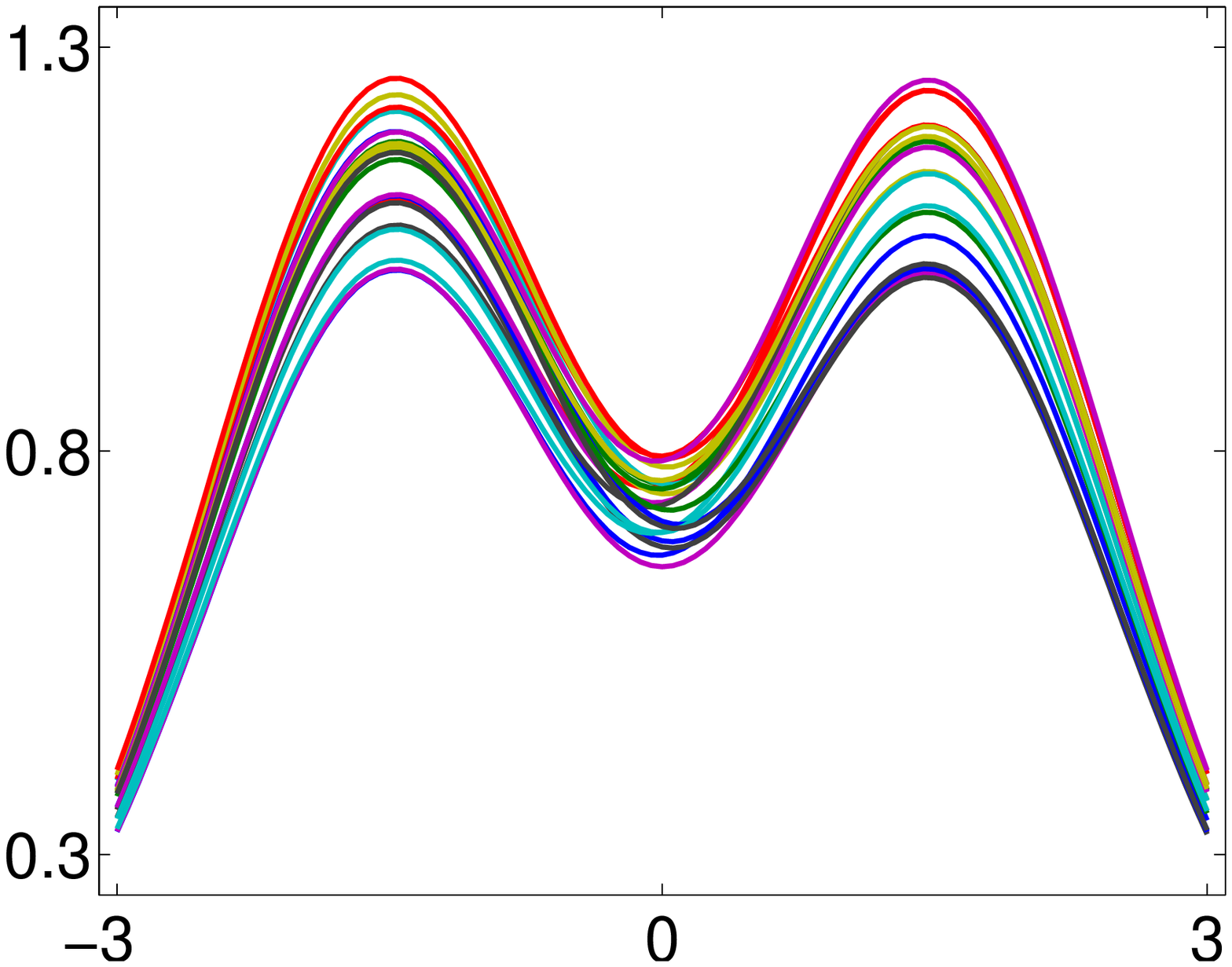} \ &
\includegraphics[height=1.0in]{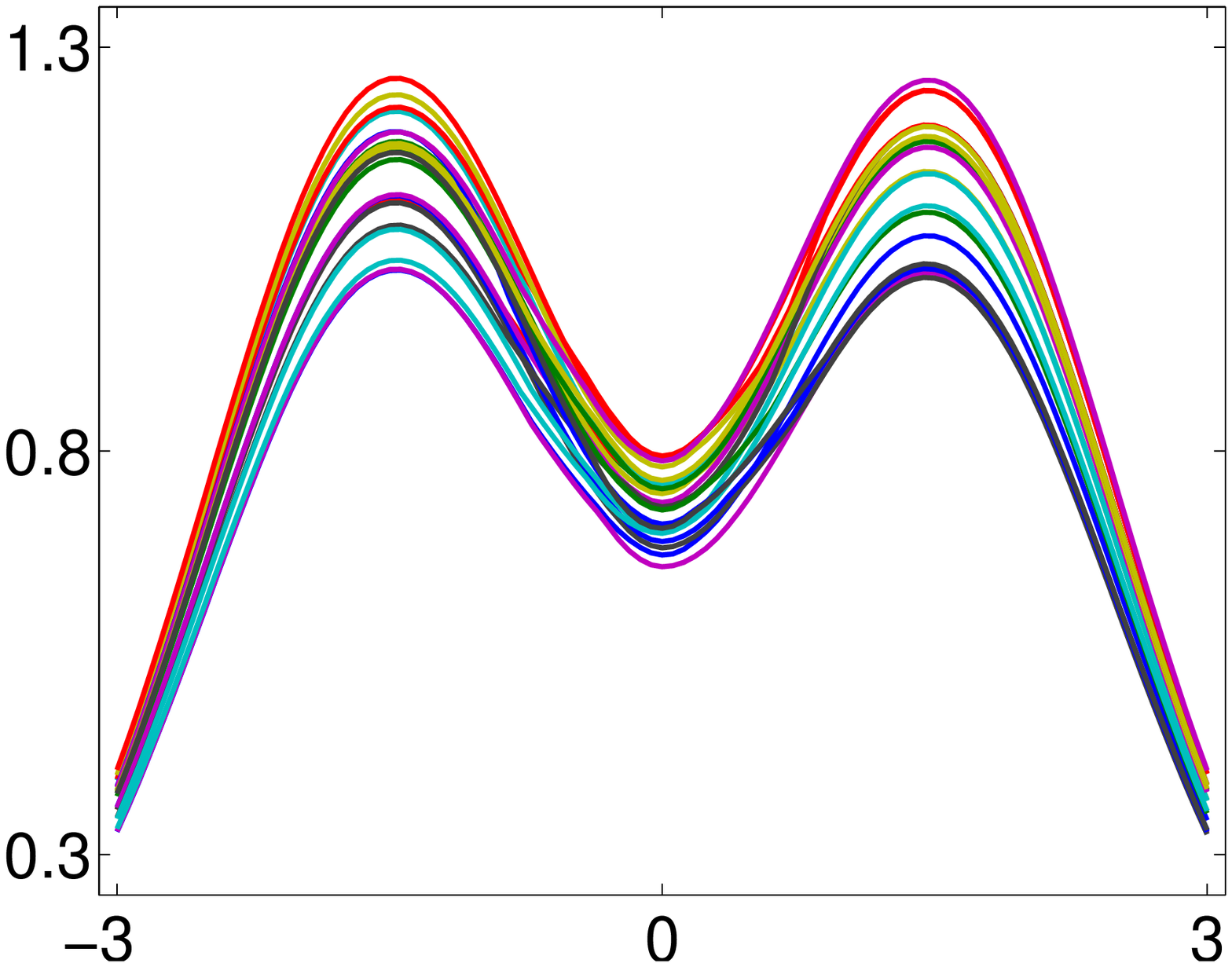} \ &
\includegraphics[height=1.0in]{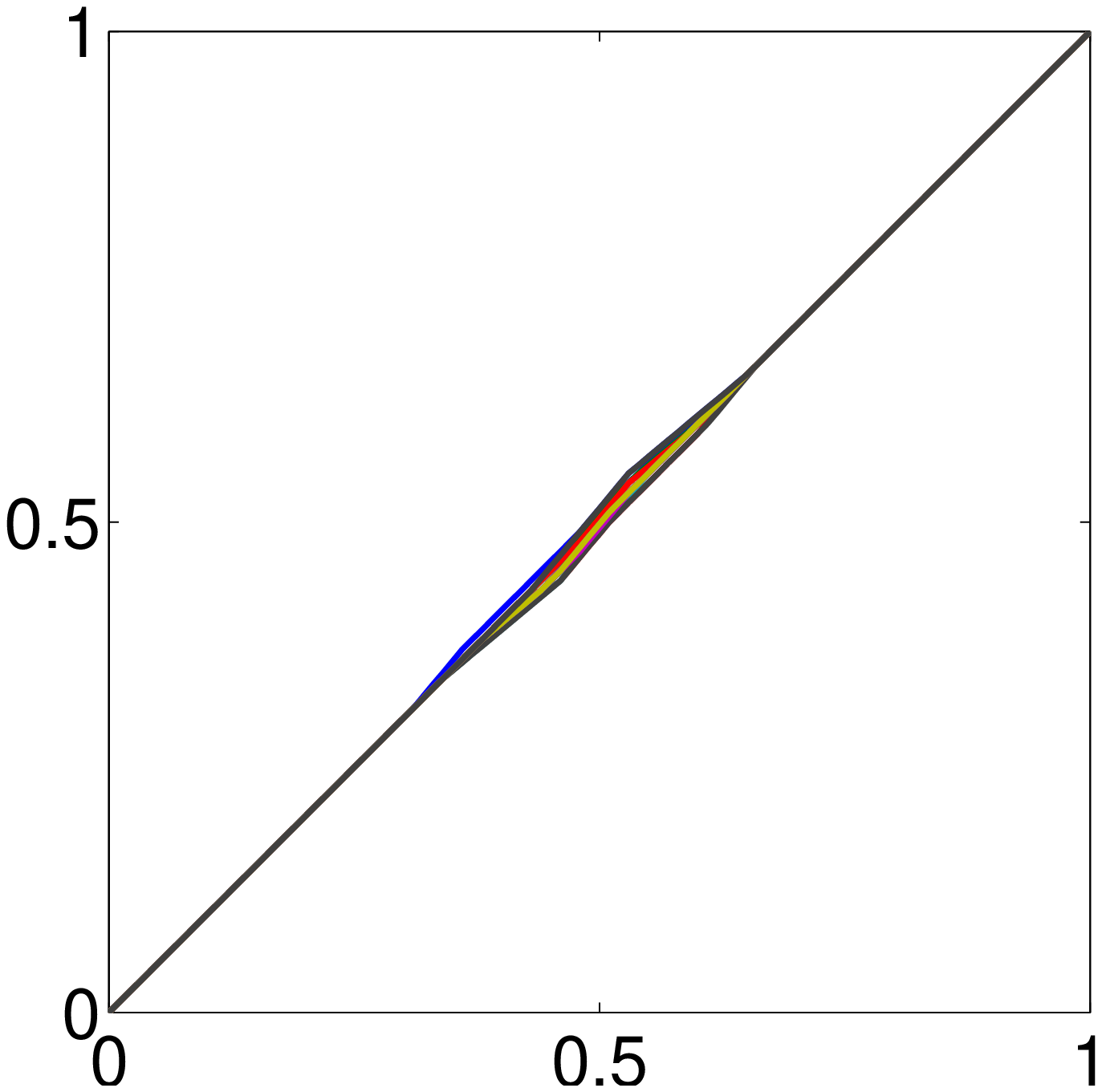} \ &
\includegraphics[height=1.0in]{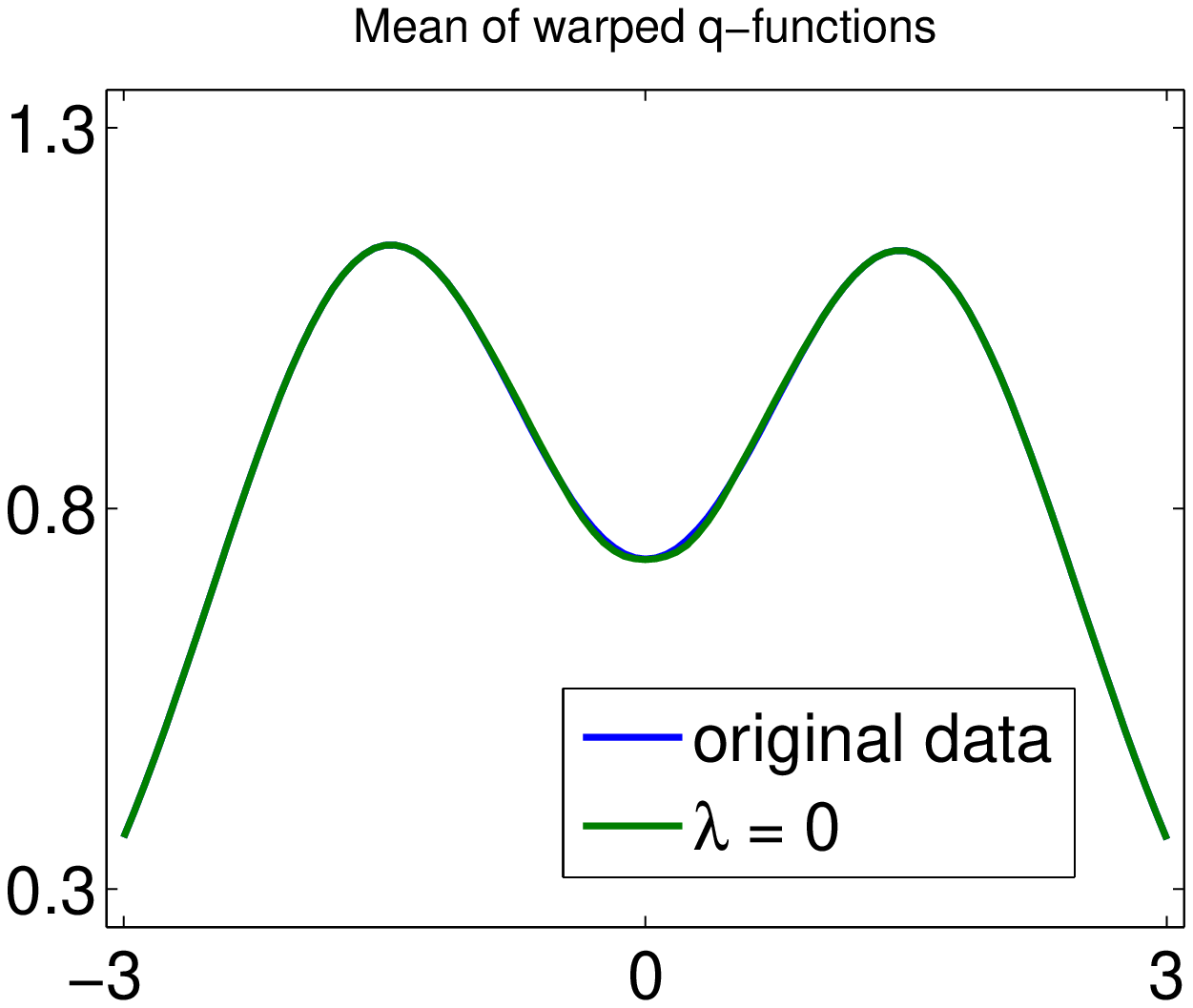}
\end{tabular}
\caption{Results on simulated data 2.} \label{fig:result-sim2}
\end{center}
\end{figure}

\item {\bf Simulated Data 3}: In this case we take a family of Gaussian kernel
functions with the same shape but with significant phase
variability, in the form of horizontal shifts, and minor
amplitude variation. Figure \ref{fig:result-sim3} shows the original
$29$ functions $\{ f_i \}$, the aligned functions $\{
\tilde{f}_i\}$, the warping functions $\{\gamma_i^*\}$, and the
before-and-after cross sectional mean and standard deviations. Once
again we notice a tighter alignment of functions with only minor
variability left in $\{ \tilde{f}_i\}$ reflecting the differing
heights in the original data. The remaining two plots show that mean 
$\pm$ standard deviation of the aligned data is far more compact than 
the raw data.

\begin{figure}
\begin{center}
\begin{tabular}{ccccc}
\hspace*{-0.3in} $\{f_i\}$ & $\{ \tilde{f}_i\}$ & $\{ \gamma_i^*\}$
& mean $\pm$ std, before  & mean $\pm$ std, after
\\
\includegraphics[height=1.0in]{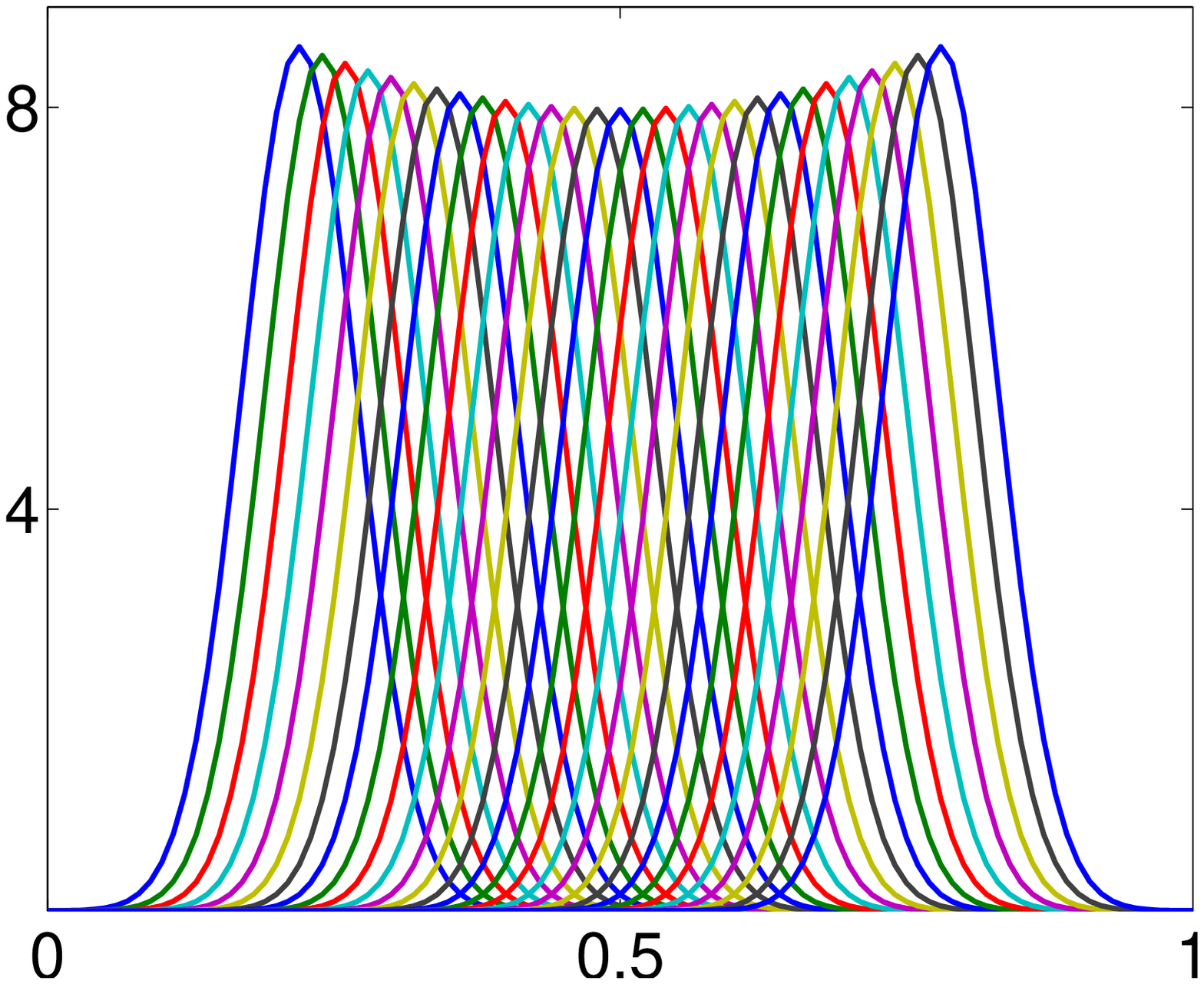} \ &
\includegraphics[height=1.0in]{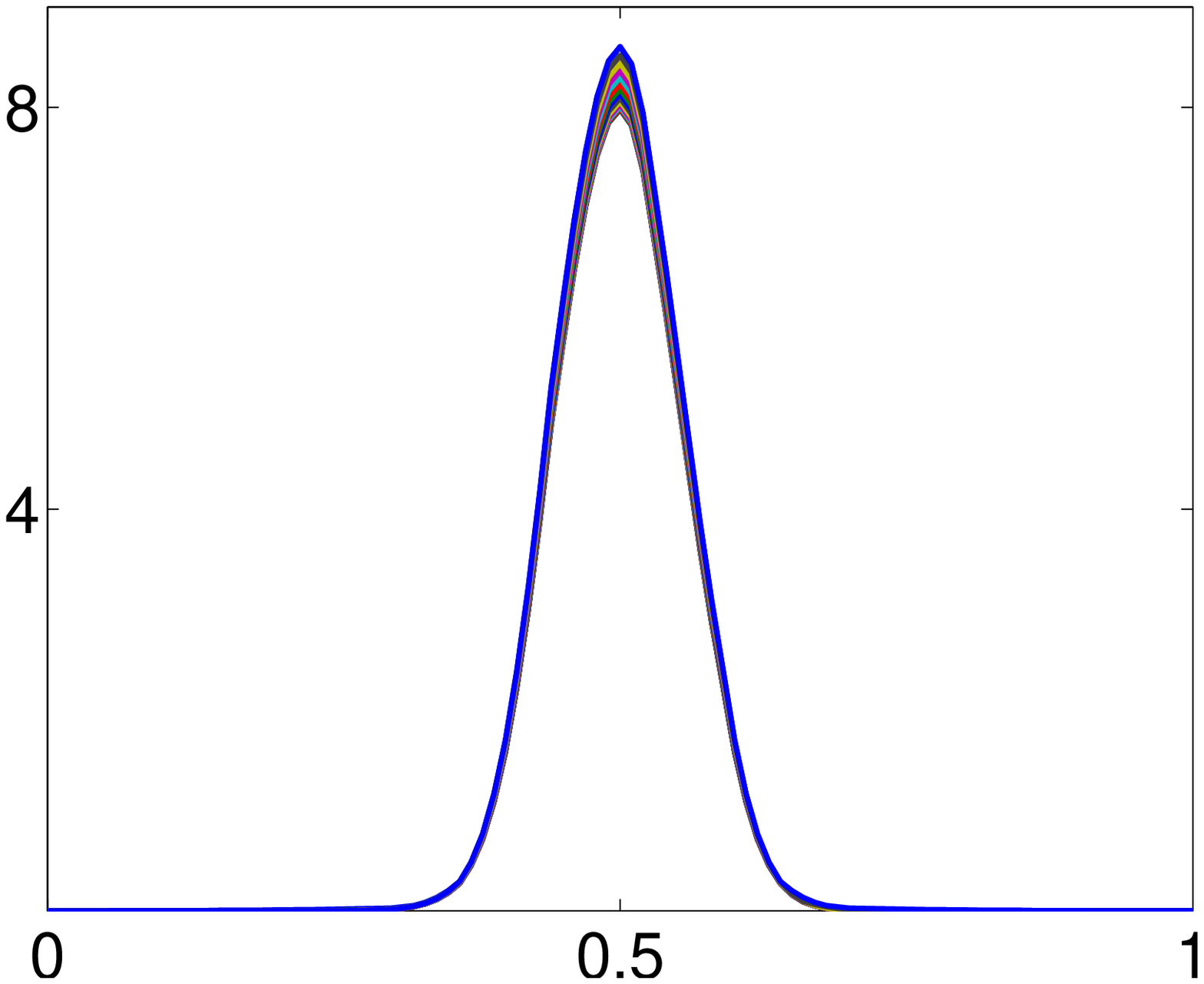} \ &
\includegraphics[height=1.0in]{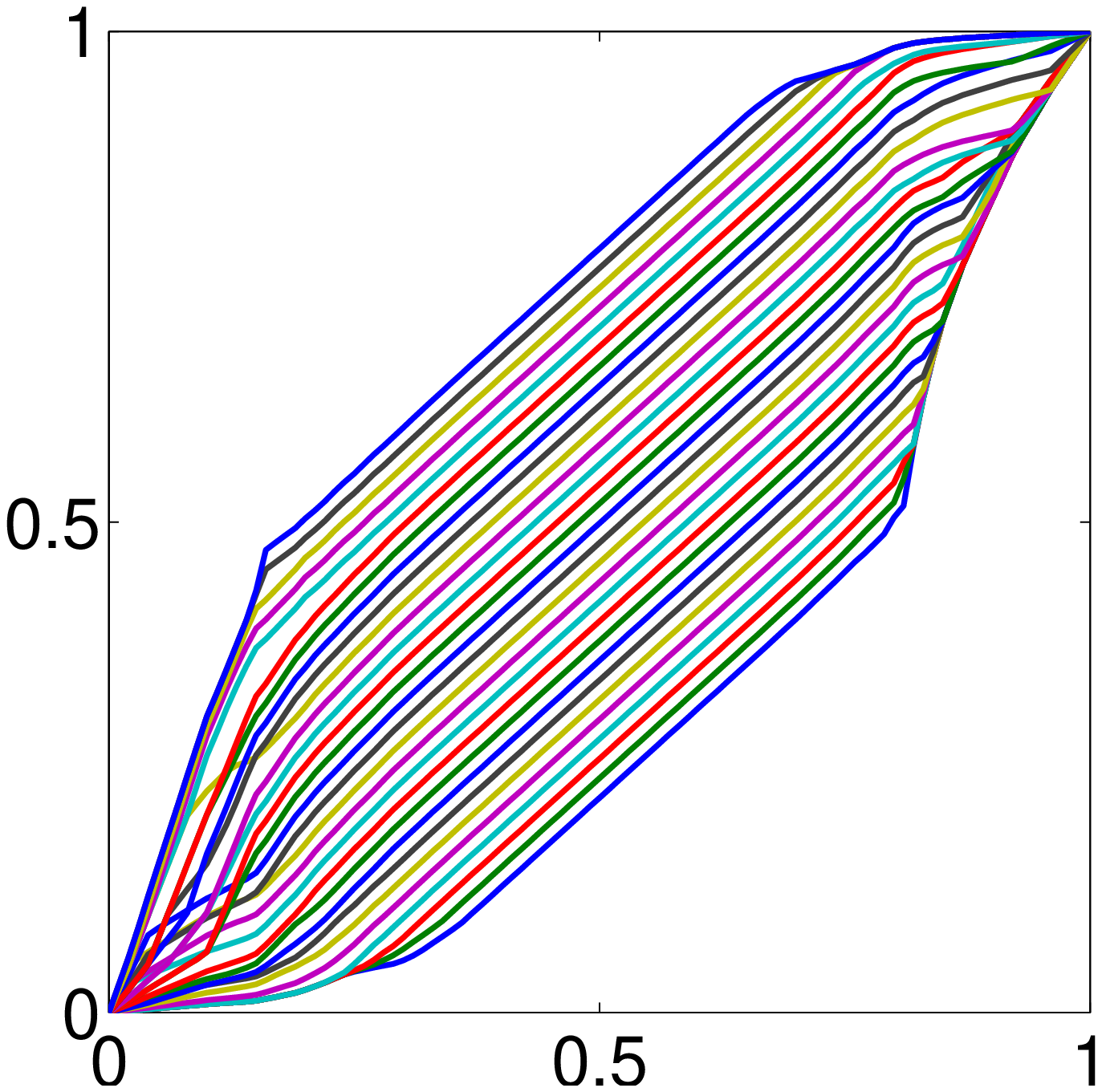} \ &
\includegraphics[height=1.0in]{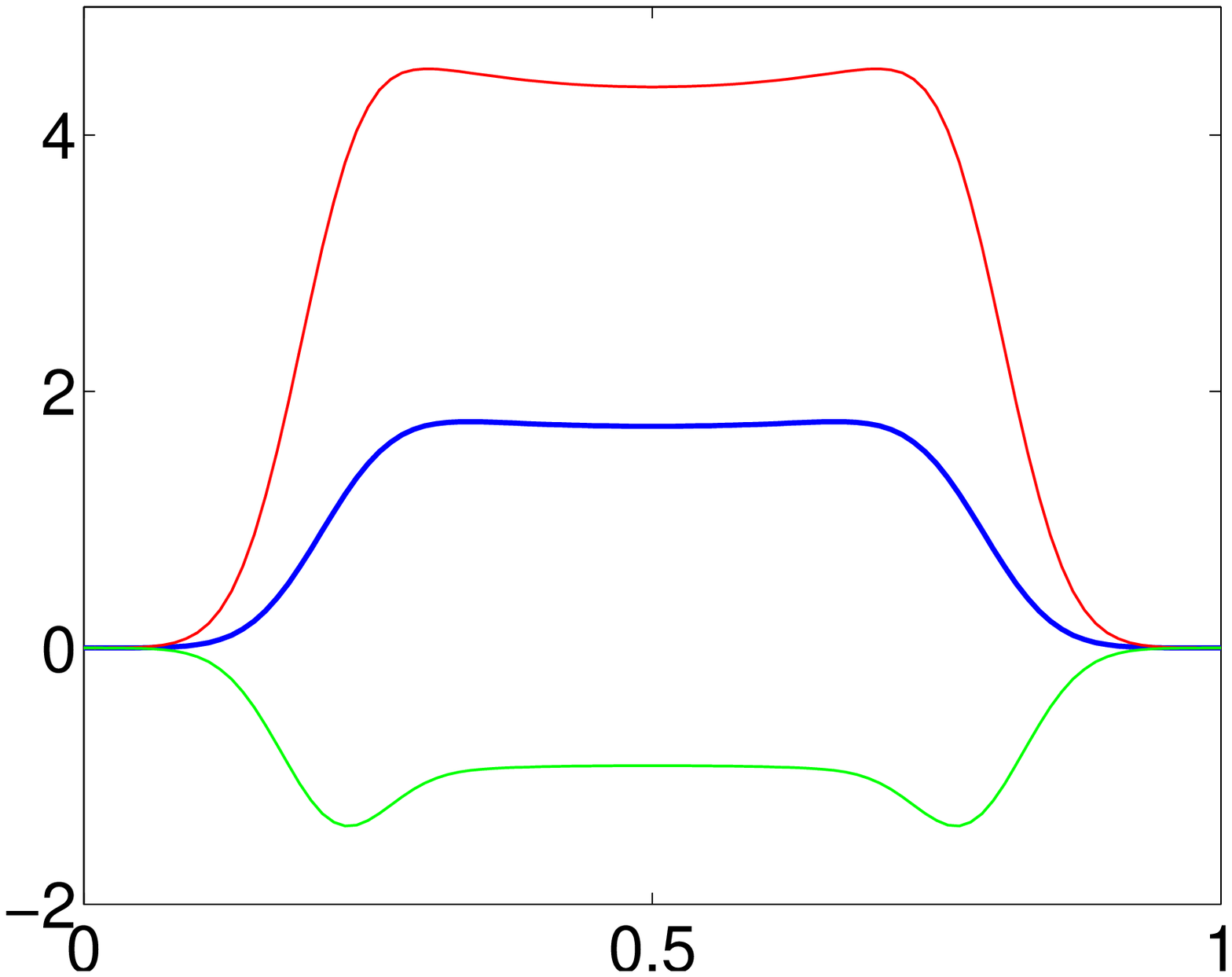} \ &
\includegraphics[height=1.0in]{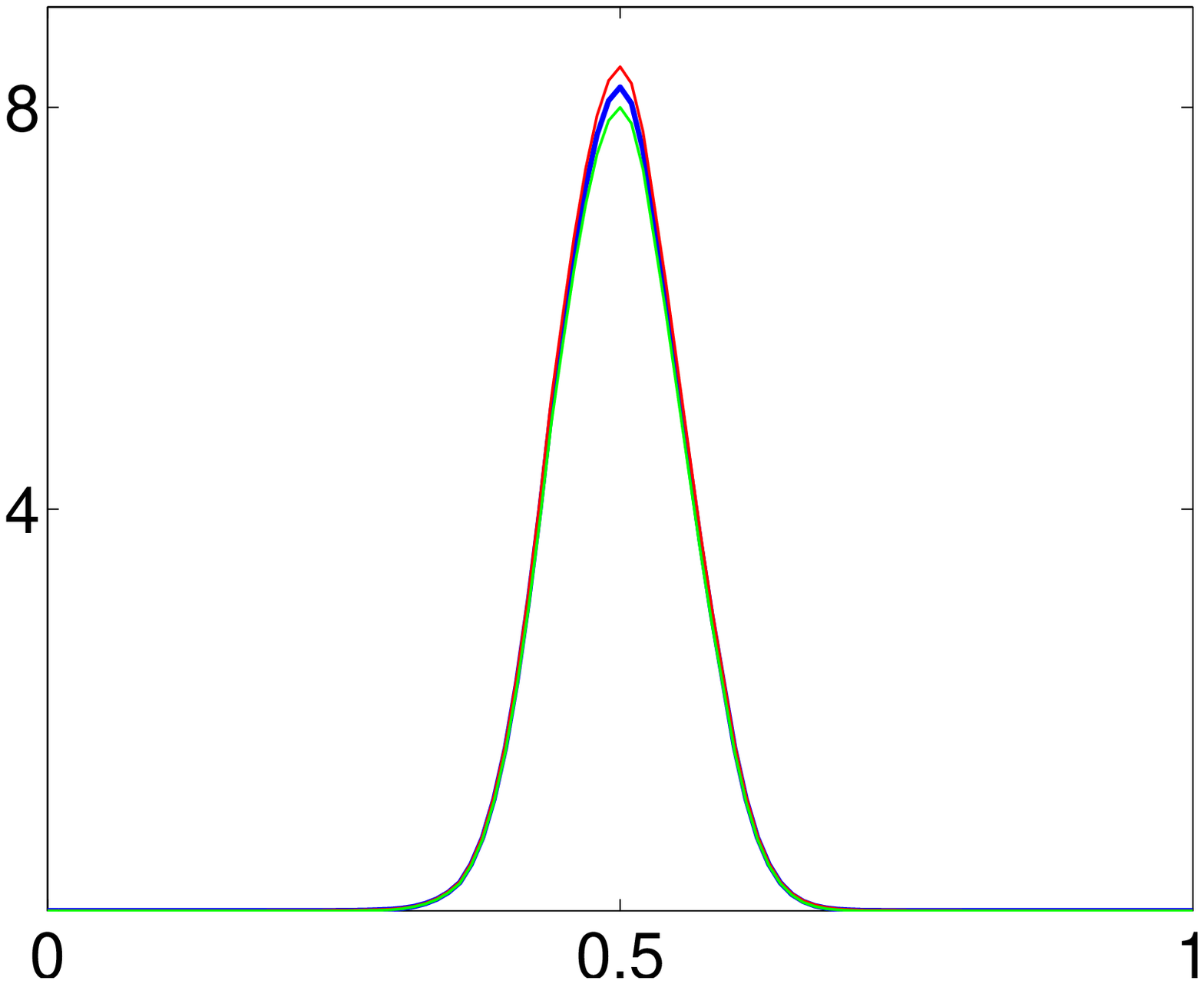}
\end{tabular}
\caption{Results on simulated data 3.} \label{fig:result-sim3}
\end{center}
\end{figure}

\item {\bf Simulated Data 4}: In this case we take a family of multimodal
wave functions with the same shape but different phase variations.
The individual functions are defined on $[0, 9]$ and given by:
$f_i(t) = (1-(\gamma_i(t)/9-0.5)^2)\sin(\pi \gamma_i(t))$,
$i=1,2,\dots, 9$, with the warping functions $\gamma_i(t) =
9({e^{a_i t/9} -1 \over e^{a_i} - 1})$ if $a_i \neq 0$, otherwise
$\gamma_i = \gamma_{id}$. Here $a_i$ are equally spaced between
$-1.5$ and $1.5$ with step size $0.375$. Figure
\ref{fig:result-sim4} shows the original $9$ functions $\{ f_i \}$,
the aligned functions $\{ \tilde{f}_i\}$ (clearly showing the common shape), the warping functions
$\{\gamma_i^*\}$, and the before-and-after cross sectional mean and
standard deviations, again showing the huge difference in apparent amplitude variation 
between aligned and unaligned functions.  In particular, with only the phase variability in the data, our
method has a perfect alignment of given functions.

\begin{figure}
\begin{center}
\begin{tabular}{ccccc}
\hspace*{-0.3in} $\{f_i\}$ & $\{ \tilde{f}_i\}$ & $\{ \gamma_i^*\}$
& mean $\pm$ std, before  & mean $\pm$ std, after
\\
\includegraphics[height=1.0in]{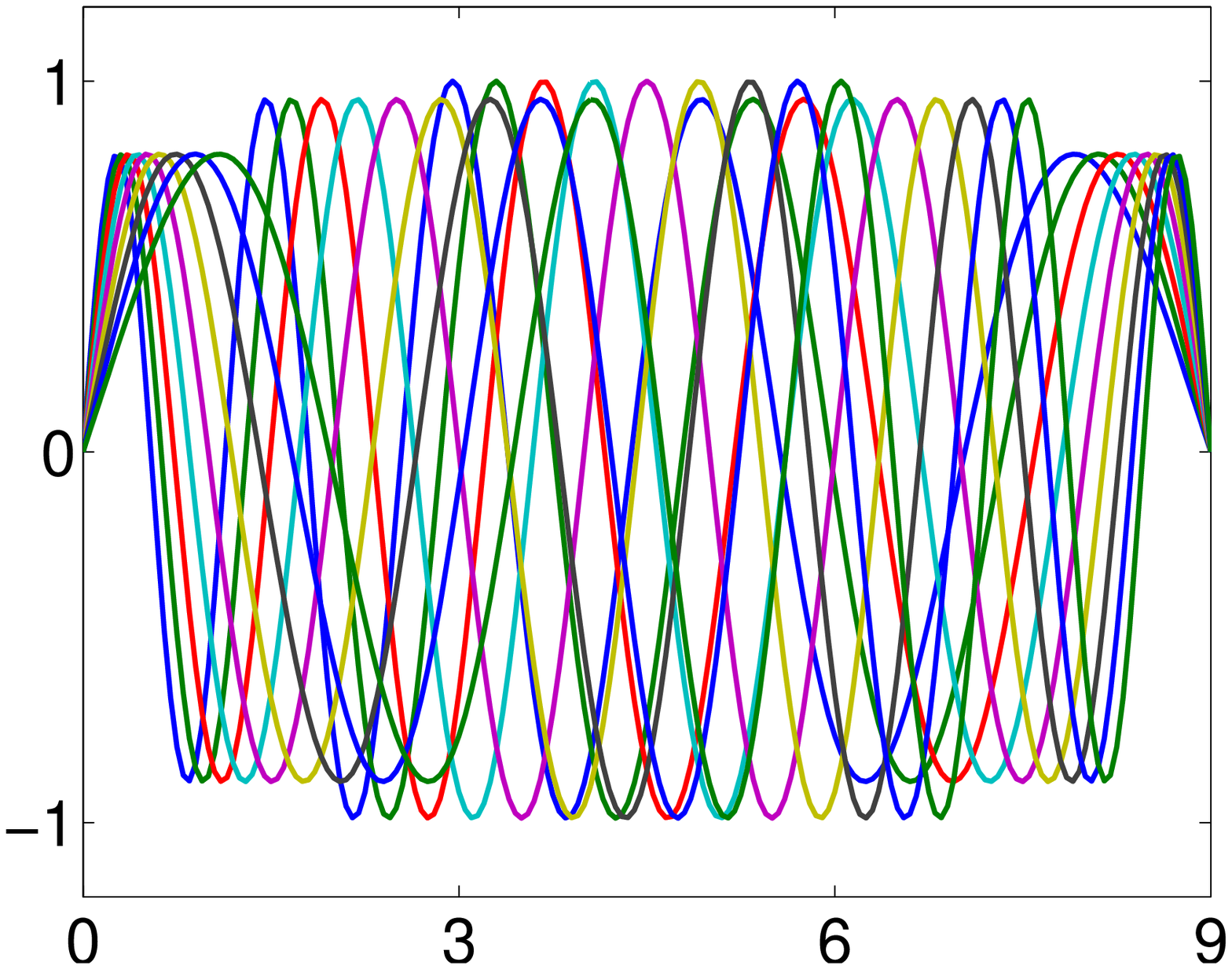} \ &
\includegraphics[height=1.0in]{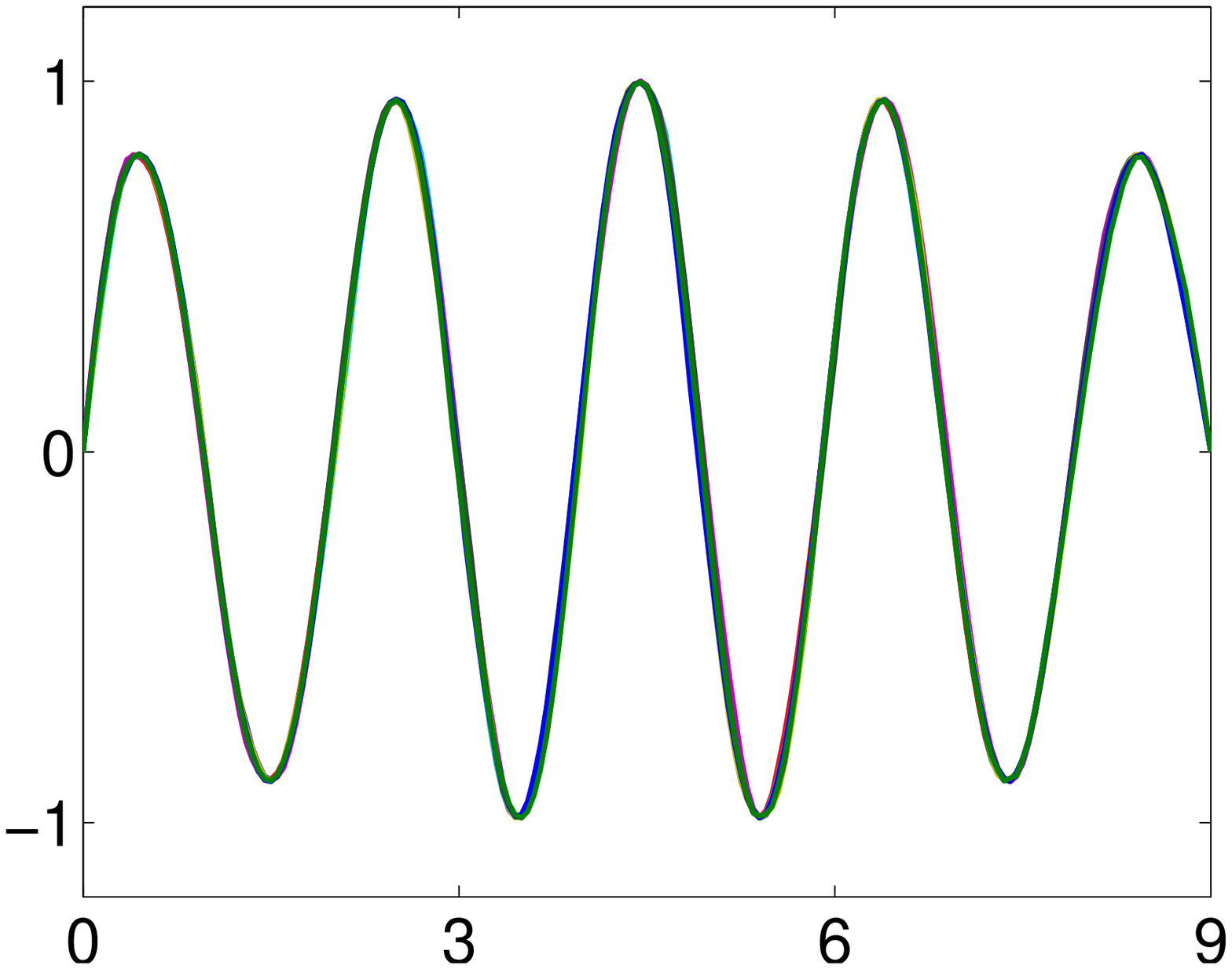} \ &
\includegraphics[height=1.0in]{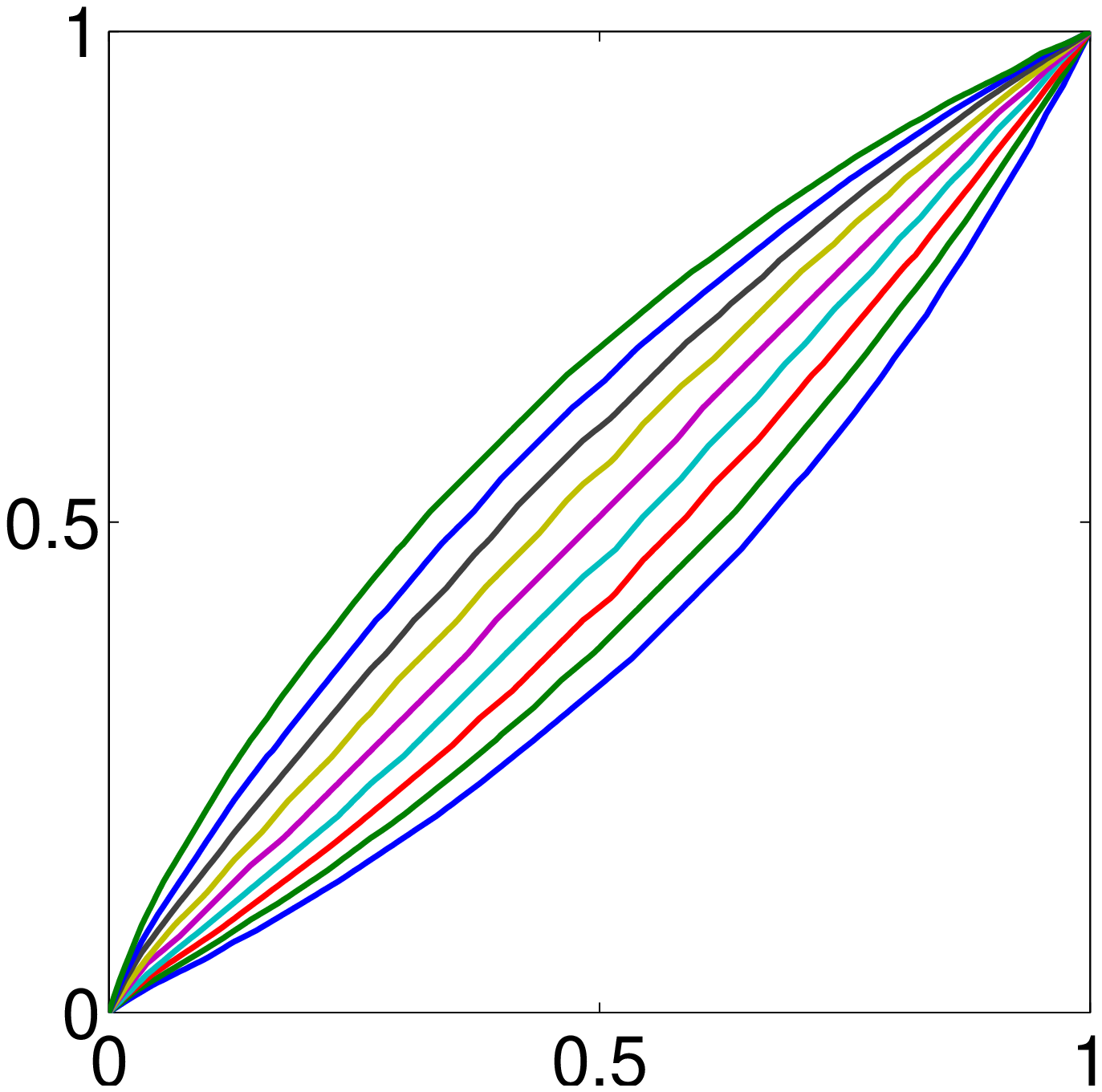} \  &
\includegraphics[height=1.0in]{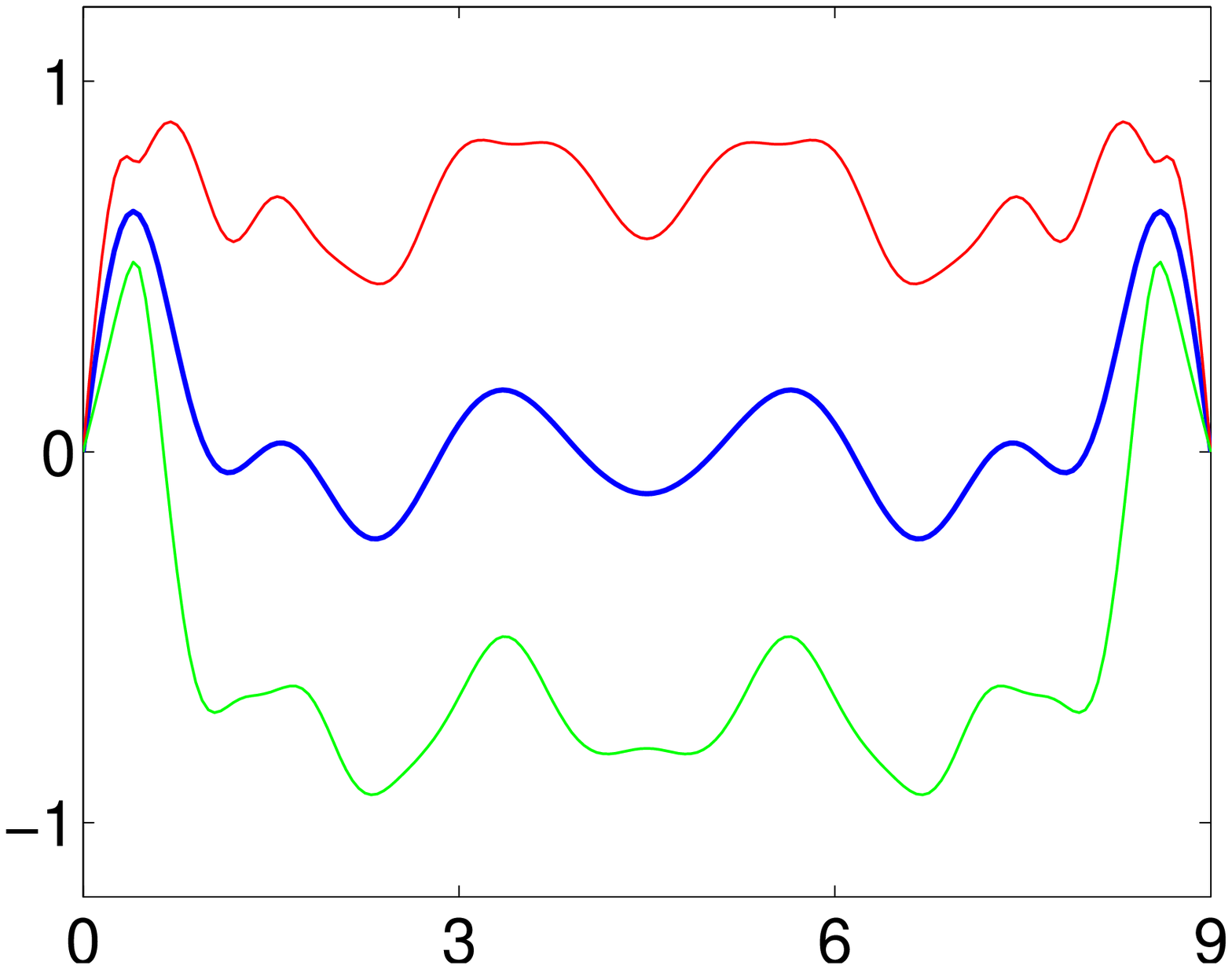} \ &
\includegraphics[height=1.0in]{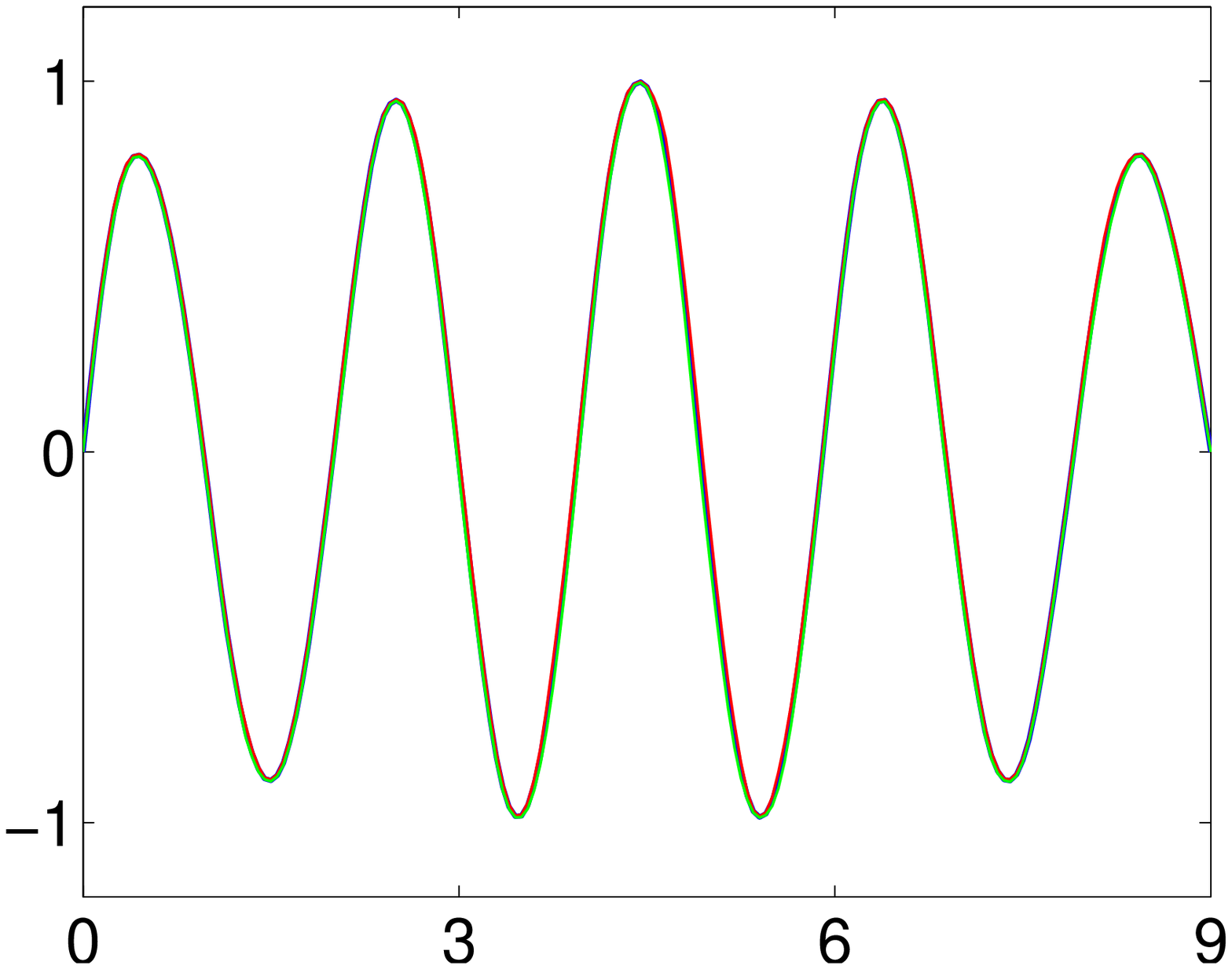}
\end{tabular}
\caption{Results on simulated data 4.} \label{fig:result-sim4}
\end{center}
\end{figure}

\end{enumerate}

\section{Signal Estimation and Estimator Consistency}
In this section we justify the proposed framework by posing and
solving a model-based estimation of alignment. Consider an
observation model $f_i = c_i (g \circ \gamma_i) + e_i$, $i=1, \dots,
n,$ where $g$ is an unknown signal, and $c_i \in \real_+$, $\gamma_i
\in \Gamma$ and $e_i \in {\cal F}$ are random. We will concentrate
on a simpler problem where the observation noise $e_i$ is set to a
constant and, given the observations $\{f_i\}$, the goal is to
estimate the signal $g$ or, equivalently, the warping functions
$\{\gamma_i\}$. This or related problems have been considered
previously by several papers, including
\cite{wang-gasser:97,ramsay-li-RSSB:98}, but we are not aware of any
formal statistical solution. Here we show that the center $\mu_n$,
resulting from the complete alignment algorithm, leads to a
consistent estimator of $g$. The proofs of Lemmas and Corollary are given 
in the appendix.

\begin{theorem} \label{theo:part1}
For a function $g$, consider a sequence of functions $f_i(t) = c_i
g(\gamma_i(t)) + e_i$, where $c_i$ is a positive constant, $e_i$ is
a constant, and $\gamma_i$ is a time warping, $i = 1, \cdots, n$.
Denote by $q_g$ and $q_i$ the SRVFs of $g$ and $f_i$, respectively,
and let $\bar{s} = { 1 \over n} \sum_{i=1}^n \sqrt{c_i}$. Then, the
Karcher mean of $\{ [q_i],i=1,2,\dots,n\}$ in ${\cal S}$ is
$\bar{s}[q_g]$. That is,
$$
[\mu]_n \equiv \argmin_{[q]} \left(\sum_{i=1}^N d^2([q_i], [q])
\right)   =  \bar{s} [q_g] =  \bar{s} \{(q_g, {\gamma}), {\gamma}
\in \Gamma\}\ .
$$
\end{theorem}
We will prove this theorem in two steps. First we establish the following
useful result.

\begin{lemma} \label{lem:sec}
For any $q_1, q_2 \in \ltwo$ and a constant $c > 0$, we have
$\argmin_{\gamma \in \Gamma} \| q_1- (q_2, \gamma)\|
 = \argmin_{\gamma \in \Gamma} \| cq_1 - (q_2, \gamma)\|$.
\end{lemma}

\begin{corollary} \label{cor:align} For any function $q \in \ltwo$ and constant $c>0$,
we have $\gamma_{id} \in \argmin_{\gamma \in \Gamma} \| cq - (q,
\gamma)\|$.  Moreover, if the set $\{t\in[0,1]| q(t)=0\}$ has
(Lebesgue) measure 0, $\gamma_{id} = \argmin_{\gamma \in \Gamma} \|
cq - (q, \gamma)\|$.
\end{corollary}

\noindent Now we get back to the proof of Theorem \ref{theo:part1}.
The SRVF of the function $f_i = c_i (g \circ \gamma_i) + e_i$ is
given by $q_i = \sqrt{c_i} (q_g, \gamma_i)$, $i = 1, \cdots, n$. For
any $q$, we get
\begin{eqnarray*}
d([q_i], [q]) &=& d([\sqrt{c_i}(q_g, \gamma_i)], [q]) =
\inf_{\gamma} \| \sqrt{c_i}(q_g, \gamma_i) - (q, \gamma) \| \\
&=& \inf_{\gamma} \| \sqrt{c_i}q_g - (q, \gamma \circ \gamma_i^{-1})
\|  = \inf_{\gamma} \| \sqrt{c_i}q_g - (q, \gamma) \|= d([\sqrt{c_i}
q_g], [q]).
\end{eqnarray*}
In the last line, the first equality is based on the isometry of the
group action of $\Gamma$ on $\ltwo$ and the second equality is based on
the group structure of $\Gamma$.

For any given $q$, let $\gamma^* \in \argmin_{\gamma \in \Gamma} \|
q_g - (q,\gamma)\|$. Then, using Lemma \ref{lem:sec}, $\gamma^* \in
\argmin_{\gamma \in \Gamma} \| \sqrt{c_i} q_g - (q,\gamma)\|$.
Therefore,
\begin{equation*}
\sum_{i=1}^n d^2([q_i], [q]) = \sum_{i=1}^n d^2([\sqrt{c_i}q_g],
[q]) = \sum_{i=1}^n \| \sqrt{c_i}q_g - (q,\gamma^*) \|^2 \ge
\sum_{i=1}^n \| \sqrt{c_i}q_g - \bar {s} q_g \|^2\ .
\end{equation*}
The last inequality comes from the fact that $\bar {s} q_g$ is
simply the mean of $\{\sqrt{c_i} q_g\}$ in $\ltwo$ space.  The
equality holds if and only if $(q, \gamma^*) = \bar {s} q_g$ or $q =
(\bar {s} q_g, {\gamma^*}^{-1})$.

Actually, for any element of $[\bar{s} q_g]$, say $(\bar{s} q_g,
\gamma)$ for any $\gamma \in \Gamma$, we have
$$\sum_{i=1}^n d^2([q_i], [(\bar{s} q_g, \gamma)]) = \sum_{i=1}^n d^2([\sqrt{c_i}q_g], [\bar{s} q_g])
= \sum_{i=1}^n \| \sqrt{c_i}q_g - \bar {s} q_g\|^2\ .$$ Therefore,
$\{(\bar{s} q_g, \gamma) | \gamma \in \Gamma\}$ is the unique
solution to the Karcher mean $[\mu]_n \equiv \argmin_{[q]}
\sum_{i=1}^n d^2([q_i], [q]). \ \ \ \ \Box$

Next, we present a simple fact about the Karcher mean of warping functions,
where the Karcher mean is given in Definition \ref{def:mean-Gamma}.
\begin{lemma} \label{lem:mean-gamma}
Given a set $\{\gamma_i \in \Gamma| i = 1, ..., n\}$ and a $\gamma_0
\in \Gamma$,  if the Karcher mean of $\{\gamma_i\}$ is $\bar
\gamma$, then the Karcher mean of $\{\gamma_i \circ \gamma_0\}$ is
$\bar \gamma \circ \gamma_0$.
\end{lemma}

Theorem \ref{theo:part1} ensures that $[\mu]_n$ belongs to the orbit of $[q_g]$ (up to a scale 
factor) but we are interested in estimating $g$ itself, rather than its orbit. 
Since we can write $g \circ \gamma_i$ as $(g \circ \gamma_a) \circ (\gamma_a^{-1} \circ \gamma_i)$, 
for any $\gamma_a \in \Gamma$, the function $g$ is not identifiable unless we impose 
an additional constraint on $\gamma_i$. The same goes for the random variables $c_i$ and $e_i$. 
Under the assumption that the population means of $\gamma_i^{-1}$, $c_i$, and 
$e_i$ are known, we will show in two steps that  
Algorithm 3, that finds the center of the orbit
$[\mu]_n$, results in a consistent estimator for $g$.
\begin{theorem} \label{theo:part2}
Under the same conditions as in Theorem \ref{theo:part1}, let $\mu =
(\bar{s} q_g,\gamma_0)$, for $\gamma_0 \in \Gamma$, denote an
arbitrary element of the Karcher mean class $[\mu]_n =
\bar{s}[q_g]$. Assume that the set $\{t\in[0,1] | \dot g(t) =0\}$ has
Lebesgue measure zero.  If the population Karcher mean of $\{\gamma_i^{-1}\}$ is
$\gamma_{id}$, then the center of the orbit $[\mu]_n$, denoted by  $\mu_n$, satisfies
$\lim_{n \rightarrow \infty} {\mu}_n = E(\bar s)q_g$.
\end{theorem}
{\bf Proof:} In Algorithm 3, we first compute
$\tilde \gamma_i = \argmin_{\gamma} \|(q_i, \gamma) - \mu \| =
\argmin_{\gamma} \|(\sqrt{c_i} (q_g, \gamma_i), \gamma) - (\bar{s}
q_g,  \gamma_0)\| 
= \argmin_{\gamma} \|(\sqrt{c_i}q_g, \gamma_i \circ \gamma \circ
\gamma_0^{-1}) - \bar{s} q_g \|$.
Since the set $\{t\in[0,1] | \dot g(t) = 0\}$ has
measure zero, the set $\{t\in[0,1]| q_g(t)=0\}$ also has measure zero.
Using Corollary 1, this above distance is uniquely minimized when
$\gamma_i \circ \tilde \gamma_i \circ \gamma_0^{-1} = \gamma_{id}$,
or $\tilde \gamma_i = \gamma_i^{-1} \circ \gamma_0$. Denote the
Karcher mean of these warping functions $\{\tilde \gamma_i\}$ by
$\bar{\gamma}_{n}$. Applying the inverse of this $\bar{\gamma}_{n}$
to $\mu$, we get $\mu_n = (\mu, \bar{\gamma}_{n}^{-1})$. As $n
\rightarrow \infty$, the Karcher mean of $\tilde \gamma_i$ converges
to its population mean which, by Lemma \ref{lem:mean-gamma},  is
$\gamma_0$.  Thus, $\mu_n \stackrel{n \rightarrow
\infty}{\longrightarrow} E(\bar{s})((q_g, \gamma_0), \gamma_0^{-1})
= E(\bar{s})q_g$.
  \ \ \ \ $\Box$

This result shows that asymptotically one can recover the SRVF of
the original signal using the Karcher mean of the SRVFs of the
observed signals. Of course, one is really interested in the signal
$g$ itself, rather than its SRVF. One can reconstruct $g$ using
aligned functions $\{\tilde f_i\}$ generated by the Alignment
Algorithm in Section \ref{sec:KM}.  As discussed above, we
further assume the population mean of $e_i$ is known.
\begin{theorem} \label{theo:part3}
Under the same conditions as in Theorem \ref{theo:part2}, let
$\gamma_i^* = \argmin_{\gamma} \|(q_i, \gamma) - \mu_n \|$ and
$\tilde f_i = f_i \circ \gamma_i^*$.  If we denote $\bar {c} =
\frac{1}{n}\sum_{i=1}^n c_i$ and $\bar {e} = \frac{1}{n}\sum_{i=1}^n
e_i$, then $\lim_{n \rightarrow \infty} \frac{1}{n}\sum_{i=1}^n
\tilde f_i = E(\bar c)g + E(\bar e)$.
\end{theorem}
{\bf Proof:} In the proof for Theorem \ref{theo:part2}, $\tilde
\gamma_i = \argmin_{\gamma} \|(q_i, \gamma) - \mu \| =
\argmin_{\gamma} \|(q_i, \gamma) - (\mu_n, \bar \gamma_n) \|$. Hence
$\gamma_i^* = \argmin_{\gamma} \|(q_i, \gamma) - \mu_n \|= \tilde
\gamma_i \circ \bar \gamma_n^{-1} = \gamma_i^{-1} \circ \gamma_0
\circ \bar \gamma_n^{-1}$. This implies that $ \tilde f_i = f_i
\circ \gamma_i^* = (c_i (g \circ \gamma_i) + e_i)\circ
(\gamma_i^{-1} \circ \gamma_0 \circ \bar \gamma_n^{-1}) = c_i (g
\circ \gamma_0 \circ \bar \gamma_n^{-1}) + e_i$. As $
\bar{\gamma}_{n} \rightarrow \gamma_0$ when $n \rightarrow \infty$,
we have
$$
\lim_{n \rightarrow \infty} \frac{1}{n}\sum_{i=1}^n \tilde f_i =
E(\bar c)(g \circ \gamma_0) \circ \gamma_0^{-1} + E(\bar e)= E(\bar
c)g + E(\bar e).
$$

\begin{figure}
\begin{center}
\begin{tabular}{ccccc}
$g$ & $\{f_i\}$ & $\{\tilde f_i\}$ & estimate of $g$ &
error w.r.t. $n$ \\
\includegraphics[height=1.0in]{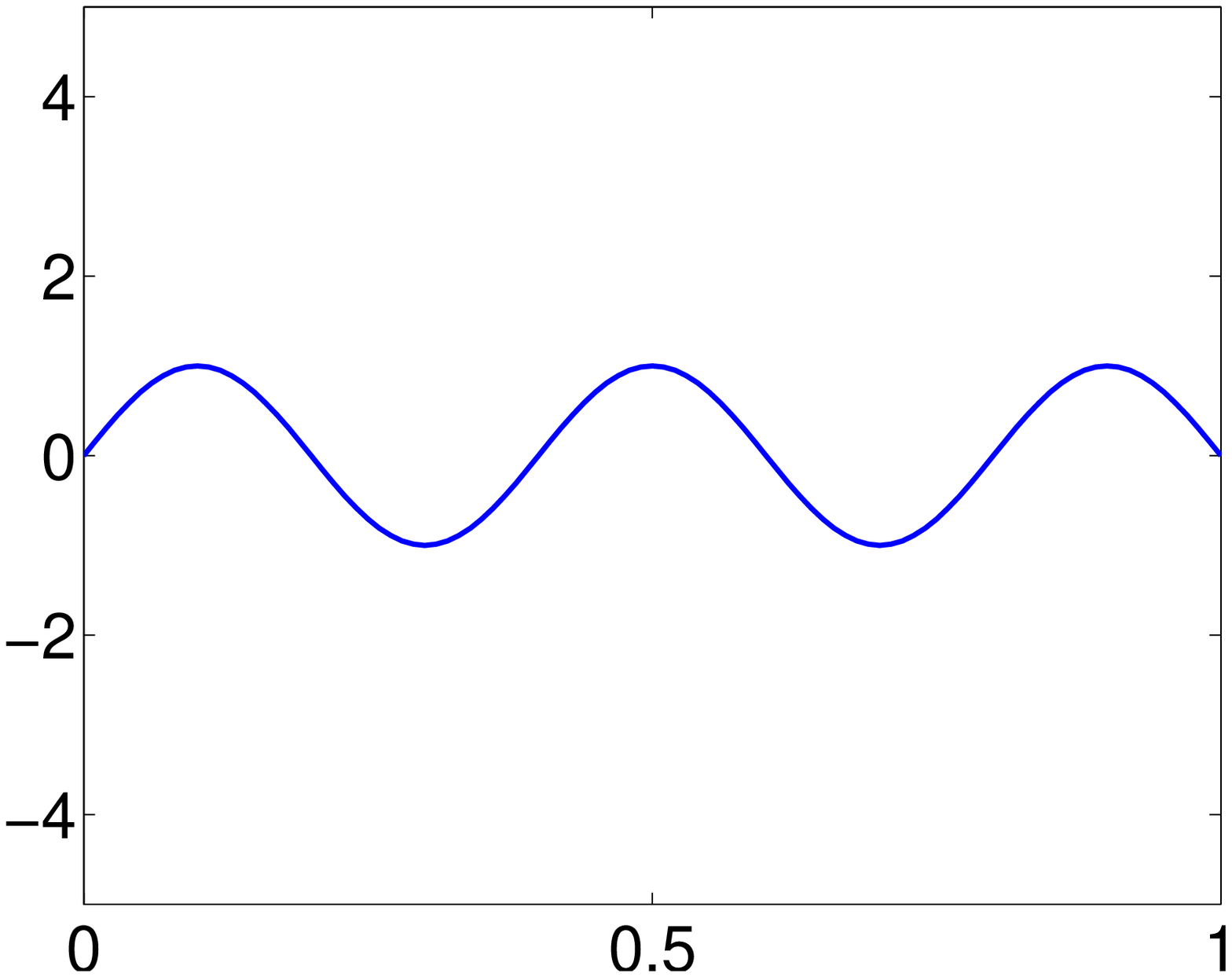} \ &
\includegraphics[height=1.0in]{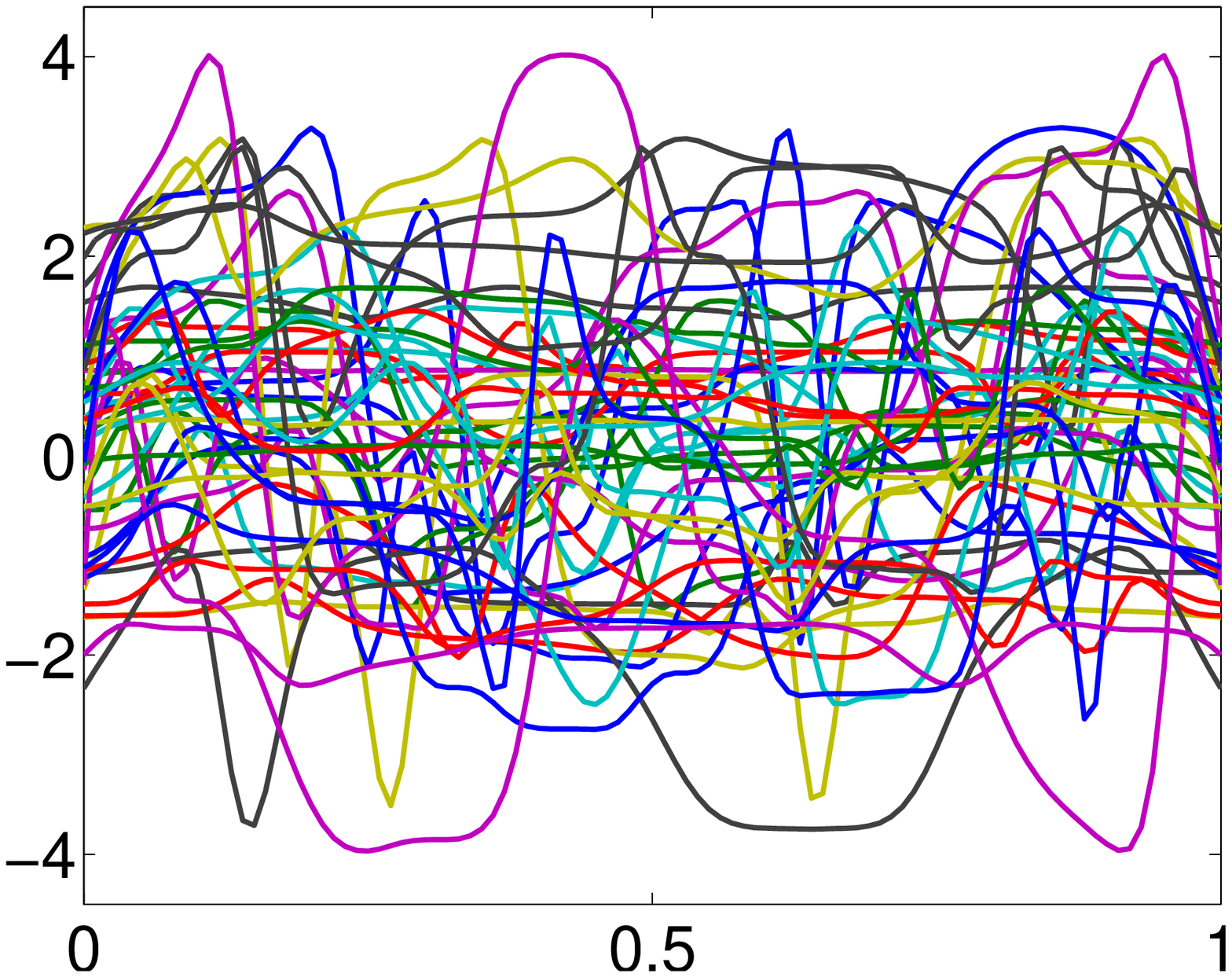} \ &
\includegraphics[height=1.0in]{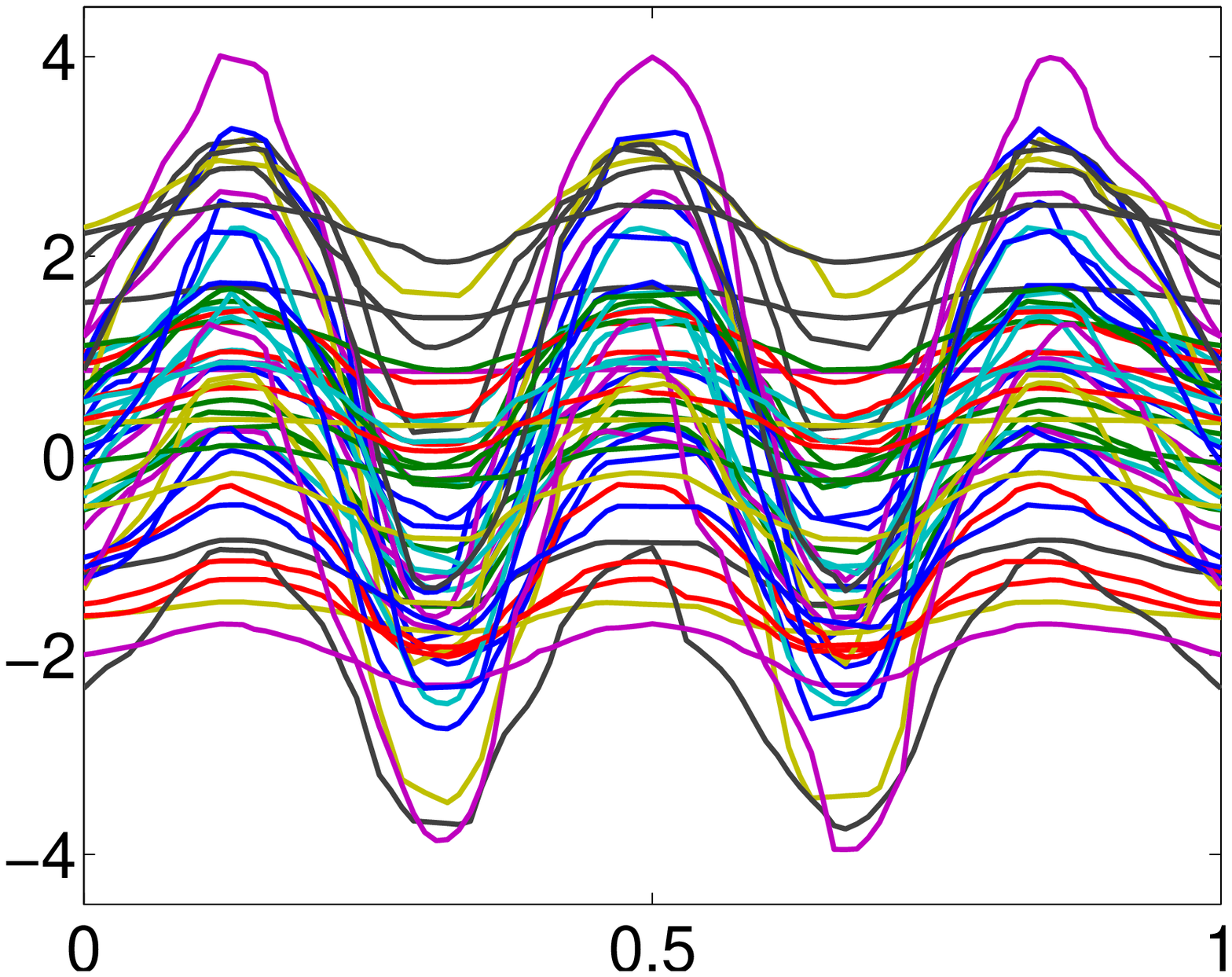} \ &
\includegraphics[height=1.0in]{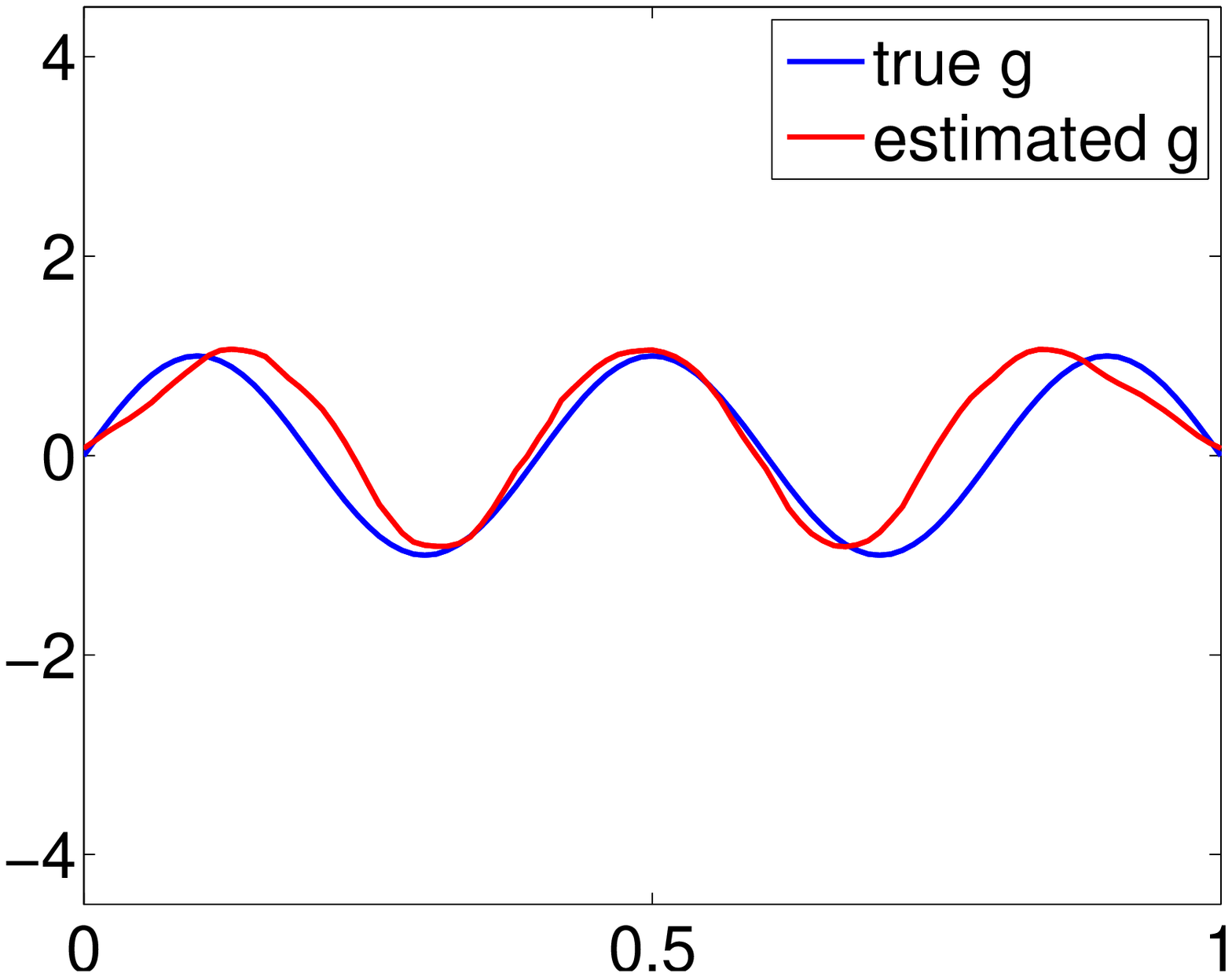} \ &
\includegraphics[height=1.0in]{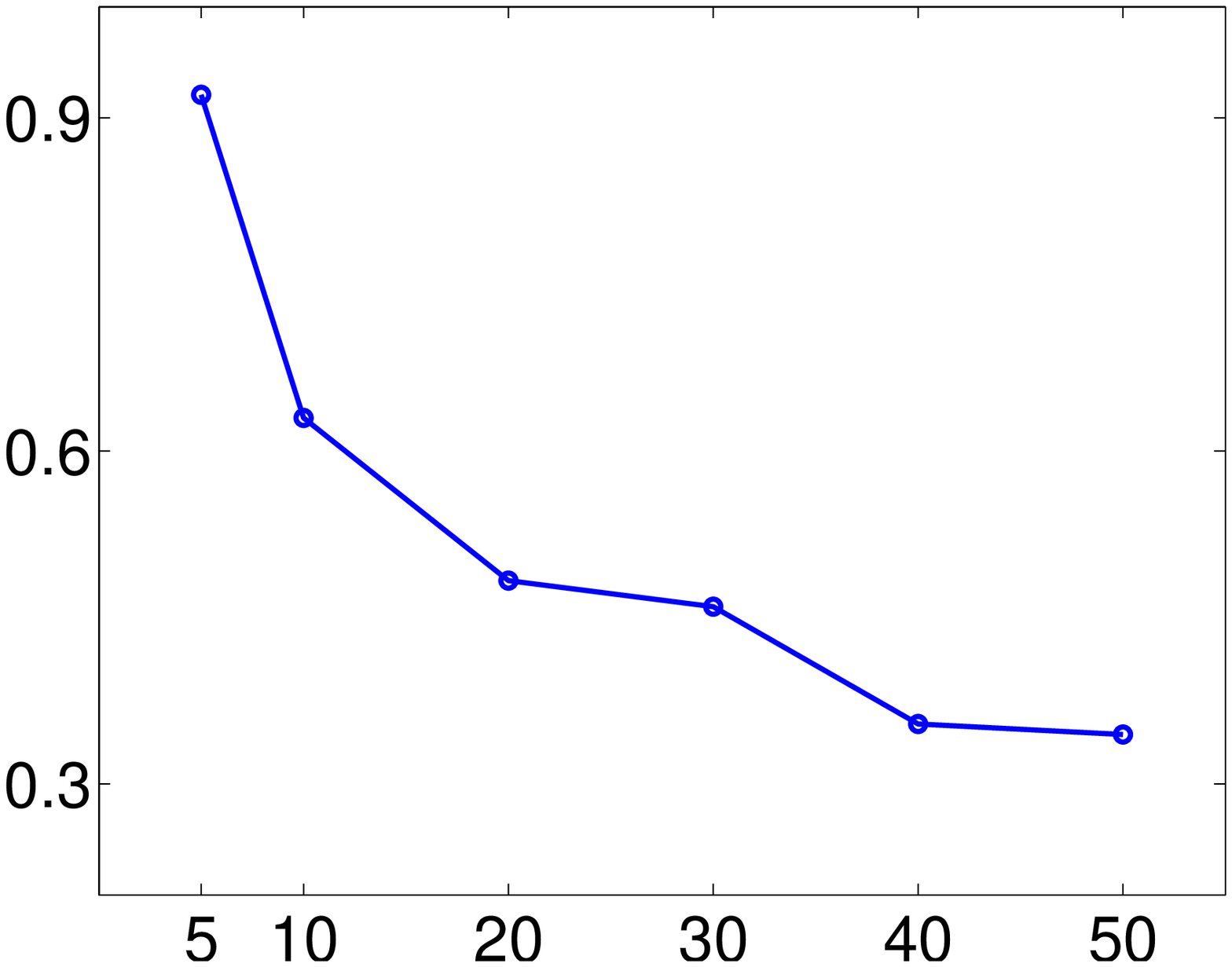}
\end{tabular}
\caption{Example of consistent estimation.} \label{fig:consistency}
\end{center}
\end{figure}

\noindent {\bf Illustration.} We illustrate the estimation process
using an example with $g(t) = \sin(5\pi t)$,  $t \in [0,1]$.  We
randomly generate $n = 50$ warping functions $\{\gamma_i\}$ such
that $\{\gamma_i^{-1}\}$ are {\it i.i.d}
 with mean $\gamma_{id}$.  We also generate
{\it i.i.d} sequences $\{c_i\}$ and $\{e_i\}$ from the exponential
distribution with mean 1 and the standard normal distribution,
respectively.  Then we compute functions $f_i = c_i (g \circ r_i) +
e_i$ to form the functional data. In Fig. \ref{fig:consistency}, the
first panel shows the function $g$, and the second panel shows the data
$\{f_i\}$. The Alignment Algorithm in Section \ref{sec:KM} results
in the aligned functions  $\{\tilde f_i = f_i \circ \gamma_i^* \}$
that are are shown  in the third panel in Fig.
\ref{fig:consistency}. Using Theorem \ref{theo:part3}, the original
signal $g$ can be estimated by $(\frac{1}{n}\sum_{i=1}^n \tilde f_i
- E(\bar e))/E(\bar c)$.  In this case, $E(\bar c)) = 1, E(\bar e) =
0$. This estimated $g$ (red) as well as the true $g$ (blue) are shown in the
fourth panel. Note that the estimate is reasonably
successful despite large variability in the raw data. Finally, we examine the performance of the
estimator with respect to the sample size, by performing this estimation for  $n$ equal to $5, 10, 20,
30$, and $40$.  The estimation
errors, computed using the $\ltwo$ norm between estimated $g$'s and the true $g$, are
shown in the last panel.
As expected from the earlier theoretical development, this
estimate converges to the true $g$ when the sample size $n$ grows large.

\section{Experimental Evaluation of Function Alignment} \label{sec:experiment}
In this section we take functional data from several
application domains and analyze them using the framework developed
in this paper.  Specifically, we will focus on function alignment
and comparison of alignment performance with some previous methods on
several datasets.

\subsection{Applications on real data}
We start with demonstrations of the proposed framework on some
well known functional data. \\

\noindent 1. {\bf Berkeley Growth Data}: As a first example, we consider the
Berkeley growth dataset for 54 female and 39 male subjects.  For
better illustrations, we have used the first derivatives of the
growth curves as the functions $\{ f_i\}$ in our analysis. (In this
case, since SRVF is based on the first derivative of $f$, we
actually end up using the second derivatives of the growth
functions.)

\begin{figure}
\begin{center}
\begin{tabular}{cc|cc}
\multicolumn{2}{c}{Female Data} & \multicolumn{2}{c}{Male Data} \\
original data & aligned functions & original data & aligned functions\\
\includegraphics[height=1.0in]{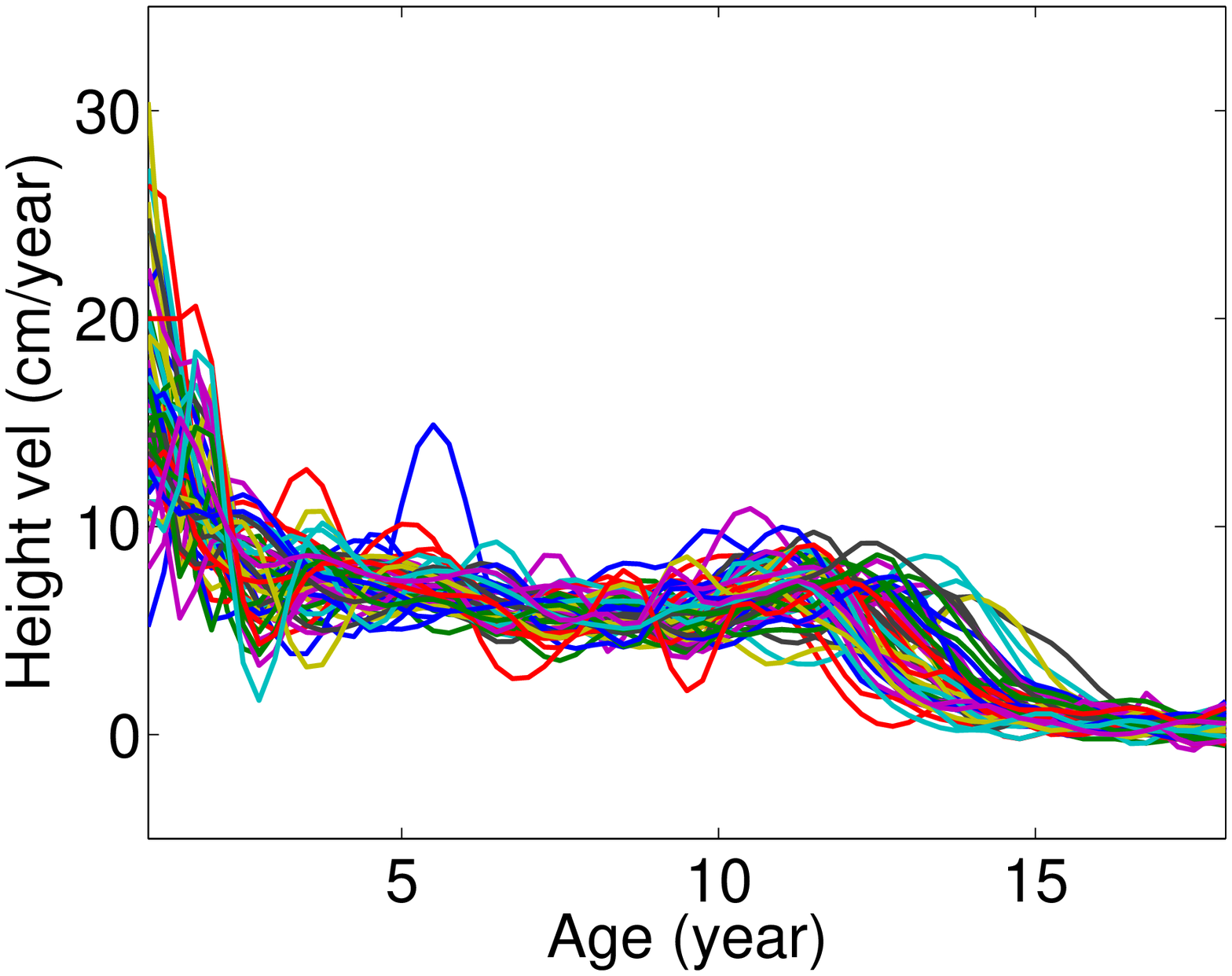} \ &
\includegraphics[height=1.0in]{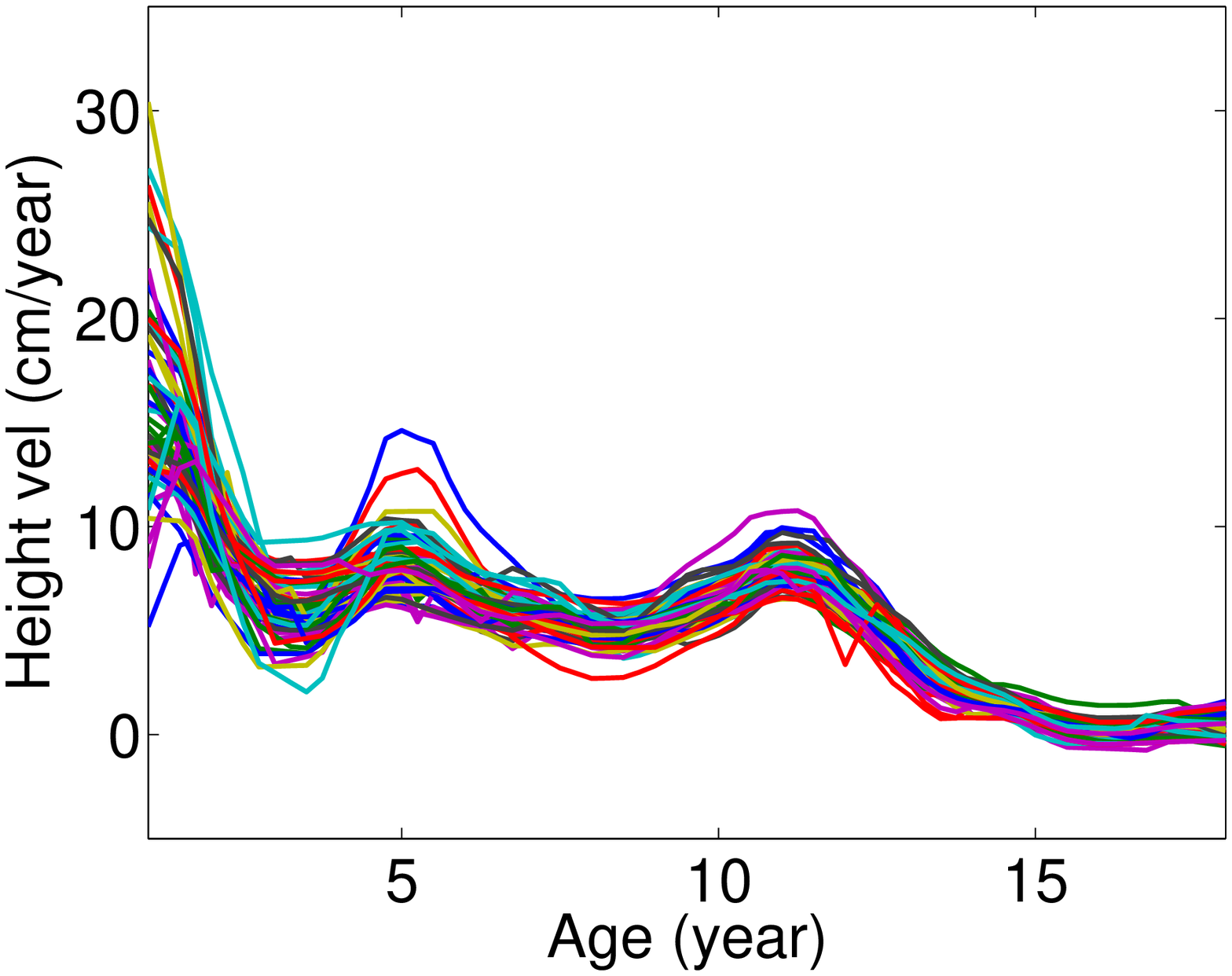} \ &
\includegraphics[height=1.0in]{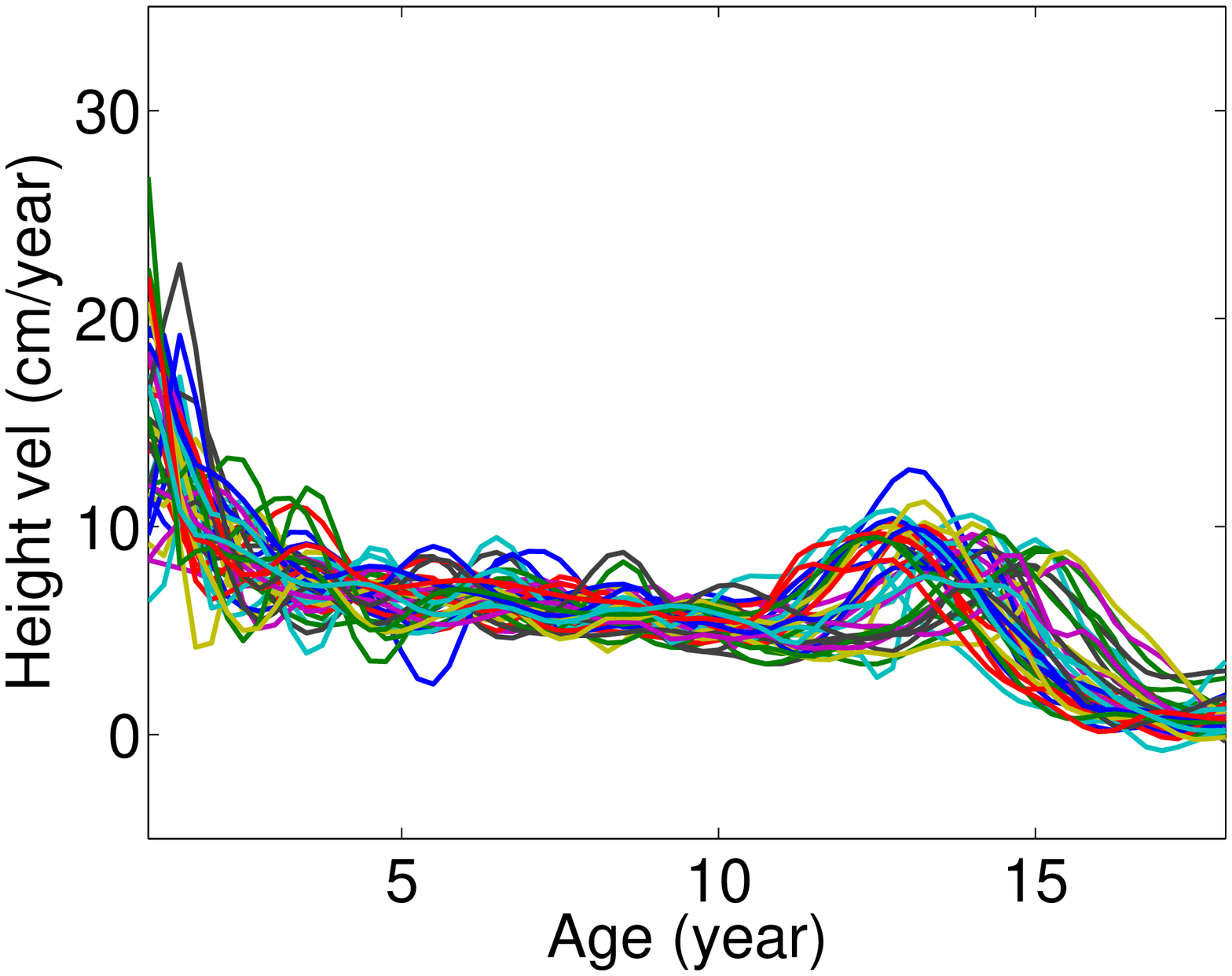} \ &
\includegraphics[height=1.0in]{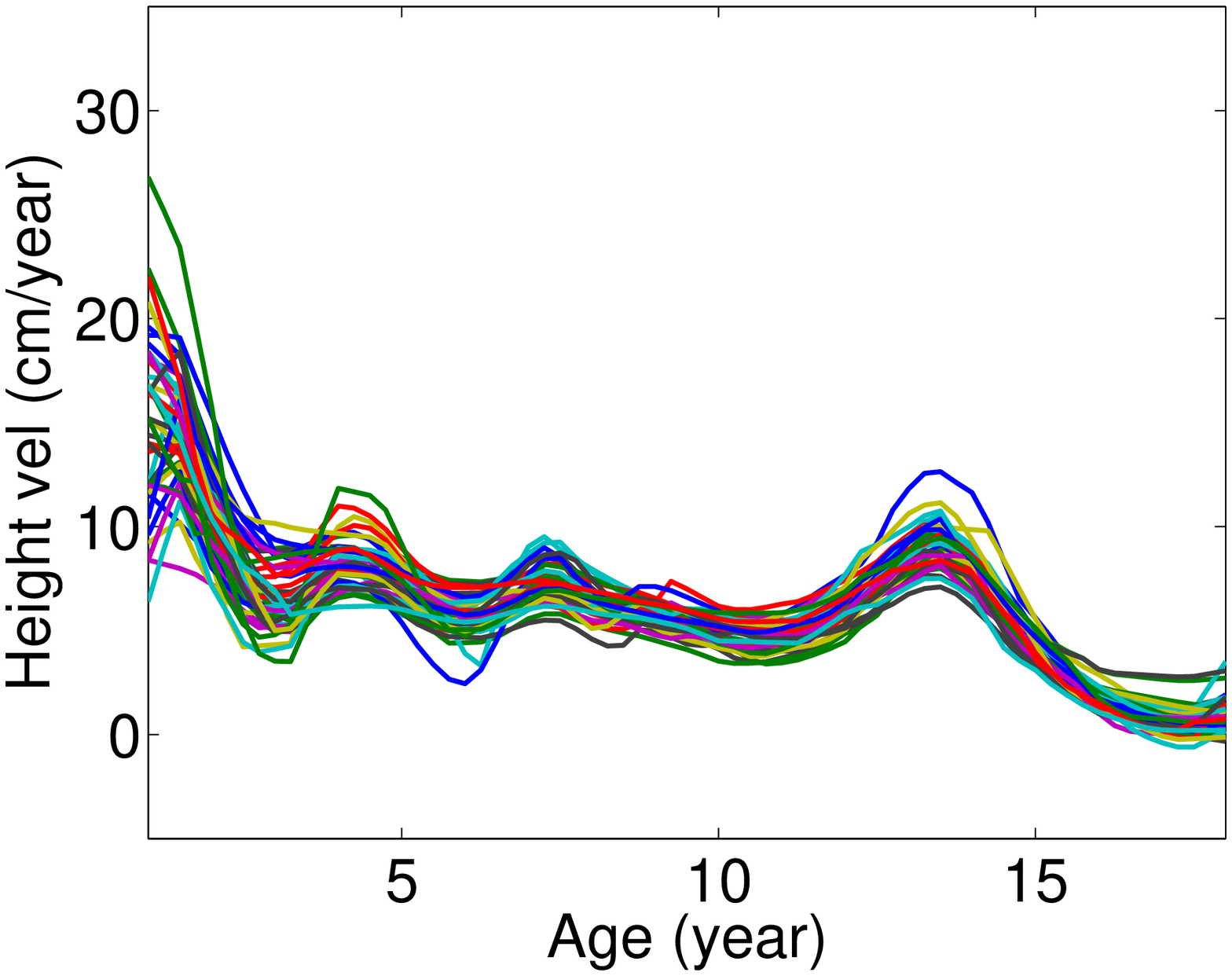}\\
mean $\pm$ std, before & mean $\pm$ std, after &  mean $\pm$ std, before & mean $\pm$ std, after \\
\includegraphics[height=1.0in]{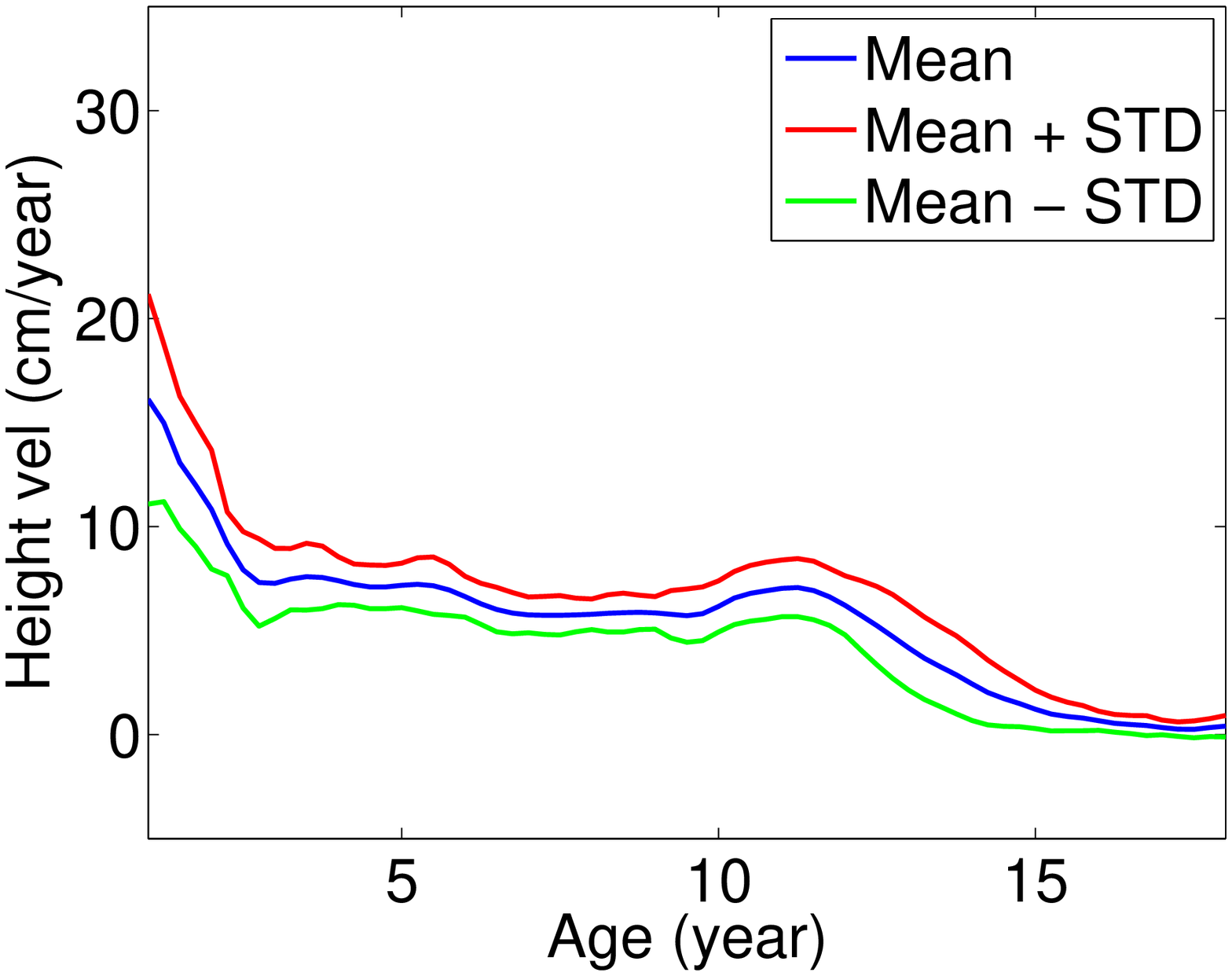} \ &
\includegraphics[height=1.0in]{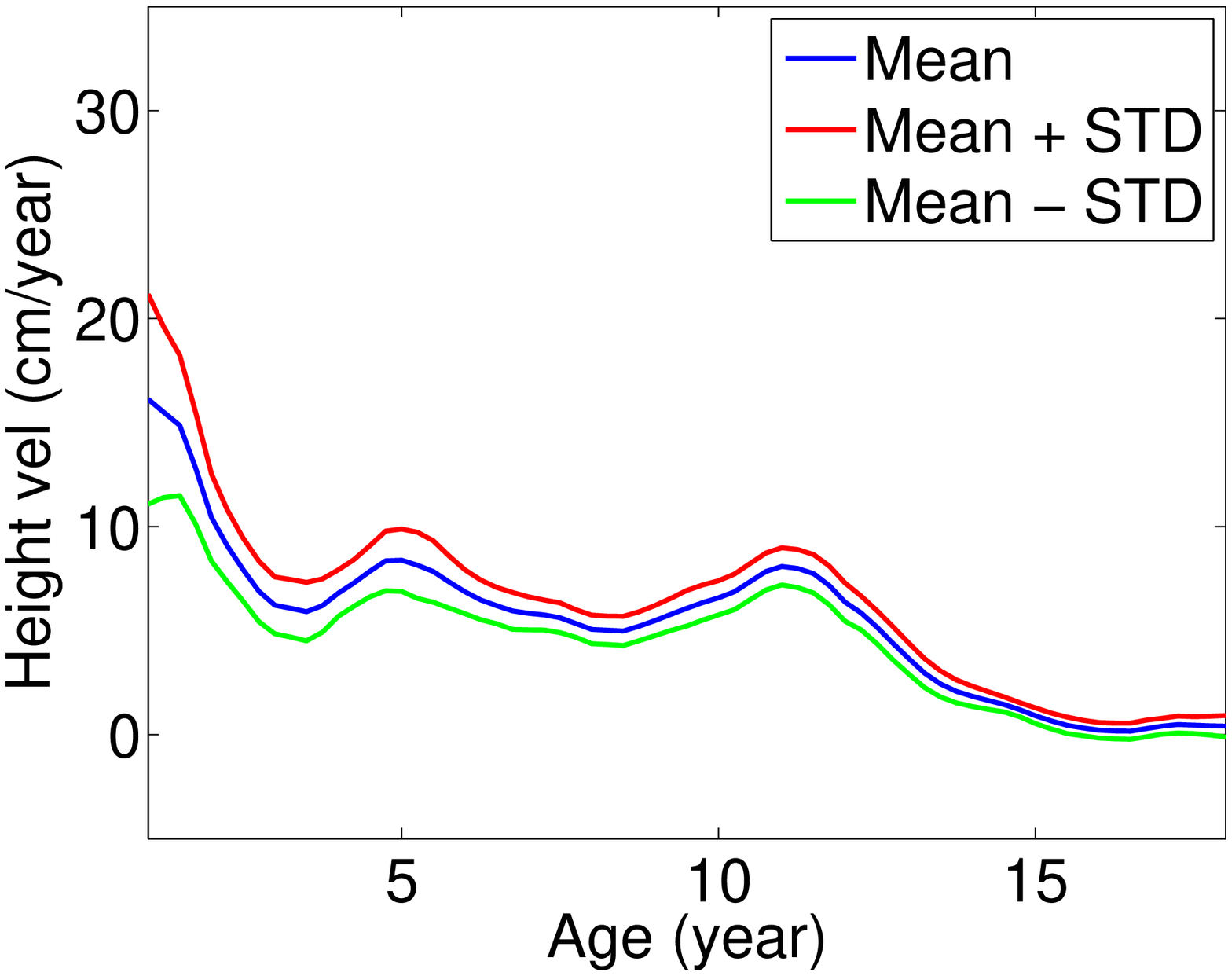} \ &
\includegraphics[height=1.0in]{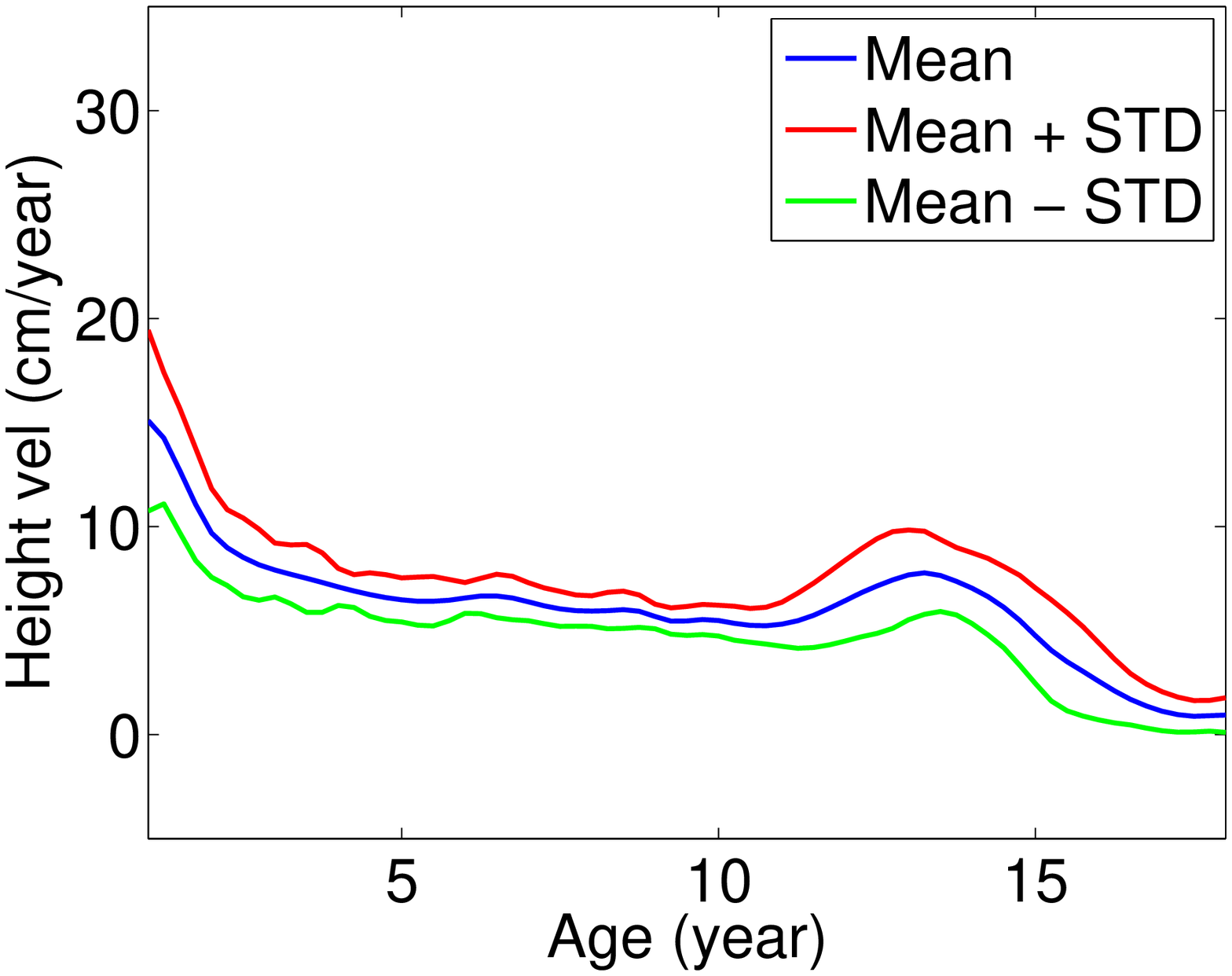} \ &
\includegraphics[height=1.0in]{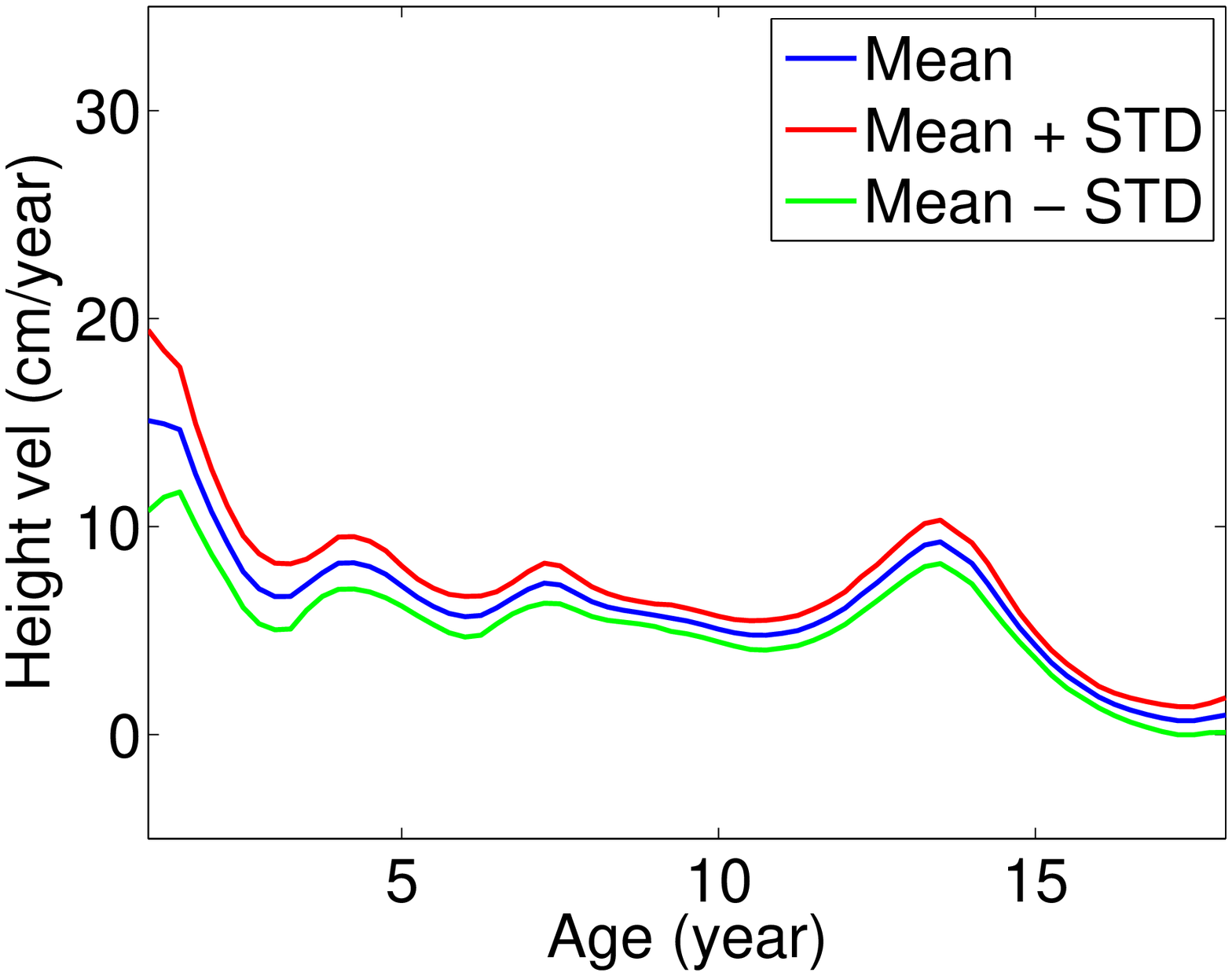} \\
\end{tabular}
\caption{Analysis of growth data. Top:  original data and the aligned functions.  Bottom: the corresponding
cross-sectional mean and mean $\pm$ standard deviations. } \label{fig:growth-results}
\end{center}
\end{figure}

The results from our elastic function analysis on the female growth
curves are shown in Fig. \ref{fig:growth-results} (left side). The
top-left panel shows the original data. It can
be seen from this plot that while the growth spurts for different
individuals occurs at slightly different times, there are some
underlying patterns to be discovered.  This can also be observed in the
cross-sectional mean and mean $\pm$ standard deviation plot in the
bottom-left panel. In the second panel of the top row we show the
aligned functions $\tilde{f}_i(t)$.
The panel below it, which shows the cross-sectional mean and mean $\pm$
standard deviation, exhibits a much tighter alignment of the functions
and, in turn, an enhancement of peaks and valleys in the aligned
mean. In fact, this mean function suggests the presence of
two growth spurts, one between 3 and 4 years, and the other between
10 and 12 years, on average.
Similar analysis is performed on the male growth curves and
Fig. \ref{fig:growth-results} (right) shows the results:  the original data (consisting of
39 derivatives of the original growth functions), the aligned
functions $\tilde{f}_i(t)$, and the corresponding cross-sectional means and means $\pm$ standard
deviations. The cross-sectional mean functions also show a much tighter
alignment of the functions and, in turn, an enhancement of peaks and
valleys in the aligned mean. This mean function
suggests the presence of several growth spurts,  between:  3 and 5, 6
and 7, and 13 and 14 years, on average.\\

\noindent 2. {\bf Handwriting Signature Data}: As another example of the 
data that can be effectively modeled using elastic functions,
we consider some handwritten signatures and the acceleration
functions along the signature curves. This application was also
considered in the paper \cite{kneip-ramsay:2008}. Let $(x(t), y(t))$
denote the $x$ and $y$ coordinates of a signature traced as a
function of time $t$. We study the acceleration
functions $f(t) = \sqrt{\ddot{x}(t)^2 + \ddot{y}(t)^2}$ for different
instances of the signatures and study their variability after
alignment.

\begin{figure}
\begin{center}
\begin{tabular}{cccc}
original data & mean $\pm$ std, before & aligned functions  & mean $\pm$ std, after \\
\includegraphics[height=1.0in]{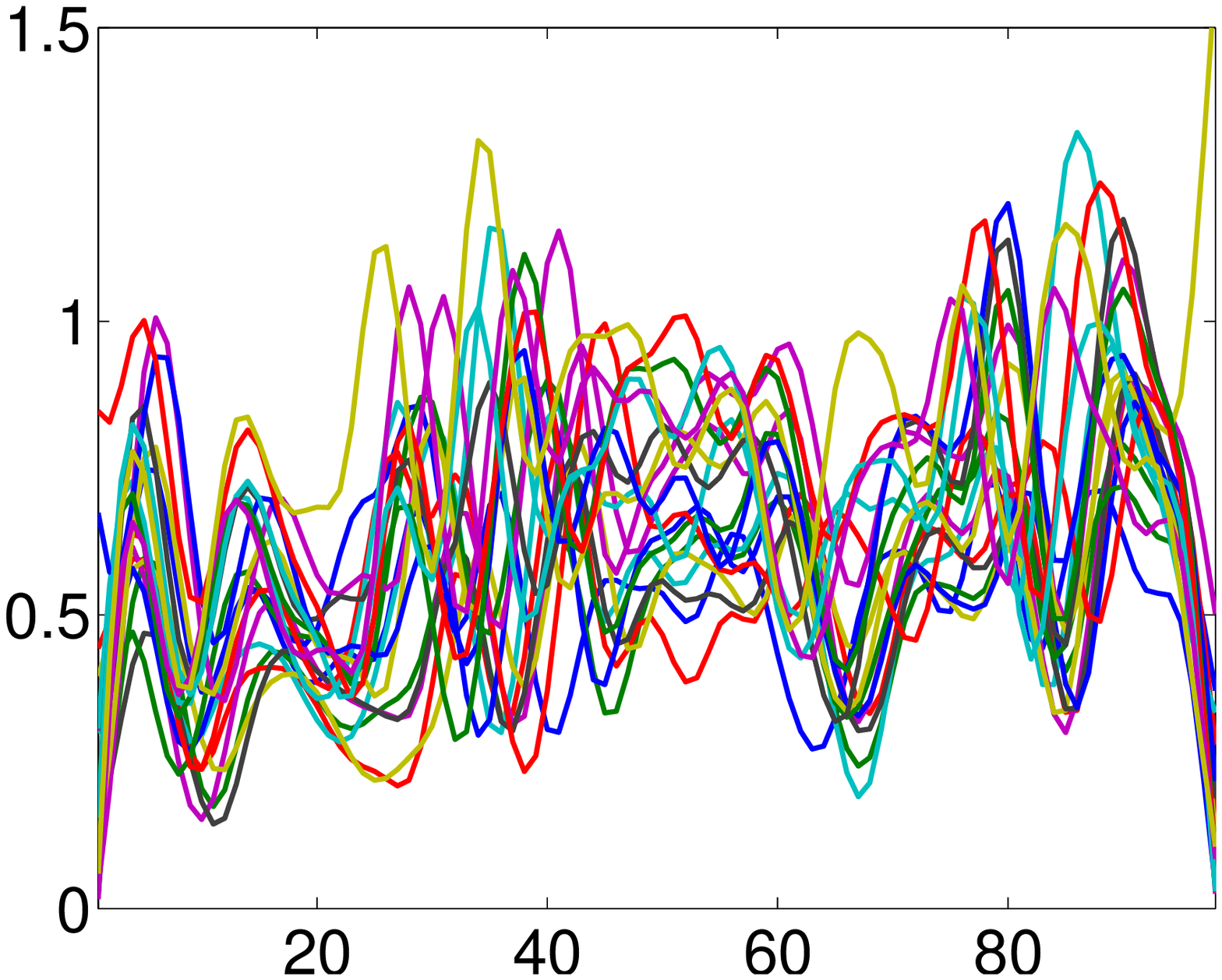} \ &
\includegraphics[height=1.0in]{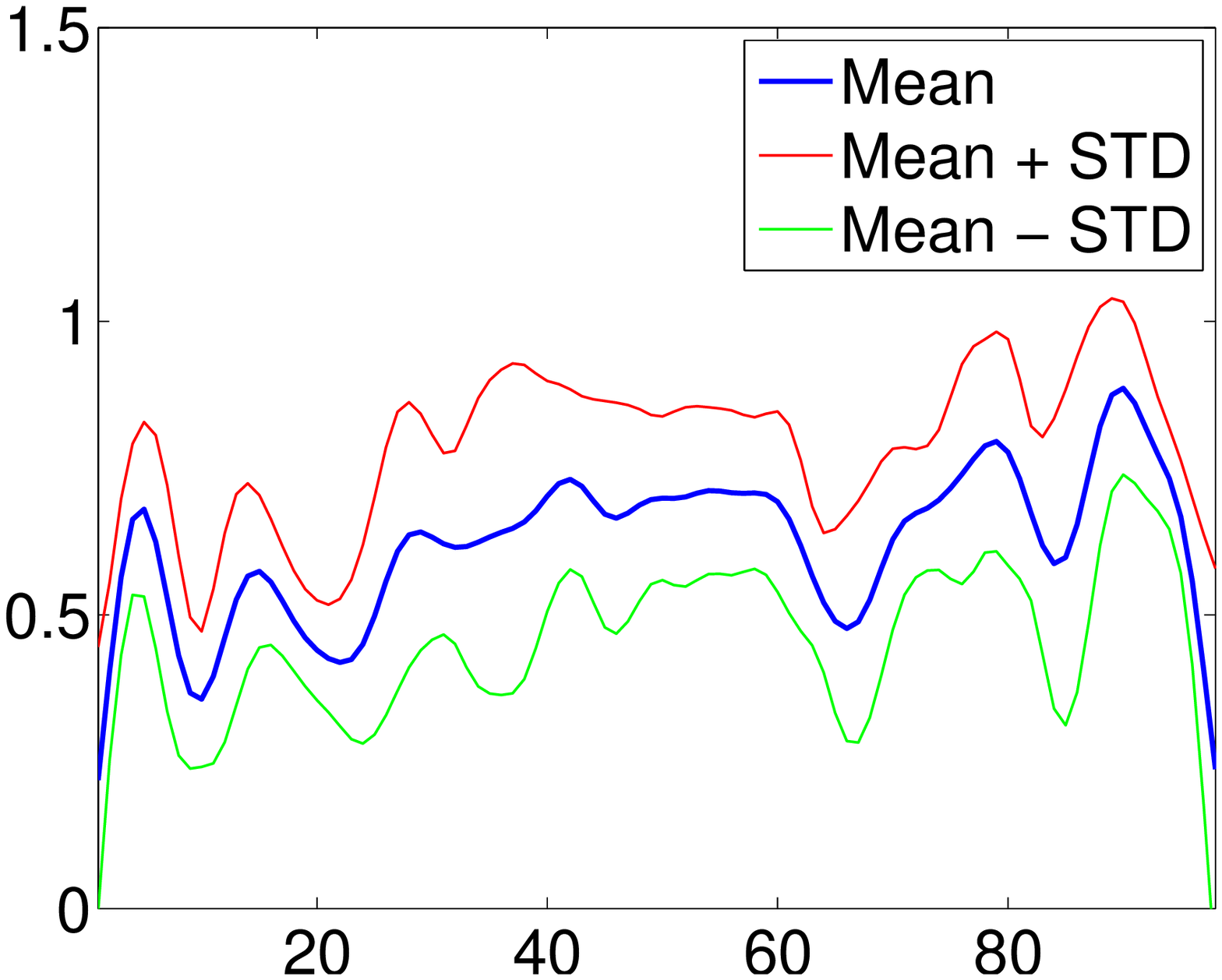} \ &
\includegraphics[height=1.0in]{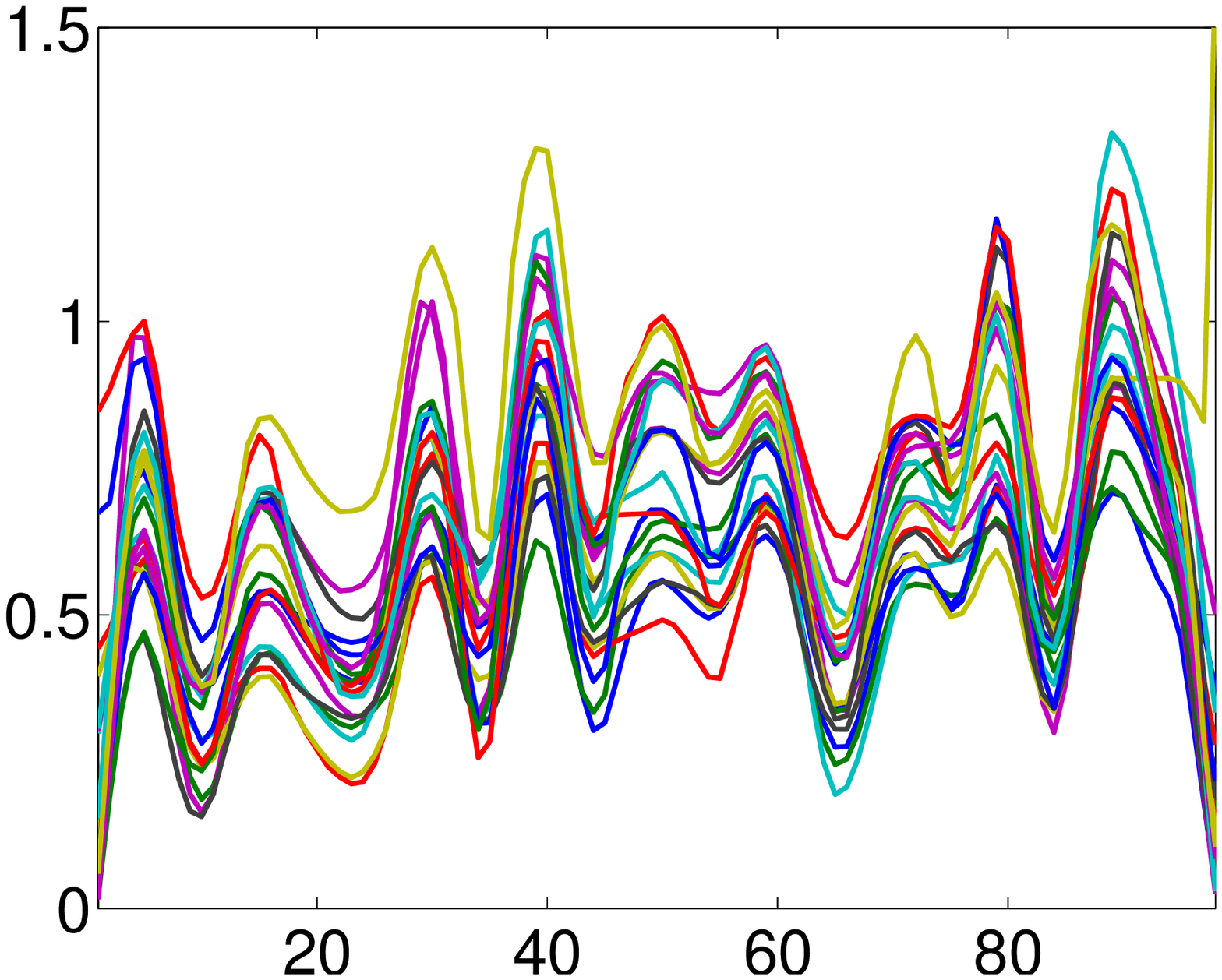} \ &
\includegraphics[height=1.0in]{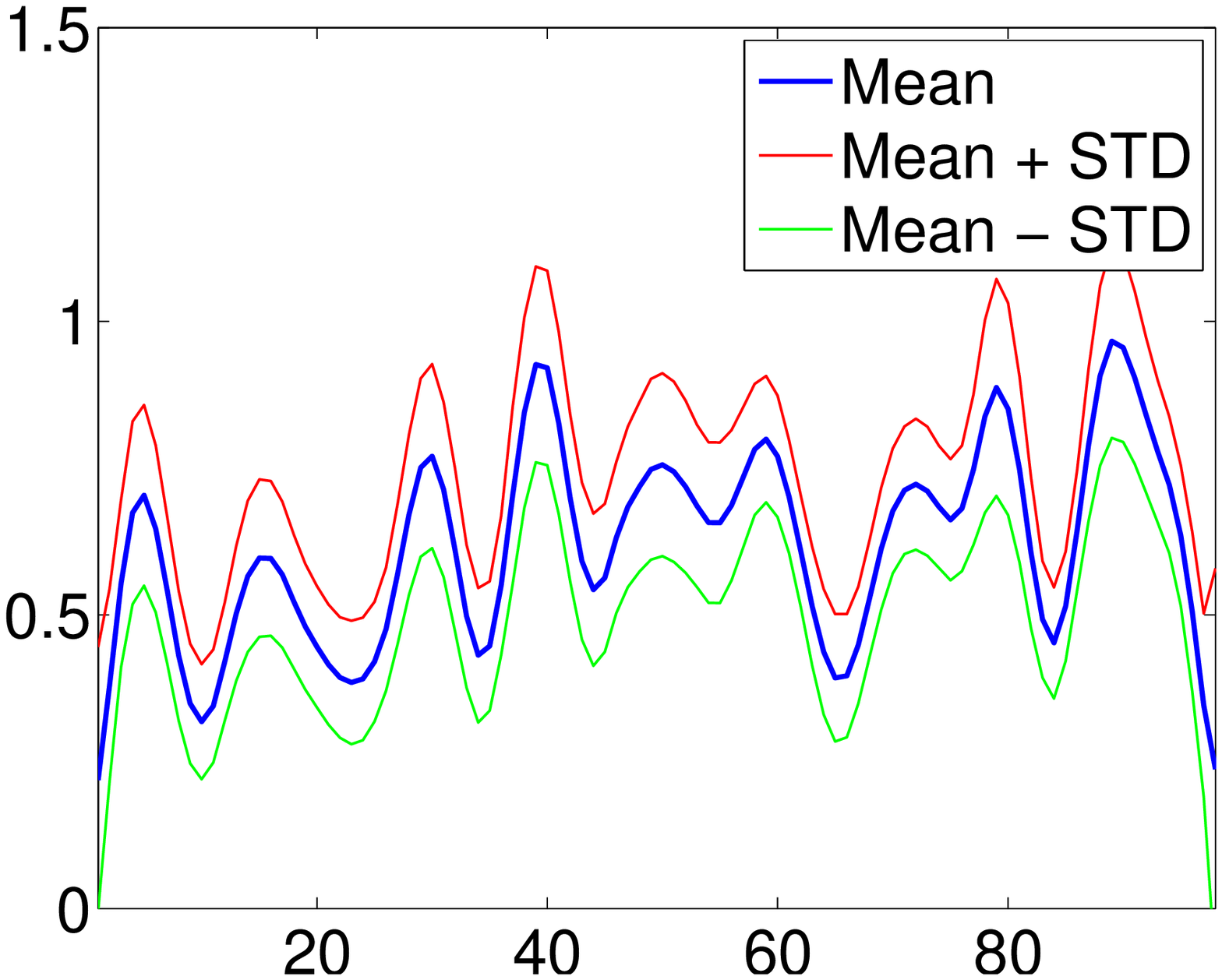}
\end{tabular}
\caption{Analysis of signature profiles.}
\label{fig:signature-results}
\end{center}
\end{figure}

The left panel in Fig. \ref{fig:signature-results} shows the 20
acceleration functions of 20 signatures that are used in our
analysis as $\{f_i\}$.  The corresponding cross-sectional mean and
mean $\pm$ standard deviations before the alignment are shown in the
next panel.  The right two panels show the aligned functions
$\tilde{f}_i$s and the corresponding mean and mean $\pm$ standard
deviations after the alignment. A look at the cross-sectional mean
functions suggests that the aligned functions have much more
exaggerated peaks and valleys, resulting from the alignment of these
features due to warping. \\

\begin{figure}
\begin{center}
\begin{tabular}{cccc}
original data & mean $\pm$ std, before & aligned functions  & mean $\pm$ std, after \\
\includegraphics[height=1.0in]{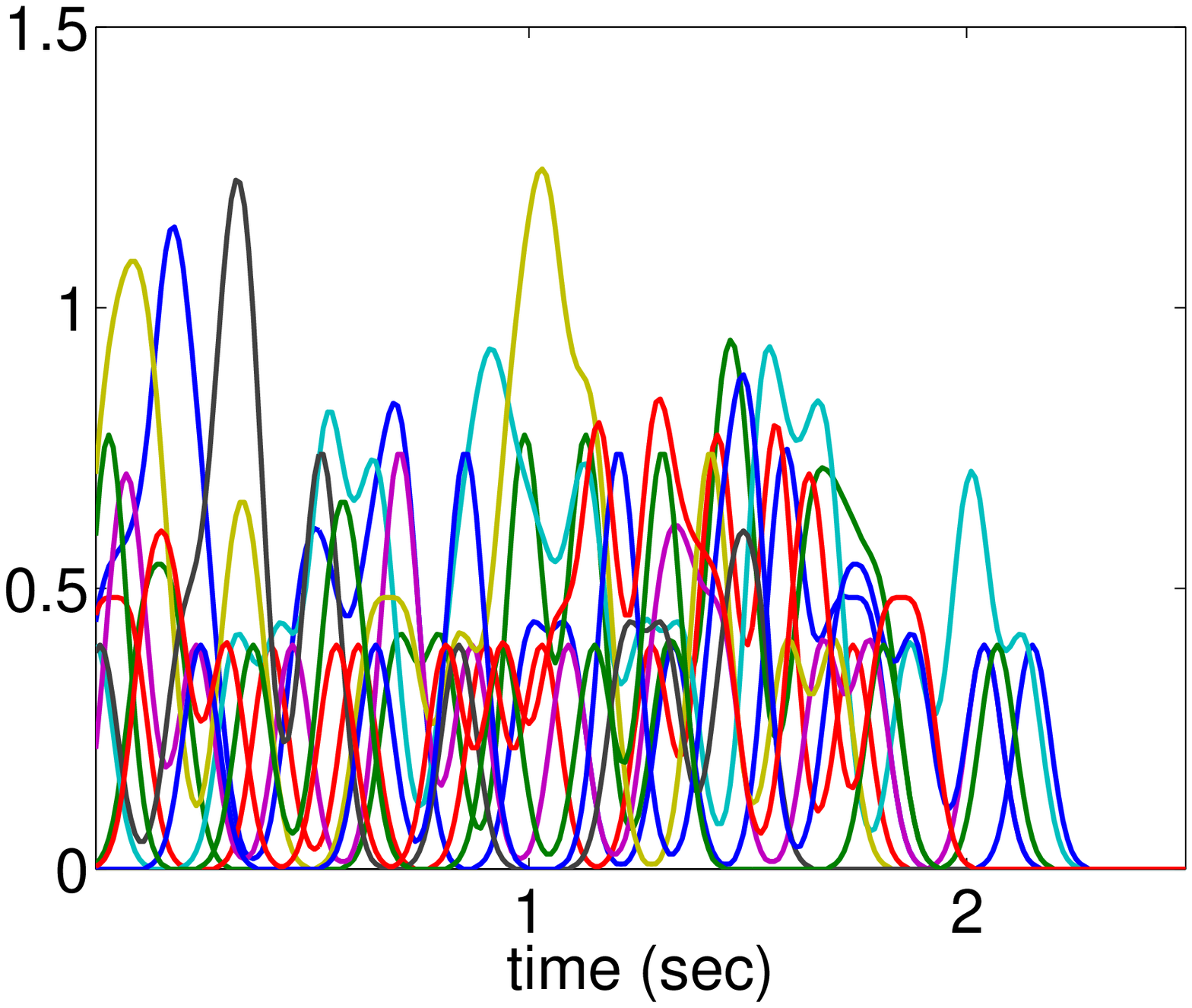} \ &
\includegraphics[height=1.0in]{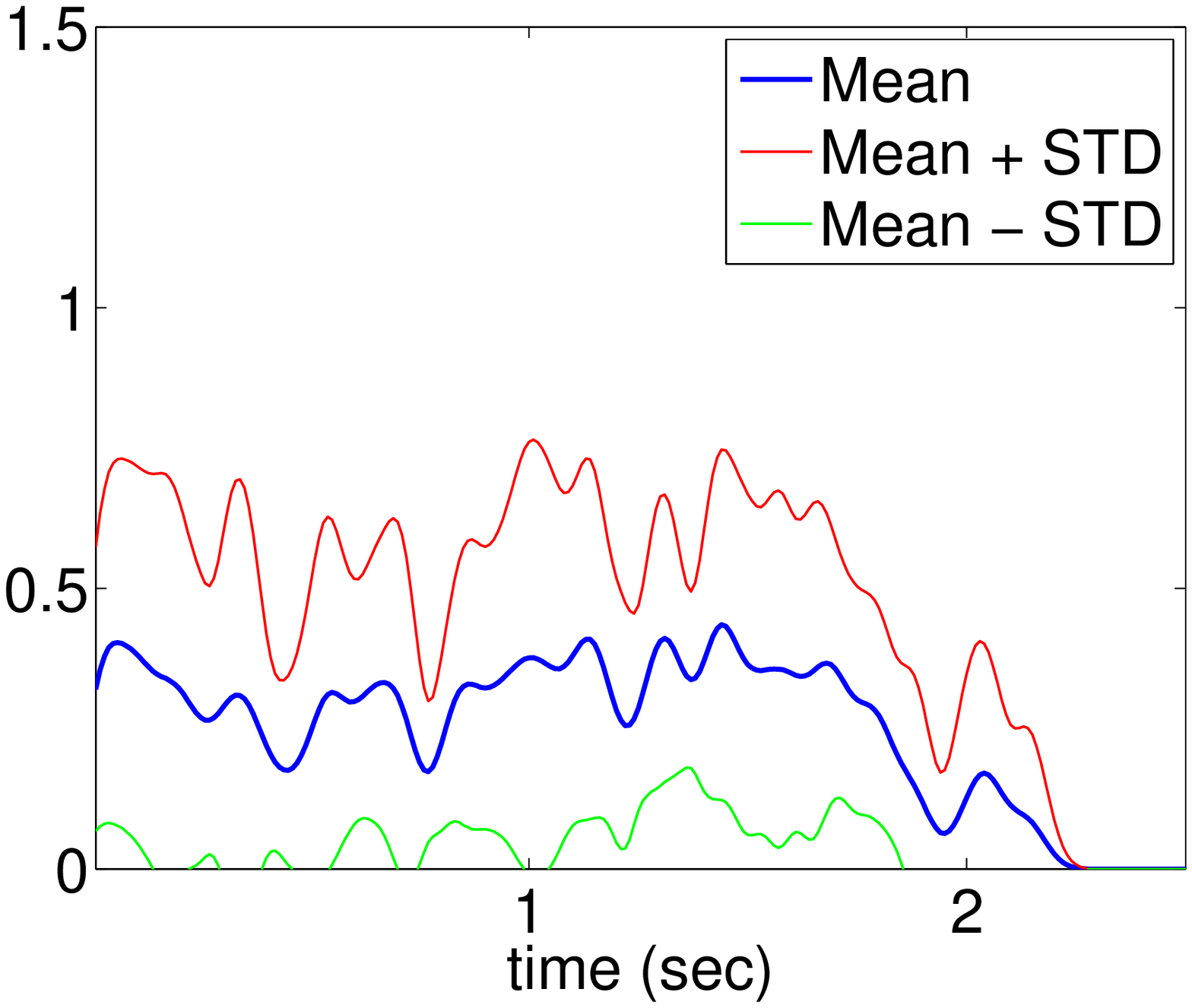} \ &
\includegraphics[height=1.0in]{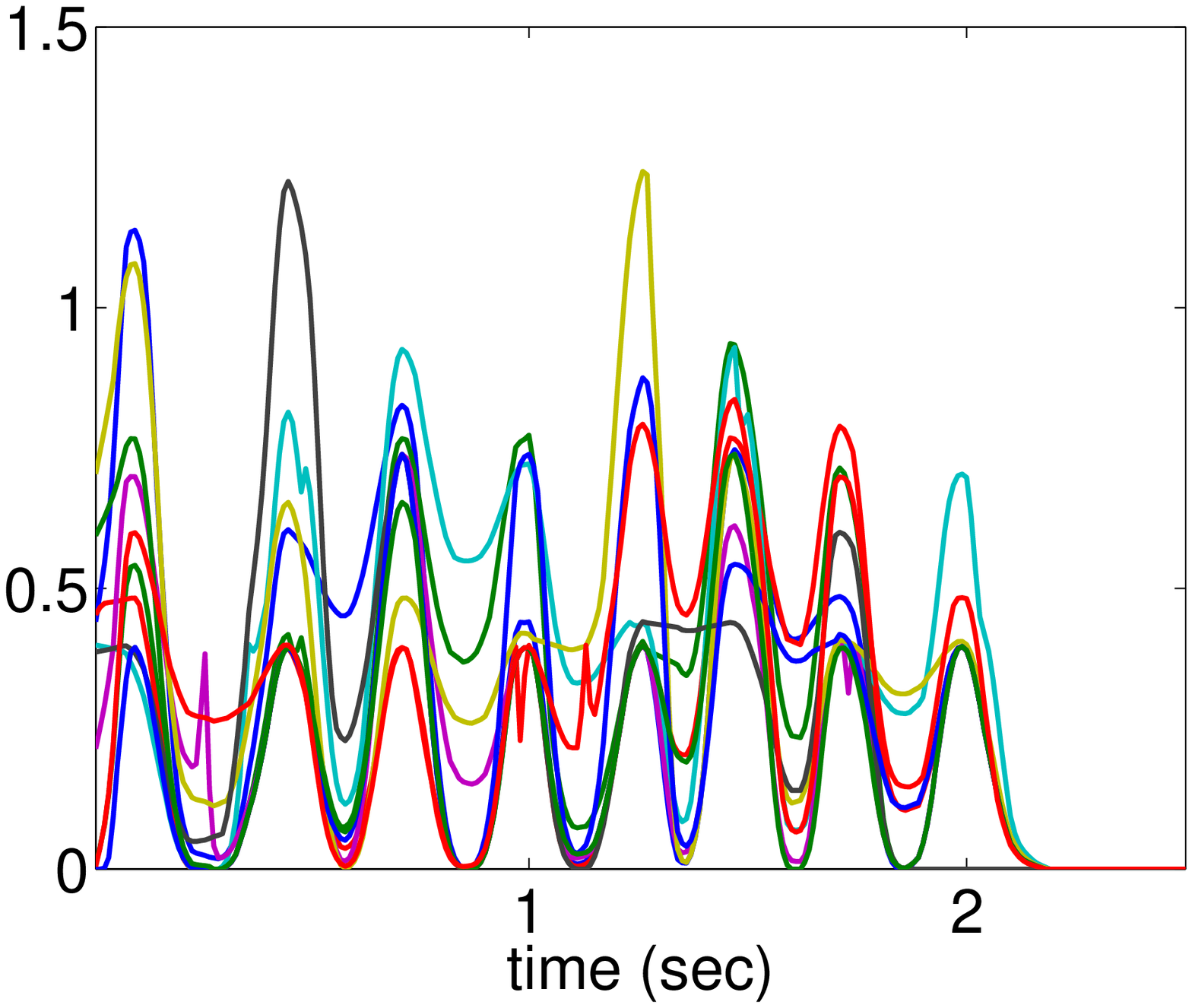} \ &
\includegraphics[height=1.0in]{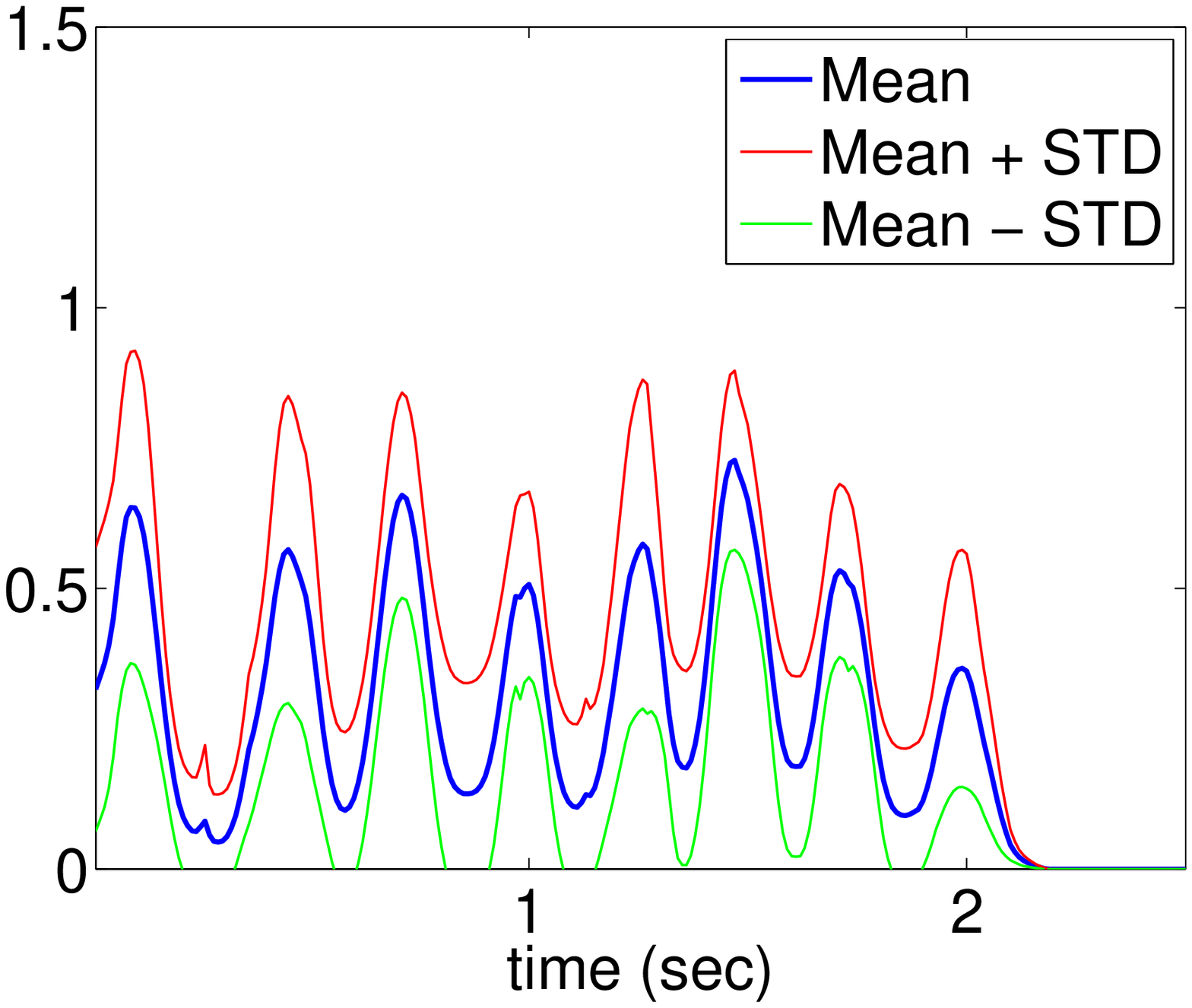}
\end{tabular}
\caption{Analysis of spike train data.} \label{fig:spike-results}
\end{center}
\end{figure}

\noindent 3. {\bf Neuroscience Spike Data}: Time-dependent information is
represented via sequences of stereotyped spike waveforms in the
nervous system.  These waveform sequences (or spike trains) have
been commonly looked as the {\it language} of the brain and are the
focus of much investigation.
Before we apply our framework, we need to convert the spike
information into functional data. Assume $s(t)$ is a spike train
with spike times $0 < t_1 < t_2 < \cdots < t_M < T$, where $[0, 1]$
denotes the recording time domain.  That is, $s(t) = \sum_{i=1}^M
\delta(t-t_i),\ \ t \in [0,1]\ ,$ where $\delta(\cdot)$ is the Dirac
delta function. One typically smooths the spike trains
to better capture the time correlation between spikes. In this paper
we use a Gaussian kernel $K(t) =
e^{-t^2/(2\sigma^2)}/(\sqrt{2\pi}\sigma)$, $\sigma \ge 0$ ($\sigma =
1ms$ here). That is, the smoothed spike train is $f(t) = (s*K)(t) =
\sum_{i=1}^M \frac{1}{\sqrt{2\pi}\sigma}
e^{-(t-t_i)^2/(2\sigma^2)}.$

Figure \ref{fig:spike-results} left panel shows one example of such
smoothed spike trains for 10 trials of one neurons in the primary motor
cortex of a Macaque monkey subject that was performing a squared-path
movement \cite{wu-srivastava-neuro:10}.  The next panel shows the
cross-sectional mean and mean $\pm$ standard deviation of the
functions in this neuron.  The third panel shows
$\{\tilde{f}_i\}$ where we see that the functions are well aligned with
more exaggerated peaks and valleys.  The next panel shows the mean
and mean $\pm$ standard deviation.  Similar to the growth data
and signature data, an increased amplitude variation and decreased
standard deviation are observed in this plot. \\

\begin{figure}
\begin{center}
\begin{tabular}{cccc}
original data & mean $\pm$ std, before & aligned functions  & mean $\pm$ std, after \\
\includegraphics[height=1.0in]{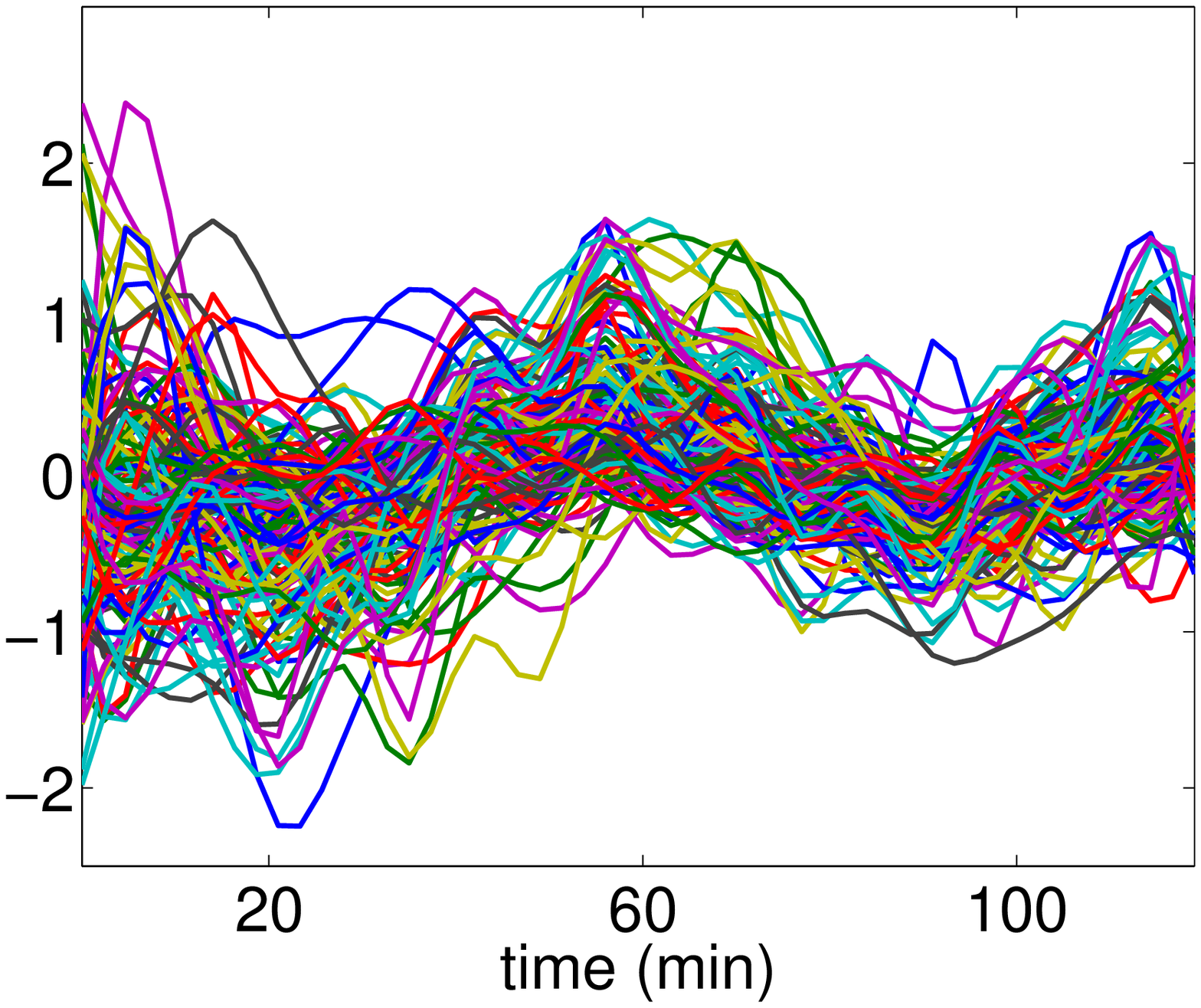} \ &
\includegraphics[height=1.0in]{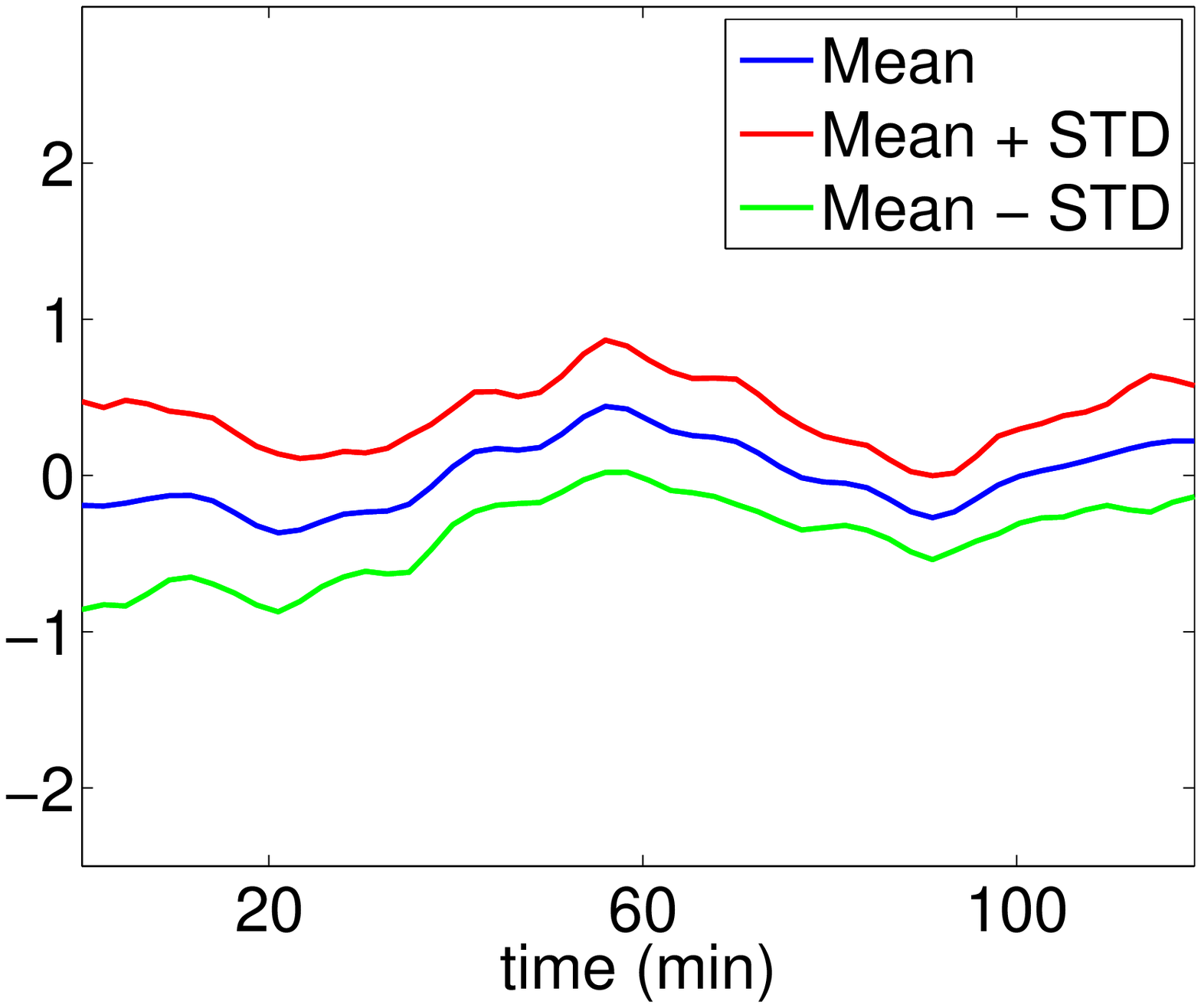} \ &
\includegraphics[height=1.0in]{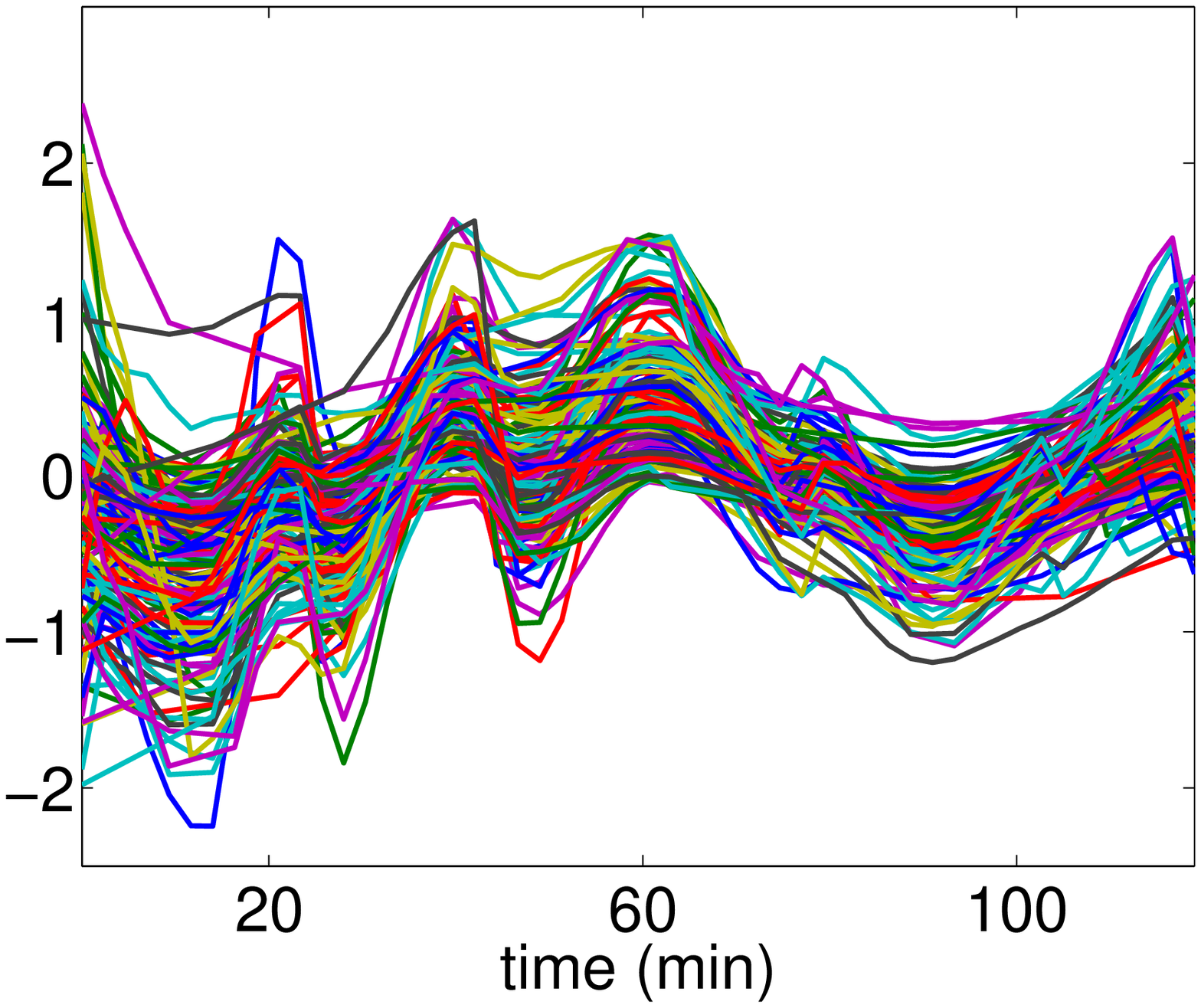} \ &
\includegraphics[height=1.0in]{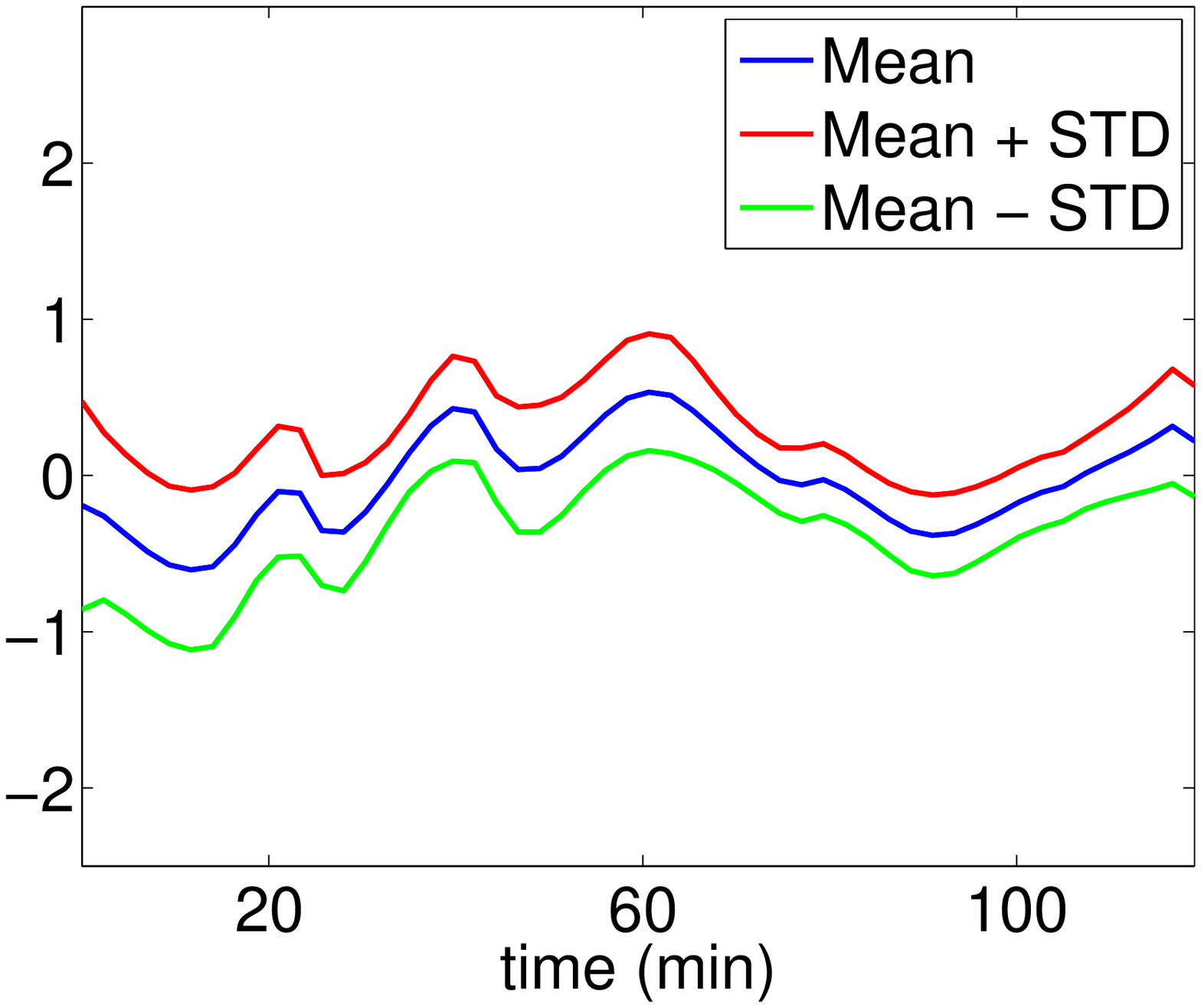}
\end{tabular}
\caption{Analysis of gene expression data.} \label{fig:gene-results}
\end{center}
\end{figure}

\noindent 4. {\bf Gene Expression data}: In this example, we consider temporal
microarray gene expression profiles.  The time ordering for yeast
cell-cycle genes was investigated in
\cite{muller-biostatistics:2006}, and we use the same data in this
study.  The expression level was measured during a period of 119
minutes for a total of 612 fully-recorded genes.  There are 5
clusters with respect to phases in these continuous function.  In
particular, 159 of these functions were known to be related to M
phase regulation of the yeast cell cycle.
These 159 functions used here are same as those used in \cite{muller-biostatistics:2006}, 
and are
shown in the left panel in Fig. \ref{fig:gene-results}. Although, in general, gene expression 
analysis has many goals and problems, we use focus only on the subproblem of
expression alignment as functional data. 
The corresponding cross-sectional mean and mean $\pm$ standard
deviations before the alignment are shown in the next panel.  The
right two panels show the aligned functions $\tilde{f}_i$s and the
corresponding mean and mean $\pm$ standard deviations after the
alignment. Once again we find a strong
alignment of functional data with improved peaks and valleys.

\subsection{Comparisons with other Methods}
In this section we compare the results from our method to some of the past
ideas where the software is available publicly. While we have compared our
framework with other published work in conceptual terms earlier, in this section
we focus on a purely empirical evaluation. In particular, we utilize several
evaluation criteria for comparing the alignments of functional data in
the several simulated and real datasets discussed in previous sections.
The choice of an evaluation criteria is not obvious, as there is no single criterion
that has been used consistently by past authors for measuring the
quality of alignment. Thus,  we use three criteria so that together they provide
a more comprehensive evaluation.
We will continue to use $f_i$
and $\tilde{f}_i$, $i = 1, ..., N$, to denote the original and the
aligned functions, respectively.

\begin{enumerate}
\item {\bf Least Squares}: A cross-validated measure of the level of synchronization \cite{james:10}:
\begin{equation}
ls = \frac{1}{N}\sum_{i=1}^N {\int (\tilde{f}_i(t) -
\frac{1}{N-1}\sum_{j\neq i} {\tilde{f}_j}(t))^2 dt  \over \int
({f}_i(t) -  \frac{1}{N-1}\sum_{j\neq i}f_j(t))^2 dt } \ ,
\end{equation}
$ls$ measures the total cross-sectional variance of the aligned
functions, relative to the original value. The smaller the value of
$ls$, the better the alignment is.
\item {\bf Pairwise Correlation}: It measures pairwise correlation between functions:
\begin{equation}
pc = {\sum_{i\neq j} cc(\tilde{f}_i(t), \tilde{f}_j(t)) \over
\sum_{i\neq j} cc({f}_i(t), {f}_j(t))}\ ,
\end{equation}
where $cc(f,g)$
is the pairwise Pearson's correlation between functions.
The larger the value of $pc$, the better the alignment between
functions in general.

\item {\bf Sobolev Least Squares}: This time we compute the
least squares using the first derivative of the functions:
\begin{equation}
sls = {\sum_{i=1}^N \int (\dot{\tilde
f}_i(t)-\frac{1}{N}\sum_{j=1}^N\dot{\tilde{f_j}})^2 dt \over
\sum_{i=1}^N \int ({\dot f}_i(t)-\frac{1}{N}\sum_{j=1}^N\dot{f_j})^2
dt} \ ,
\end{equation}
This criterion measures the total cross-sectional variance of the derivatives  of the aligned
functions, relative to the original value,
and is an alternative measure of the synchronization.  The smaller
the value of $sls$, the better synchronization the method achieves.
\end{enumerate}

\begin{figure}
\begin{center}
\begin{tabular}{|c|c|c|c|c|c|}
\hline
Original & AUTC  \cite{muller-JASA:2004} & PACE
\cite{muller-biometrika:2008}
   & SMR \cite{gervini-gasser-RSSB:04} & MBM \cite{james:10} & F-R \\
\hline
\includegraphics[height=0.85in]{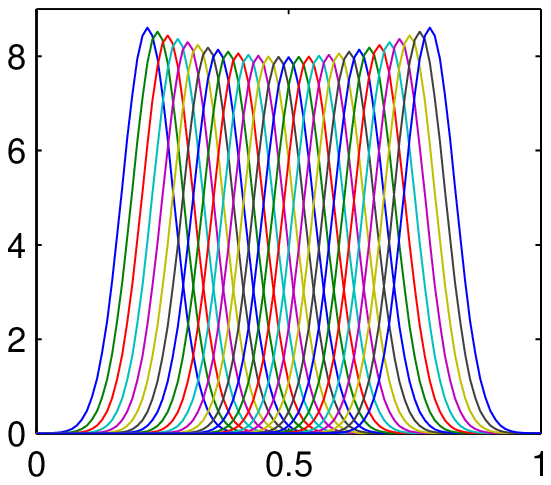} &
\includegraphics[height=0.85in]{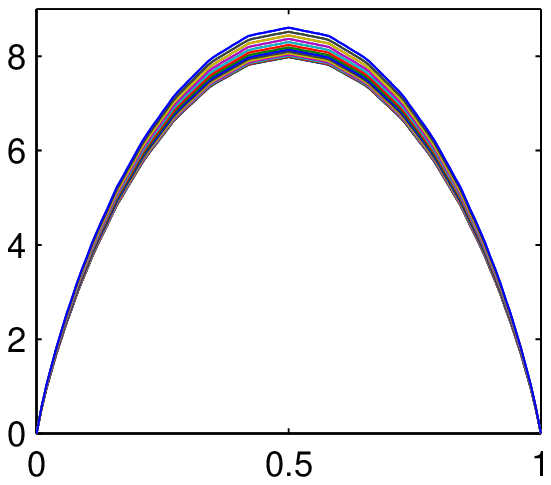} &
\includegraphics[height=0.85in]{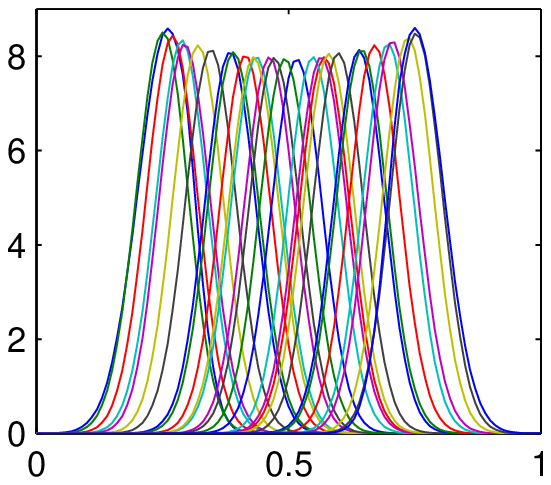} &
\includegraphics[height=0.85in]{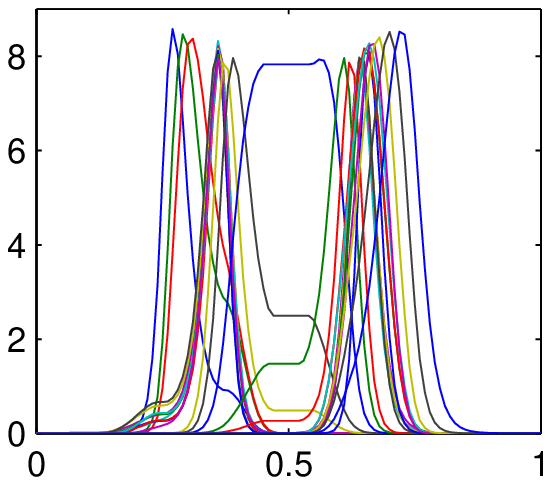} &
\includegraphics[height=0.85in]{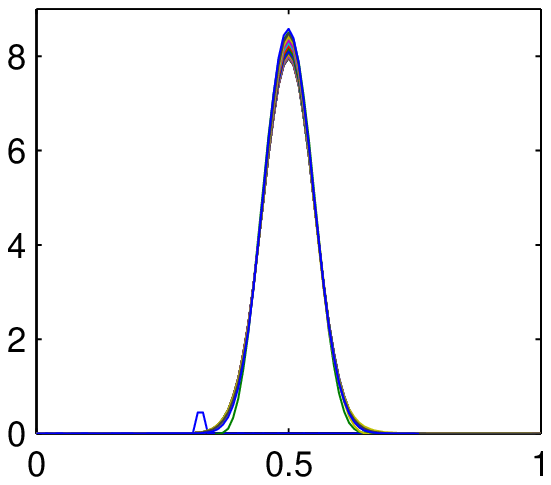} &
\includegraphics[height=0.85in]{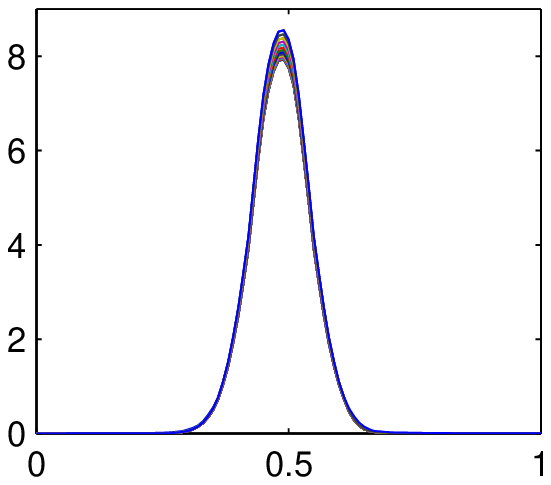}
\\
{\scriptsize Gaussian kernel} & {\scriptsize (0.01, 11.0, {\bf
0.00})} & {\scriptsize (0.98, 1.18, 0.98)} & {\scriptsize (0.58,
2.86, 1.18)} & {\scriptsize (0.03,
10.3, 0.04)} & {\scriptsize {\bf (0.00, 11.0, 0.00)}} \\
\hline
\includegraphics[height=0.85in]{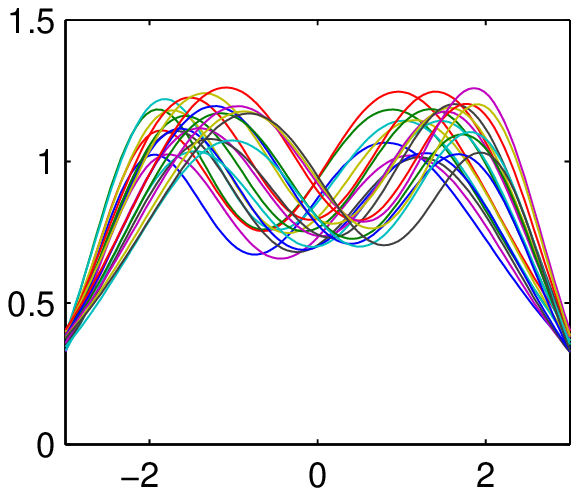} &
\includegraphics[height=0.85in]{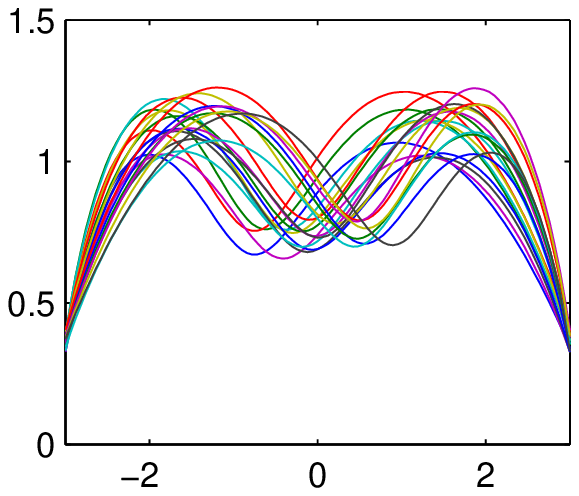} &
\includegraphics[height=0.85in]{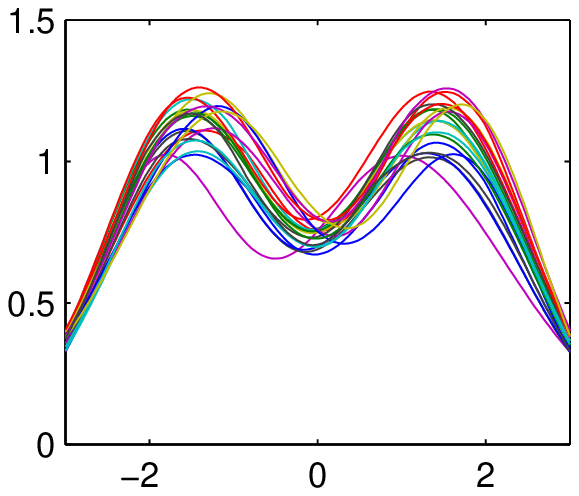} &
\includegraphics[height=0.85in]{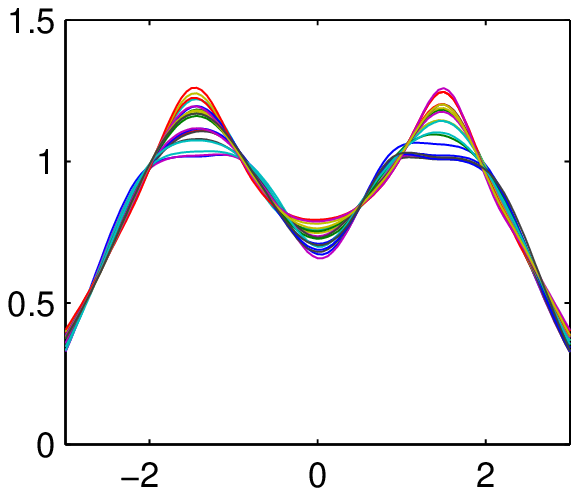} &
\includegraphics[height=0.85in]{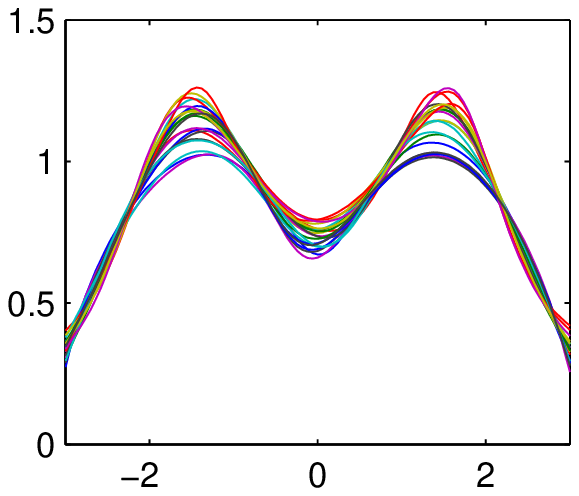} &
\includegraphics[height=0.85in]{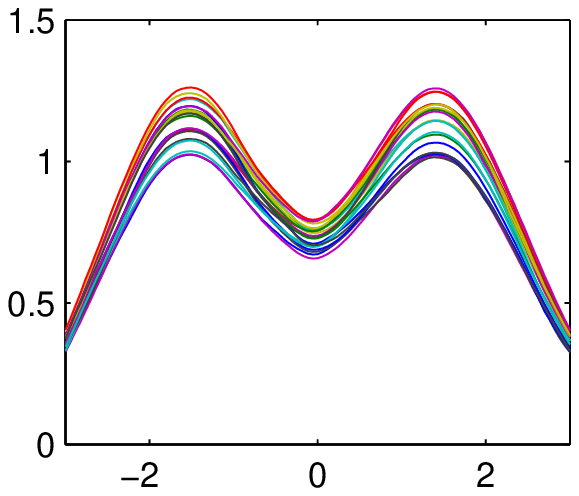}
\\
{\scriptsize Bimodal}& {\scriptsize (0.89, 1.01, 0.82)} &
{\scriptsize (0.49, 1.22, 0.23)} & {\scriptsize ({\bf 0.11,  1.27},
0.30)} & {\scriptsize (0.16, {\bf
1.27}, 0.29)} &{\scriptsize  (0.29, {\bf  1.27, 0.04})} \\
\hline
\includegraphics[height=0.85in]{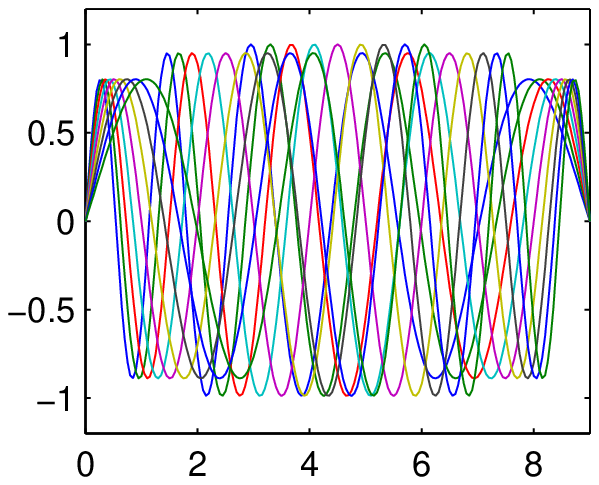} &
\includegraphics[height=0.85in]{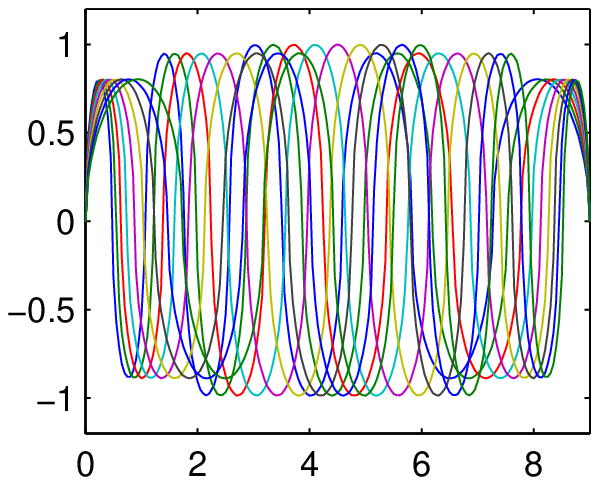} &
\includegraphics[height=0.85in]{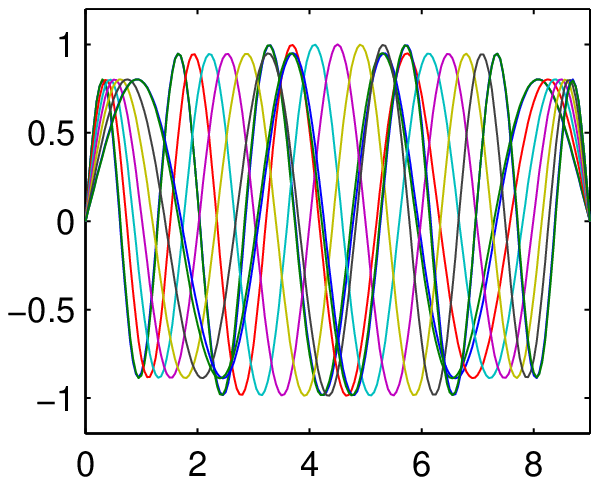} &
\includegraphics[height=0.85in]{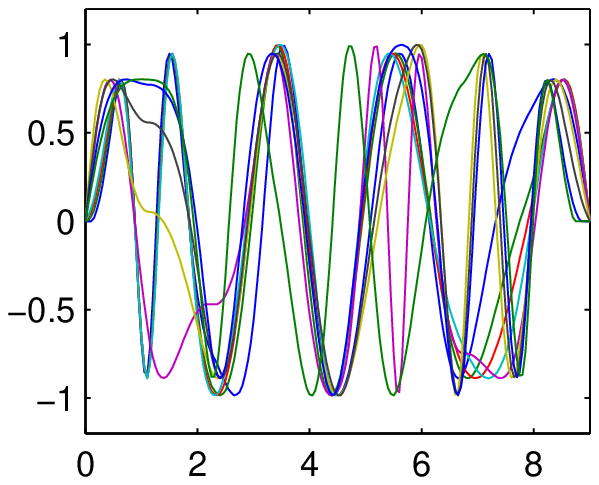} &
\includegraphics[height=0.85in]{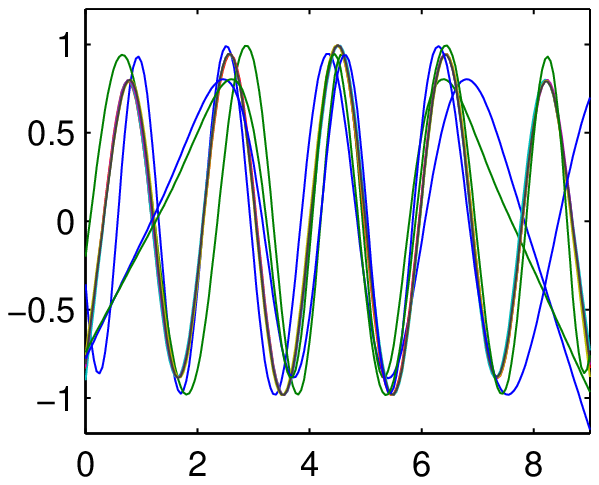} &
\includegraphics[height=0.85in]{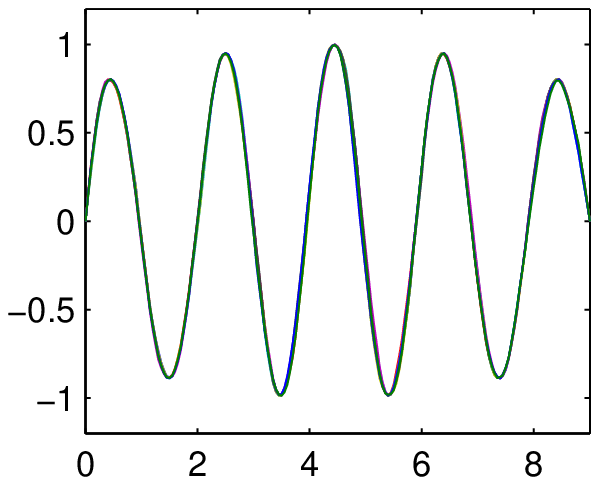}
\\
{\scriptsize Wave function}& {\scriptsize (1.33, 2.12, 1.73)} &
{\scriptsize (0.92, 13.0, 0.92)} & {\scriptsize (0.48, 83.7, 0.66)}
& {\scriptsize (0.30,
124, 0.24)} &{\scriptsize ({\bf 0.00, 175, 0.00})} \\
\hline
\includegraphics[height=0.85in]{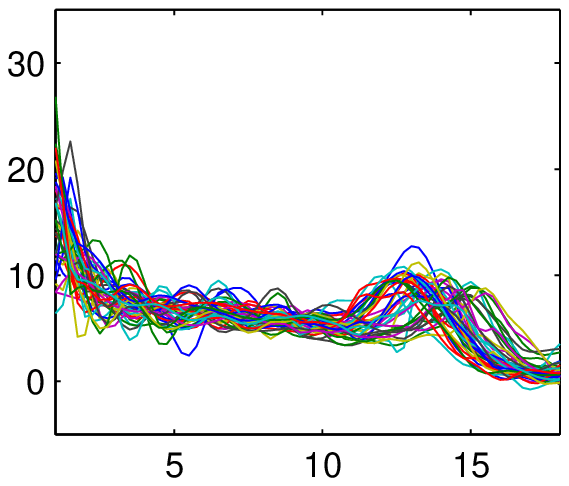} &
\includegraphics[height=0.85in]{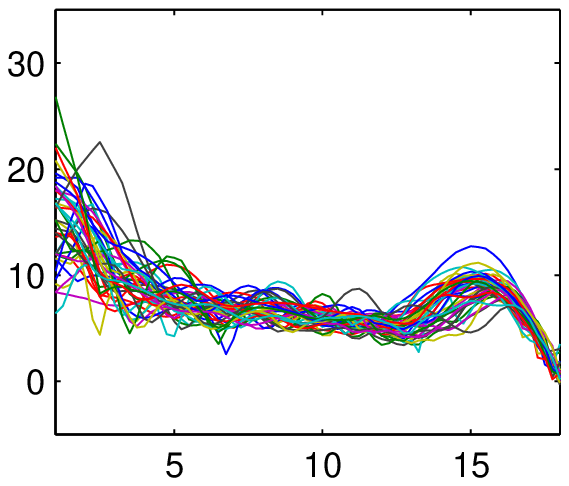} &
\includegraphics[height=0.85in]{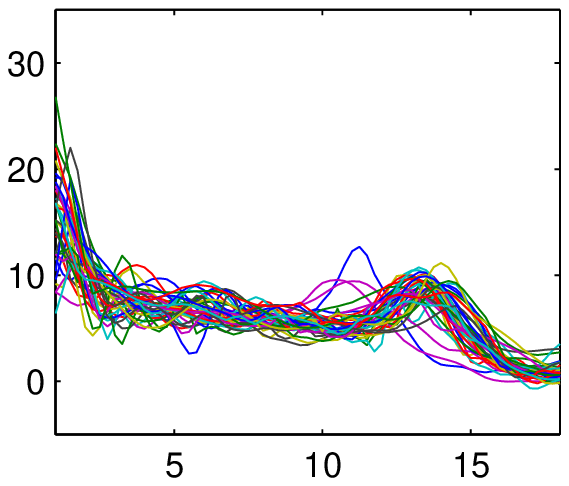} &
\includegraphics[height=0.85in]{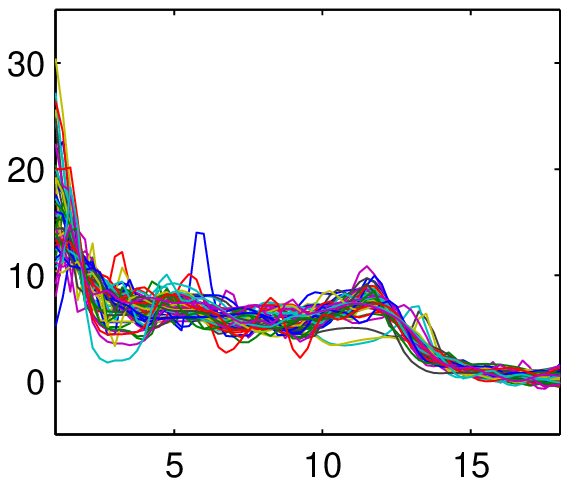} &
\includegraphics[height=0.85in]{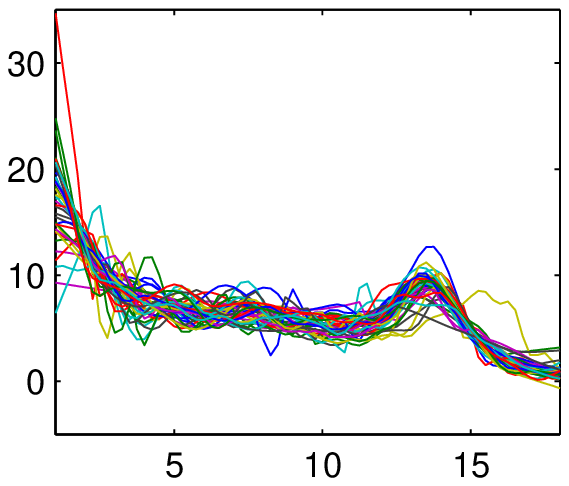} &
\includegraphics[height=0.85in]{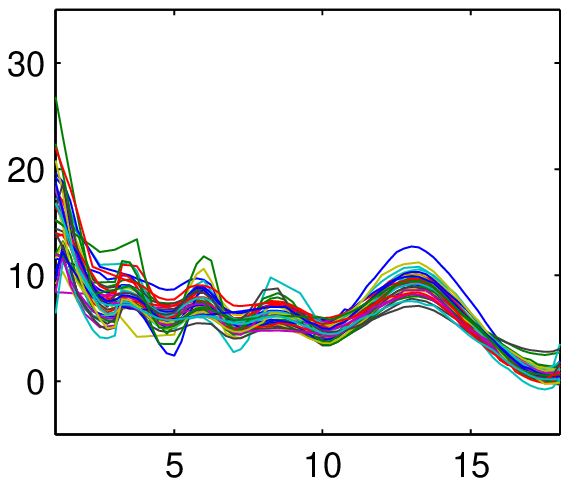}
\\
{\scriptsize Growth-male}& {\scriptsize (1.10, 1.05, 0.55)} &
{\scriptsize (0.91, 1.09, 0.68)} & {\scriptsize ({\bf 0.45}, 1.17,
0.77)} & {\scriptsize (0.70, 1.17,
0.62)} & {\scriptsize (0.64, {\bf 1.18, 0.31})} \\
\hline
\includegraphics[height=0.85in]{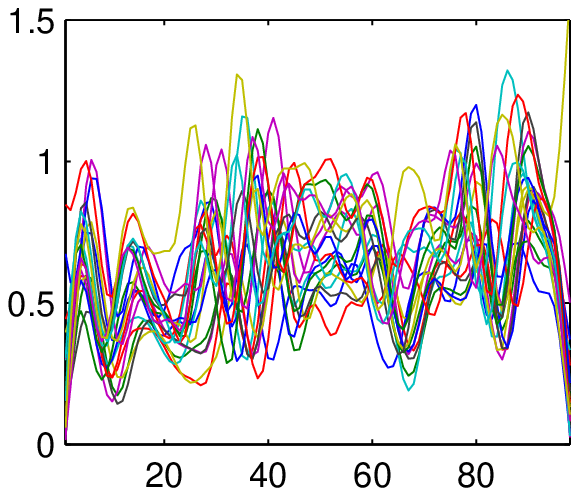} &
\includegraphics[height=0.85in]{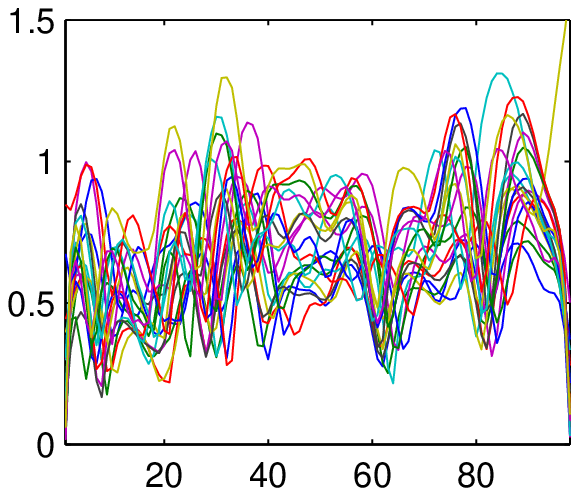} &
\includegraphics[height=0.85in]{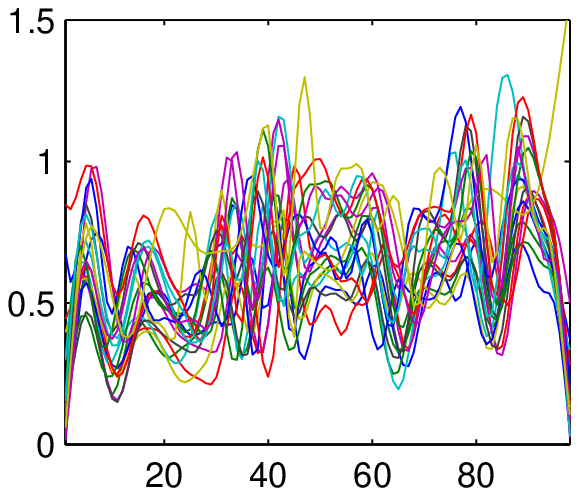} &
\includegraphics[height=0.85in]{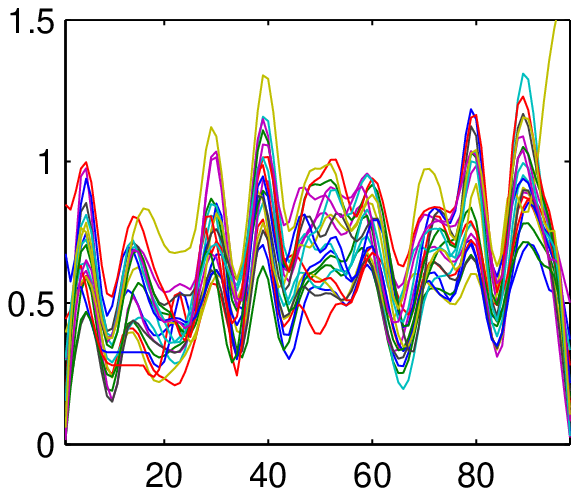} &
\includegraphics[height=0.85in]{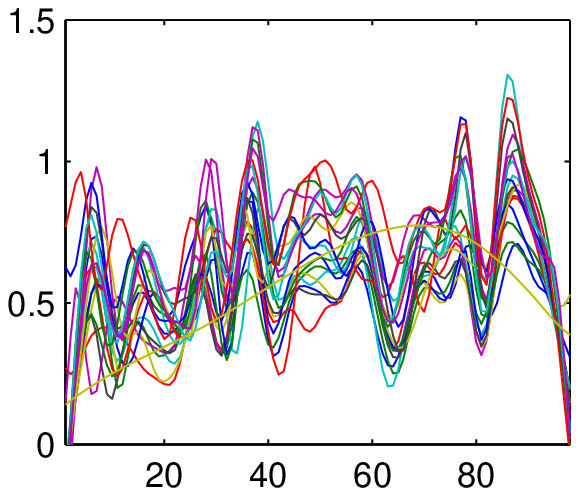} &
\includegraphics[height=0.85in]{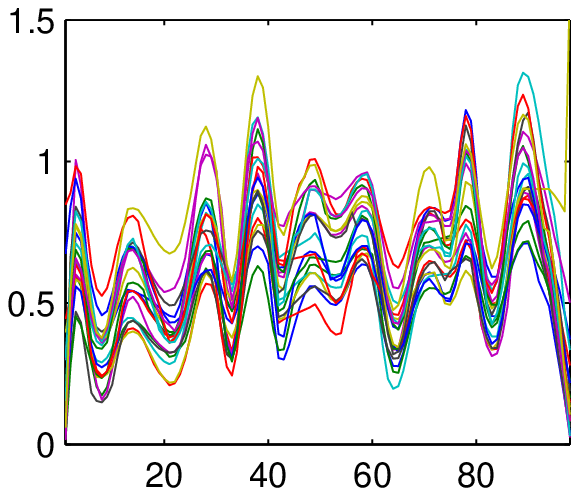}
\\
{\scriptsize Signature}& {\scriptsize (1.00, 1.02, 0.99)} &
{\scriptsize (0.91, 1.18, 0.84)} & {\scriptsize (0.62, 1.59, {\bf
0.31})} & {\scriptsize (0.64,
1.57, 0.46)} &{\scriptsize ({\bf 0.56, 1.79, 0.31})} \\
\hline
\includegraphics[height=0.85in]{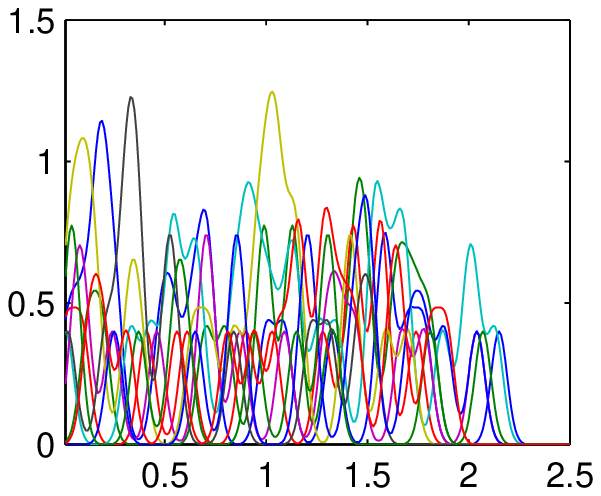} &
\includegraphics[height=0.85in]{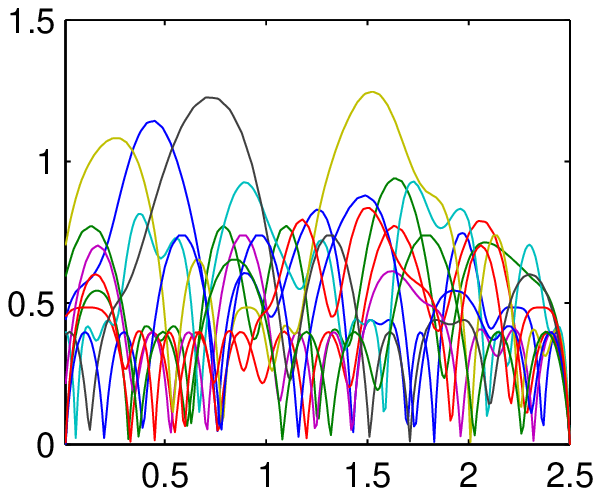} &
\includegraphics[height=0.85in]{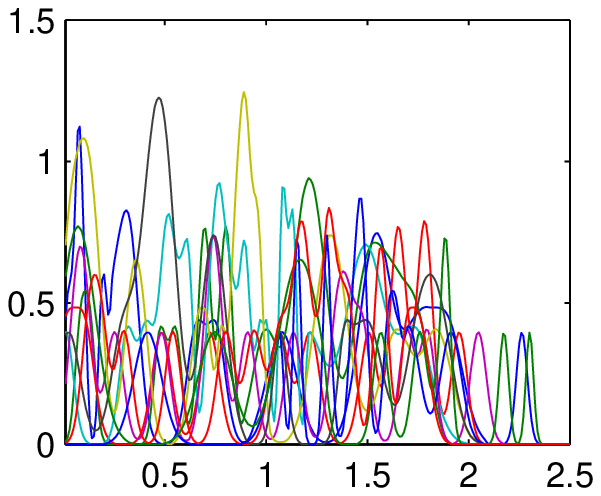} &
\includegraphics[height=0.85in]{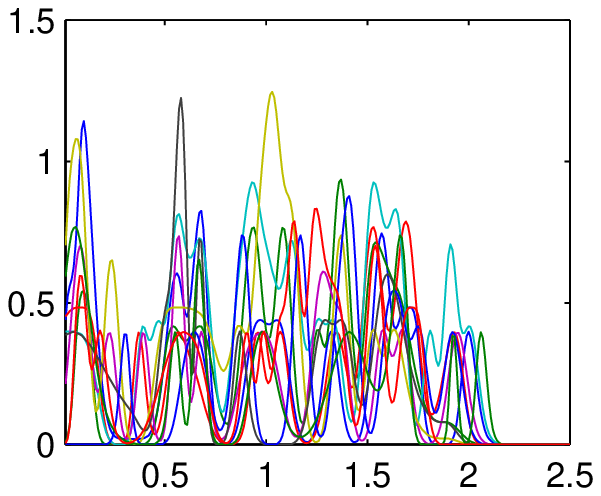} &
\includegraphics[height=0.85in]{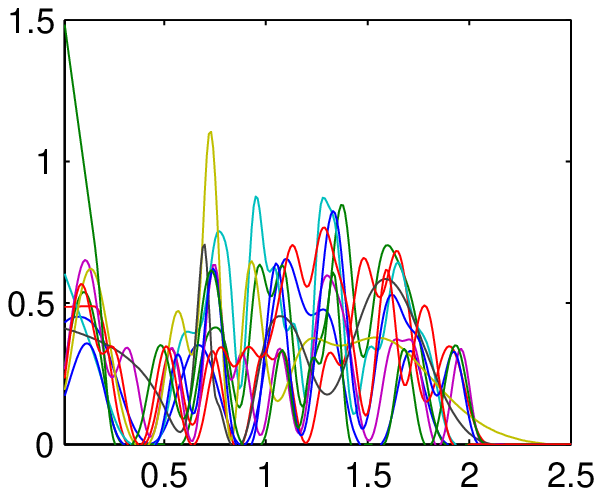} &
\includegraphics[height=0.8in]{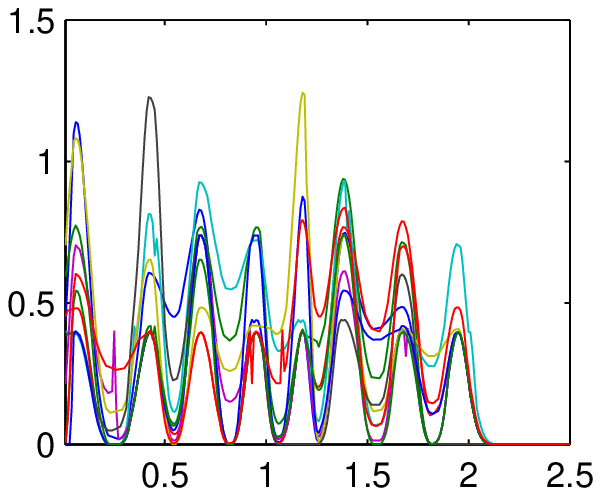}
\\
{\scriptsize Neural data}& {\scriptsize (1.05, 0.27, 0.70)}
&{\scriptsize (0.87, 1.35, 1.10)} & {\scriptsize (0.69, 2.54, 0.95)}
& {\scriptsize (0.48,
3.06, 0.40)} &{\scriptsize ({\bf 0.40, 3.77, 0.28})}\\
\hline
\includegraphics[height=0.85in]{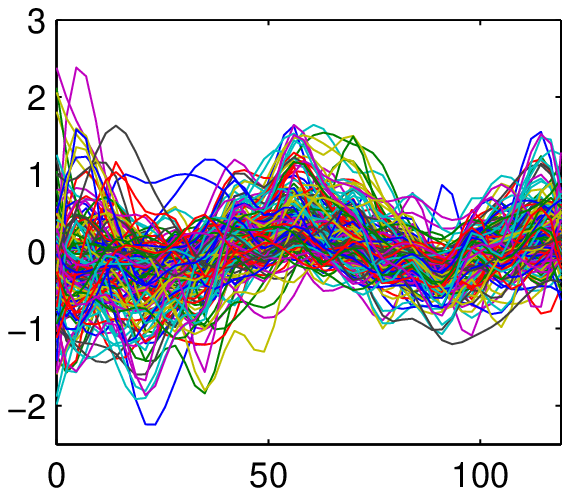} &
\includegraphics[height=0.85in]{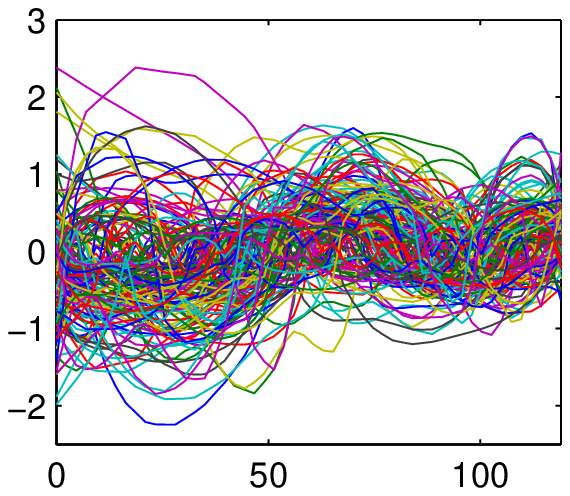} &
\includegraphics[height=0.85in]{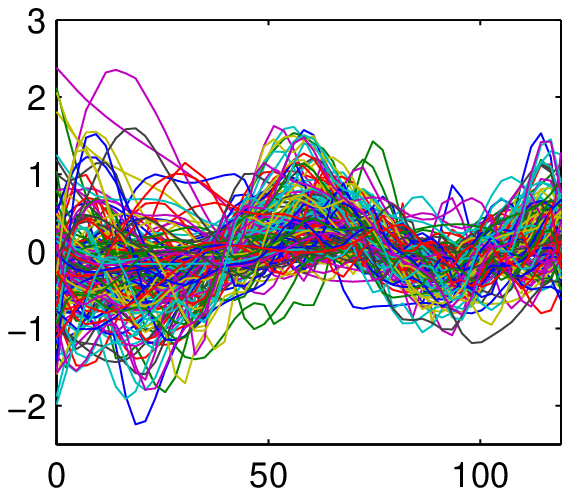} &
\includegraphics[height=0.85in]{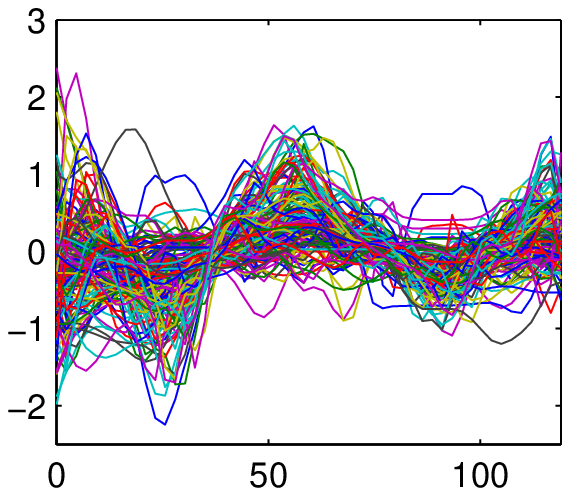} &
\includegraphics[height=0.85in]{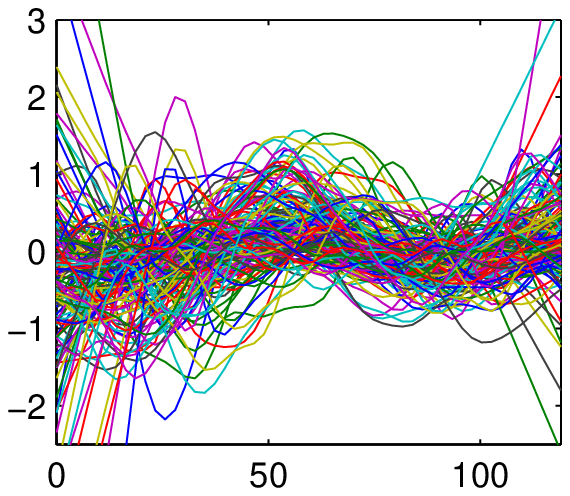} &
\includegraphics[height=0.85in]{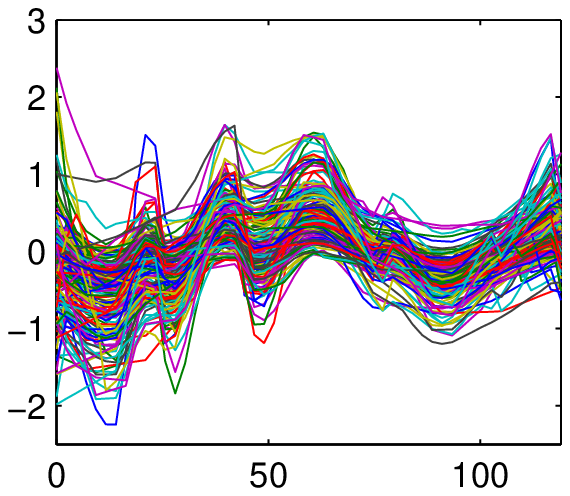}
\\
{\scriptsize Gene expression data}& {\scriptsize (1.81, 0.65, 0.96)}
&{\scriptsize (1.00, 1.14, 0.87)} & {\scriptsize ({\bf 0.77}, 1.53,
0.85)} & {\scriptsize (1.79, 0.72,
0.74)} & {\scriptsize (0.80, {\bf 2.36, 0.48})} \\
\hline
\end{tabular}
\caption{Empirical evaluation of five methods on 3 simulated
datasets and 4 real datasets, with the alignment performance
computed using three criteria $(ls, pc, sls)$. The best cases are
shown in boldface.} \label{fig:comparison}
\end{center}
\end{figure}

We compare our Fisher-Rao (F-R) method with the area under the curve
(AUTC) method presented in \cite{muller-JASA:2004}, the Tang-M\"{u}ller
method \cite{muller-biometrika:2008} provided in principal analysis by conditional expectation (PACE) package, the
self-modeling registration (SMR) method presented in
\cite{gervini-gasser-RSSB:04}, and the moment-based matching (MBM)
technique presented in \cite{james:10}. Fig. \ref{fig:comparison}
summarizes the values of $(ls, pc, sls)$ for these five methods
using 3 simulated and 4 real datasets. From the results, we can see
that the F-R method does uniformly well in functional alignment
under all the evaluation metrics. We have found that the $ls$
criterion is sometimes misleading in the sense that a low value can
result even if the functions are not very well aligned. This is the
case, for example, in the male growth data under SMR method. Here the
$ls = 0.45$, while for our method $ls = 0.64$, even though it is
easy to see that latter has performed a better alignment. On the
other hand, the $sls$ criterion seems to best correlate with a
visual evaluation of the alignment. Sometimes all three criteria
fail to evaluate the alignment performance properly, especially when
preserving the shapes of the original signals are considered. This
is the case in the first row of the figure where the AUTC method has
the same values of $ls$, $pc$, and $sls$ as our method but shapes of the
individual functions have been significantly distorted. The wave
function simulated data is the most challenging and no other method
except ours does a good job. Another point of evaluation is the
number of parameters used by different methods. While our method
does not have any parameter to choose, the other methods involve
choosing at least two but often more parameters which makes it
challenging for a user to apply them in different scenarios.
The computational costs associated with the different methods are given in 
Table \ref{tab:comp-cost}. This table is for some of the datasets used in our 
experiments and are representative of the general complexities of these methods. 
\begin{table} \label{tab:comp-cost}
\begin{center}
\begin{tabular}{|c|c|c|c|c|c|}
\hline Method & AUTC  \cite{muller-JASA:2004} & PACE
\cite{muller-biometrika:2008}
   & SMR \cite{gervini-gasser-RSSB:04} & MBM \cite{james:10} & F-R \\
\hline
Software & Matlab  & Matlab & Matlab & R & Matlab \\
\hline Gaussian kernel & 0.07sec & 68sec & 7.7sec & 101sec & 25sec \\
\hline Bimodal  & 0.02sec & 80sec & 4.5sec & 150sec & 17sec \\
\hline Growth-male & 0.03sec & 254sec & 14.5sec & 175sec & 22sec \\
\hline Signature  & 0.02sec & 145sec & 4.2sec & 117sec & 27sec \\
\hline
\end{tabular}
\caption{Computational cost of the five methods some datasets used 
in Fig. \ref{fig:comparison}. } \label{fig:cost}
\end{center}
\end{table}

\section{Discussion}
In this paper we have described a parameter-free approach for an automated alignment
of given functions using time warping. The basic idea is to use the Fisher-Rao
Riemannian metric and the resulting geodesic distance to define a proper distance,
called elastic distance, between
warping orbits of functions. This distance is used to compute a Karcher mean of the orbits,
and a template is selected from the mean orbit using an additional condition that the mean
of the warping functions is identity. Then, individual functions are aligned to the template
using the elastic distance and a natural separation of the amplitude and phase variability
of the function data is achieved.
One interesting application of this framework is in estimating a signal observed under random
time warpings. We propose the Karcher mean template as an estimator of the original signal
and prove that it is a consistent estimator of the signal under some basic conditions on the
random warping functions.

Important future directions in this work include: (1) the development of joint statistical models for
the amplitude and phase components of the data, and (2) the use of such models for classification
of observed functions into pre-determined classes. While the techniques for modeling the amplitude
variability are quite common, e.g. using functional principal component analysis, the corresponding
ideas for the phase component are relatively limited. The main reason is that $\Gamma$ is a nonlinear
manifold and one cannot directly apply FPCA here.
We mention that some solutions to this problem have been presented in
\cite{srivastava-etal-Fisher-Rao-CVPR:2007,srivastava-jermyn-PAMI:09,ashok-srivastava-etal-TIP:09}
albeit in different contexts.

\bibliography{bibfile}
\bibliographystyle{plain}

\appendix

\section{Proofs of Lemmas  \ref{lemma:transform} and \ref{lemma:isometry}}

\noindent
{\bf Proof of Lemma \ref{lemma:transform}}:
The mapping from $f$ to $q$ is as follows:
$f(t) \stackrel{{d \over dt}}{\rightarrow} \dot{f}(t) \stackrel{Q}{\rightarrow} q(t)$. For any
$v \in T_f({\cal F})$, the differential of this mapping is:
$v(t) \stackrel{{d \over dt}}{\rightarrow} \dot{v}(t) \stackrel{Q_{*,f(t)}}{\rightarrow} w(t)$.
To evaluate the expression for $w$, we need the expression for $Q_*$.
In case $x > 0$, we have
$Q(x) = \sqrt{x}$ and its directional derivative in the direction of
$y \in \real$ is $y/(2\sqrt{x})$. In case  $x < 0$, we have $Q(x) =
-\sqrt{-x}$ and its directional derivative in the direction of $y
\in \real$ is $y/(2\sqrt{-x})$. Combining the two, the directional
derivative of $Q$ is $Q_{*,x}(y) = y/(2\sqrt{|x|})$.
Now, to apply
this result to our situation, consider  two tangent vectors $v_1,
v_2 \in T_f({\cal F})$, and define their mappings under $Q_*$ as
$w_i(t) \equiv Q_{*,\dot{f}(t)}(\dot{v}_i(t)) =
\dot{v}_i(t)/(2\sqrt{|\dot{f}(t)|})$. Taking the $\ltwo$
inner-product between the resulting tangent vectors, we get:
$\inner{w_1(t)}{w_2(t)} = \int_0^1 w_1(t) w_2(t) dt = {1 \over 4} \int_0^1 \dot{v}_1(t) \dot{v}_2(t) {1 \over |\dot{f}(t)|} dt$.
The RHS is compared with Eqn. \ref{eq:def-FR} to complete the proof. $\Box$

\noindent
{\bf Proof of Lemma \ref{lemma:isometry}}:
For an arbitrary element $\gamma\in\Gamma$, and $q_1,\ q_2 \in \ltwo$,
we have: $\| (q_1,\gamma) - (q_2,\gamma) \|^2 =
 \int_0^1 (q_1(\gamma(t)) \sqrt{\dot{\gamma}(t)} -
q_2(\gamma(t)) \sqrt{\dot{\gamma}(t)})^2 dt  =
 \int_0^1 (q_1(\gamma(t))  -
q_2(\gamma(t)) )^2  \dot{\gamma}(t) dt = \| q_1 - q_2\|^2\ . \Box$

\section{Proofs of Lemma \ref{lem:sec}, Corollary \ref{cor:align}, and Lemma \ref{lem:mean-gamma}}

{\bf Proof of Lemma \ref{lem:sec}}: Using the definition:
$ \| cq_1 - (q_2, \gamma)\|^2 = \int_0^1 (cq_1(t) - (q_2, \gamma)(t))^2 dt
 = c^2 \| q_1\|^2 + \|q_2\|^2 - 2c\int_0^1 q_1(t)(q_2, \gamma)(t) dt$. Note
that we have used $\| (q_2, \gamma)\|^2 = \|q_2\|^2$, an important fact,  in the last
equality. Thus,
$$\argmin_{\gamma \in \Gamma} \| q_1 - (q_2, \gamma)\| =
\argmax_{\gamma \in \Gamma} \int_0^1 q_1(t)(q_2, \gamma)(t) dt
=
\argmin_{\gamma \in \Gamma} \| cq_1 - (q_2, \gamma)\|. \ \ \ \ \Box $$

\noindent
{\bf Proof of Corollary \ref{cor:align}}: $\gamma_{id} \in \argmin_{\gamma \in \Gamma} \| cq -
(q, \gamma)\|$ follows directly from Lemma \ref{lem:sec} since
 $\gamma_{id}$ minimizes $\| q- (q,\gamma)\|$.  Next we
show that this minimizer is unique if the set $\{t\in [0,1]|
q(t)=0\}$ has measure 0.  In this case, if we define $F(t) =
\int_0^t q(s)^2 ds$, then $F$ is a strictly increasing function on
$[0, 1]$.

Using Lemma \ref{lem:sec} again, we only need to show that
$\gamma_{id}$ is the unique minimizer for $\| q - (q,\gamma)\|$. For
any $\gamma^* \in \Gamma$ that minimizes $\| q- (q,\gamma)\|$, we
have $\|q - (q,\gamma^*)\| = \| q- (q,\gamma_{id})\| = 0$.
Therefore, $q(t) = q(\gamma^*(t)) \sqrt{\dot \gamma^*(t)}$ (almost
everywhere), and $F(t) = \int_0^t q(s)^2 ds = \int_0^t
q(\gamma^*(s))^2 \dot \gamma^*(s) ds = \int_0^{\gamma^*(t)} q(r)^2
dr = F(\gamma^*(t))$. As $F$ is strictly
increasing, we must have $\gamma^* = \gamma_{id}$. \ \ \ $\Box$\\

\noindent {\bf Proof of Lemma \ref{lem:mean-gamma}}: First we observe that for any two $\gamma_1, \gamma_2 \in \Gamma$, we have
$d_{FR}(\gamma_1, \gamma_2) = d_{FR}(\gamma_1 \circ \gamma, \gamma_2 \circ \gamma)$
for any $\gamma \in \Gamma$. This comes directly from the isometry of the group action of
$\Gamma$ on itself (proof is similar to that of Lemma \ref{lemma:isometry}).
This implies that:
$\argmin_{\gamma} \sum_{i=1}^n d_{FR}(\gamma_i \circ \gamma_0, \gamma)^2 =
\argmin_{\gamma} \sum_{i=1}^n d_{FR}(\gamma_i, \gamma \circ \gamma_0^{-1})^2$.
Let $\gamma^*$ denote the optimal $\gamma$ in the last term. Since
$\argmin_{\gamma} \sum_{i=1}^n d_{FR}(\gamma_i, \gamma )^2 = \bar
\gamma$, this implies that $\gamma^* \circ \gamma_0^{-1} = \bar
\gamma$ or $\gamma^* = \bar \gamma \circ \gamma_0$.
\ \ \ \ $\Box$

\end{document}